\documentclass[11pt,oneside]{article}

\usepackage{epsfig,lscape}
\usepackage{amssymb,amsfonts,amsmath}
\usepackage{rotating}

\usepackage{amsthm}
\def\me{\mathrm e}

\def\dif{\mathrm d}

\def\var{\mathrm{var}}

\def\N{\mbox{N}}

\def\T{ {\mathrm{\scriptscriptstyle T}} }

\def\one{ {\scriptscriptstyle 1} }

\def\zero{ {\scriptscriptstyle 0} }

\def\myeta{ {h} }

\newenvironment{prf}
{\noindent \textbf{Proof.}}{\hfill $\Box$ \vspace{.1in}}

\newtheorem{thm}{Theorem}
\newtheorem{lem}{Lemma}
\newtheorem{pro}{Proposition}

\newtheorem{ass}{Assumption}

\theoremstyle{definition}

\theoremstyle{definition}
\newtheorem{rem}{Remark}

\begin{document}

\begin{titlepage}

\begin{center}
{\Large Model-assisted inference for treatment effects using regularized calibrated estimation with high-dimensional data}

\vspace{.1in} Zhiqiang Tan\footnotemark[1]

\vspace{.1in}
\today
\end{center}

\footnotetext[1]{Department of Statistics \& Biostatistics, Rutgers University. Address: 110 Frelinghuysen Road,
Piscataway, NJ 08854. E-mail: ztan@stat.rutgers.edu. The research was supported in part by PCORI grant ME-1511-32740.
The author thanks Cun-Hui Zhang and Zijian Guo for helpful discussions.}

\paragraph{Abstract.}
Consider the problem of estimating average treatment effects when a large number of covariates are used to adjust for possible confounding
through outcome regression and propensity score models. The conventional approach of model building and fitting iteratively can be difficult to implement,
depending on ad hoc choices of what variables are included. In addition, uncertainty from the iterative process of model selection is complicated and often ignored
in subsequent inference about treatment effects. We develop new methods and theory
to obtain not only doubly robust point estimators for average treatment effects, which remain consistent if either the propensity score model
or the outcome regression model is correctly specified,
but also model-assisted confidence intervals, which are valid when
the propensity score model is correctly specified but the outcome regression model may be misspecified. With a linear outcome model, the
confidence intervals are doubly robust, that is, being also valid when the outcome model is correctly specified but the propensity score model may be misspecified.
Our methods involve regularized calibrated estimators with Lasso penalties, but carefully chosen loss functions, for fitting propensity score and outcome regression models.
We provide high-dimensional analysis to establish the desired properties of our methods under comparable conditions to previous results, which give valid confidence intervals when both
the propensity score and outcome regression are correctly specified.
We present a simulation study and an empirical application which confirm the advantages of the proposed methods compared with related methods based on
regularized maximum likelihood estimation.

%\vspace{-.2in}

\end{titlepage}

\section{Introduction} \label{sec:intro}

Drawing inferences about effects of treatments or interventions is constantly desired from observational studies in social and medical sciences, when
randomized experiments are either infeasible or difficult for practical constraints. This subject, broadly known as causal inference in statistics, is often
based on the framework of potential outcomes (Neyman 1923; Rubin 1974).
For observational studies, causal inference inevitably involves statistical modeling and estimation of population properties and
associations from empirical data (e.g., Tsiatis 2006). In particular, as the main problem to be tackled in the paper,
estimation of average treatment effects typically requires building and fitting outcome regression or propensity score models (e.g., Tan 2007).
The fitted outcome regression functions or propensity scores can then be used in various estimators for the average treatment effects,
notably inverse probability weighted (IPW) estimators or augmented IPW estimators (Robins et al.~1994).

For building and fitting outcome regression or propensity score models,
it is possible to follow the usual process of
model specification, fitting, and checking in a cyclic manner (e.g., McCullagh \& Nelder 1989).
In fact, a conventional approach for propensity score estimation as demonstrated in Rosenbaum \& Rubin (1984) involves
fitting a propensity score model (often logistic regression) by maximum likelihood, check covariate balance, and then modify and refit the propensity score model until reasonable balance is achieved.
However, this approach can be work intensive and difficult to implement, depending on ad hoc choices of what variables are included and whether nonlinear terms or
interactions are used among others. The situation can be especially challenging when there are a large number of
potentially confounding variables (or covariates) that need to be adjusted for in outcome regression or propensity score models.
%In addition, another statistical issue facing various IPW-based methods is that
%these methods can perform poorly, due to instability to small propensity scores estimated for few treated subjects,
%even when the propensity score model appears to be ``nearly correct" as judged by standard diagnostic techniques (e.g., Kang \& Schafer 2007).
In addition, another statistical issue is that uncertainty from the iterative process of model selection is complicated and often ignored in
subsequent inference (that is, confidence intervals or hypothesis testing) about treatment effects.

In this article, we develop new methods and theory for fitting logistic propensity score models and generalized linear outcome models and then using the fitted values in
augmented IPW estimators to estimates average treatment effects, in high-dimensional settings where the number of covariates $p$ is close to or even greater than
the sample size $n$. There are two main elements in our approach.
First, we employ regularized estimation with a Lasso penalty (Tibshirani 1992)
when fitting the outcome regression and propensity score models to deal with the large number of covariates under a sparsity assumption that
only a small but unknown subset (relative to the sample size) of covariates are associated with nonzero coefficients in the propensity score and outcome regression models.
Second, we carefully choose the loss functions for regularized estimation, different from least squares or maximum likelihood,
such that the resulting augmented IPW estimator and Wald-type confidence intervals possess the following properties (G1) and at least one of (G2)--(G3)
under suitable conditions:
\begin{itemize}\addtolength{\itemsep}{-.05in}
\item[(G1)] The point estimator is doubly robust, that is, remains consistent if either the propensity score model or the outcome regression model
is correctly specified.

\item[(G2)] The confidence intervals are valid if the propensity score model is correctly specified but the outcome regression model may be misspecified.

\item[(G3)] The confidence intervals are valid if the outcome regression model is correctly specified but the propensity score model may be misspecified.
\end{itemize}
If either property (G2) or (G3) is satisfied, then the confidence intervals are said to be model-assisted,
borrowing the terminology from the survey literature (Sarndal et al.~1992).
If properties (G2)--(G3) are satisfied, then the confidence intervals are doubly robust.

Combining the two foregoing elements leads to a regularized calibrated estimator, denoted by $\hat\gamma^1_{\mbox{\tiny RCAL}}$, for the coefficients in the propensity score model
and a regularized weighted likelihood estimator, denoted by $\hat\alpha^1_{\mbox{\tiny RWL}}$, for the coefficients in the outcome model within the treated subjects.
See the loss functions in (\ref{loss-CAL}) and (\ref{loss-WL}) or (\ref{loss-WL2}).
The regularized calibrated estimator  $\hat\gamma^1_{\mbox{\tiny RCAL}}$ has recently been proposed in Tan (2017) as an alternative to the regularized maximum likelihood estimator for fitting
logistic propensity score models, regardless of outcome regression models. As shown in Tan (2017), minimization of the underlying expected calibration loss implies reduction of
not only the expected
likelihood loss for logistic regression but also a measure of relative errors of limiting propensity scores that controls the mean squared errors of IPW estimators,
when the propensity score model may be misspecified.
In a complementary manner, our work here shows that $\hat\gamma^1_{\mbox{\tiny RCAL}}$ can be used in conjunction with $\hat\alpha^1_{\mbox{\tiny RWL}}$
to yield an augmented IPW estimator with valid confidence intervals if the propensity score model is correctly specified but
the outcome regression model may be specified.

We provide high-dimensional analysis of the regularized weighted likelihood estimator $\hat\alpha^1_{\mbox{\tiny RWL}}$
and the resulting augmented IPW estimator with possible model misspecification, while building on related results about $\hat\gamma^1_{\mbox{\tiny RCAL}}$ in Tan (2017).
In fact, a new strategy of inverting a quadratic inequality is developed to tackle the technical issue that
the weighted likelihood loss for $\hat\alpha^1_{\mbox{\tiny RWL}}$ is defined depending on the estimator $\hat\gamma^1_{\mbox{\tiny RCAL}}$.
As a result, we obtain the convergence of $\hat\alpha^1_{\mbox{\tiny RWL}}$ to a target value in the $L_1$ norm at the rate $(|S_\gamma| + |S_\alpha|) \{\log(p)/n\}^{1/2} $
and the symmetrized weighted Bregman divergence
at the rate $(|S_\gamma| + |S_\alpha|)  \log(p)/n$
under comparable conditions to those for high-dimensional analysis of standard Lasso estimators (e.g., Buhlmann \& van de Geer 2011),
where $|S_\gamma|$ denotes the size of nonzero coefficients of the propensity score model
and $ |S_\alpha|$ denotes that of the outcome model.
Furthermore, we establish an asymptotic expansion of the augmented IPW estimator based on $\hat\gamma^1_{\mbox{\tiny RCAL}}$ and $\hat\alpha^1_{\mbox{\tiny RWL}}$,
and show that property (G1) is achieved provided  $(|S_\gamma| + |S_\alpha|)  (\log p)^{1/2} = o(n^{1/2})$
and property (G2) is achieved provided $(|S_\gamma| + |S_\alpha|)  (\log p) = o(n^{1/2})$ with a nonlinear outcome model.
With a linear outcome model, we obtain stronger results:
property (G1) is achieved provided  $(|S_\gamma| + |S_\alpha|)  \log (p) = o(n)$
and both (G2) and (G3) are achieved provided $(|S_\gamma| + |S_\alpha|)  \log (p) = o(n^{1/2})$.
These sparsity conditions are as weak as in previous works (Belloni et al.~2014; van de Geer et al.~2014).

\vspace{.1in}
\noindent{\bf Related works.}
We compare and connect our work with related works in several areas. Non-penalized calibrated estimation for propensity score models have been studied, sometimes
independently (re)derived, in causal inference, missing-data problems, and survey sampling
(e.g., Folsom 1991; Tan 2010; Graham et al.~2012; Hainmueller 2012; Imai \& Ratovic 2014; Kim \& Haziza 2014; Vermeulen \& Vansteelandt 2015; Chan et al.~2016).
The non-penalized version of the estimator $\hat\alpha^1_{\mbox{\tiny RWL}}$ for outcome regression models have also been proposed in Kim \& Haziza (2014) and Vermeulen \& Vansteelandt (2015),
where one of the motivations is to circumvent the need of accounting for variation of such estimators of nuisance parameters and hence
simplify the computation of confidence intervals based on augmented IPW estimators.
Our work generalizes these ideas to achieve statistical advantages in high-dimensional settings,
where model-assisted or doubly robust confidence intervals would not be obtained without using
regularized calibration estimation. See Section~\ref{sec:discussion} for further discussion.

For high-dimensional causal inference, Belloni et al.~(2014) and Farrell (2015) employed augmented IPW estimators based on regularized
maximum likelihood estimators in outcome regression and propensity score models, and obtained Wald-type confidence intervals that are valid
when both the outcome regression and propensity score models are correctly specified, provided $(|S_\gamma| + |S_\alpha|)  \log (p) = o(n^{1/2})$.
Our main contribution is therefore to provide model-assisted or doubly robust confidence intervals using differently configured augmented IPW estimators for treatment effects.
As a secondary difference, Belloni et al.~(2014) and Farrell (2015) used post-Lasso estimators, that is, refitting outcome regression and propensity scores models
only including the variables selected from Lasso estimation. In contrast, our estimators $\hat\gamma^1_{\mbox{\tiny RCAL}}$ and $\hat\alpha^1_{\mbox{\tiny RWL}}$ are directly
Lasso penalized M-estimators.
%Our theoretical results are obtained by showing convergence in the $L_1$ norm and the symmetrized Bregman divergence, but without involving consistency of variable selection.

Another related work is Athey et al.~(2016), where  valid confidence intervals are obtained for the sample treatment effects such as
$n_1^{-1} \sum_{i: T_i=1} \{m_1^*(X_i) - m_0^*(X_i)\} $,
if a linear outcome model is correctly specified. No propensity score model is explicitly used.

Our work is also connected to the literature of confidence intervals and hypothesis testing for a single or lower-dimensional coefficients
in high-dimensional regression models (Zhang \& Zhang 2014; van de Geer et al.~2014; Javanmard \& Montanari 2014).
Model-assisted inference does not seem to be addressed in these works, but can potentially be developed.

\section{Setup} \label{sec:setup}

Suppose that %a simple random sample of $n$ subjects is obtained from a population under study.
the observed data consist
of independent and identically distributed observations $\{(Y_i, T_i, X_i): i=1,\ldots,n\}$ of $(Y,T,X)$, where
$Y$ is an outcome variable, $T$ is a treatment variable taking values 0 or 1, and $X$ is a vector of measured
covariates. In the potential outcomes framework for causal inference (Neyman 1923; Rubin 1974), let
$(Y^0, Y^1)$ be potential outcomes that would be observed under treatment 0 or 1 respectively. By
consistency, assume that $Y$ is either $Y^0$ if $T=0$ or $Y^1$ if $T=1$, that is,
$Y= (1-T) Y^0 + T Y^1$.
There are two causal parameters commonly of interest: the average treatment effect (ATE), defined as
$E(Y^1- Y^0) = \mu^1 - \mu^0$ with $\mu^t= E(Y^t)$, and the average treatment effect on the treated (ATT),
defined as $E(Y^1 - Y^0 |T=1) = \nu^1 - \nu^0$ with $\nu^t = E(Y^t | T=1)$ for $t=0,1$.
For concreteness, we mainly discuss estimation of $\mu^1$ until Section~\ref{sec:discussion} to discuss ATE and ATT.

Estimation of ATE is fundamentally a missing-data problem: only one potential outcome,
$Y^0_i$ or $Y^1_i$, is observed and the other one is missing for each subject $i$. For identification of
$(\mu^0,\mu^1)$ and ATE, we make the following two assumptions throughout:
\begin{itemize}\addtolength{\itemsep}{-.05in}
\item[(i)] Unconfoundedness: $T \perp Y^0 |X$ and $T \perp Y^1|X$, that is, $T$ and $Y^0$ and, respectively, $T$ and $Y^1$ are
conditionally independent given $X$ (Rubin 1976);

\item[(ii)] Overlap: $ 0 < \pi^*(x) < 1$ for all $x$, where $\pi^*(x)=P(T=1|X=x)$ is called the propensity score (PS) (Rosenbaum \& Rubin 1983).
\end{itemize}
Under these assumptions, $(\mu^0,\mu^1)$ and ATE are often estimated by imposing additional modeling (or dimension-reduction)
assumptions on the outcome regression function $m^*_t (X) = E(Y|T=t,X)$ or
the propensity score $\pi^*(X) = P(T=1|X)$.

Consider a conditional mean model for outcome regression (OR),
\begin{align}
E (Y | T=1, X) = m_1(X; \alpha^1 ) = \psi\{ \alpha^{\one\T} g^1(X) \}, \label{model-OR}
\end{align}
where $\psi()$ is an inverse link function, assumed to be increasing, $g^1(x) = \{1,g^1_1(x),\ldots, g^1_d(x)\}^\T$ is a vector of known functions,
and $\alpha^1  =(\alpha^1_0, \alpha^1_1,\ldots,\alpha^1_d)^\T$ is a vector of unknown parameters.
For example, model (\ref{model-OR}) can be deduced from a generalized linear model
with a canonical link (McCullagh \& Nelder 1989). Then the average negative log-(quasi-)likelihood function can be written (after dropping
any dispersion parameter) as
\begin{align}
\ell_{\mbox{\tiny ML}} (\alpha^1 ) = \tilde E \left( T \left[ -Y \alpha^{\one\T} g^1(X) + \Psi\{ \alpha^{\one\T} g^1(X)\} \right] \right), \label{loss-ML-OR}
\end{align}
where $\Psi (u) = \int_0^u \psi(u^\prime) \,\dif u^\prime$, which is convex in $u$. Throughout, $\tilde E()$ denotes the sample average.
With high-dimensional data, a regularized maximum likelihood estimator, $\hat\alpha^1_{\mbox{\tiny RML}}$, can be defined by minimizing the loss $\ell_{\mbox{\tiny ML}} (\alpha^1 )$ with
the Lasso penalty (Tibshirani 1992),
\begin{align}
\ell_{\mbox{\tiny RML}} (\alpha^1 ) &=\ell_{\mbox{\tiny ML}} (\alpha^1) + \lambda  \|\alpha^1_{1:d}\|_1, \label{reg-ml-loss-OR}
\end{align}
where $\|\alpha^1_{1:d}\|_1= \sum_{j=1}^d| \alpha^1_j |$ is the $L_1$ norm of $\alpha^1_{1:d}= ( \alpha^1_1, \ldots,\alpha^1_d)^\T$ excluding $\alpha^1_0$,
and $\lambda \ge 0$ is a tuning parameter.
The resulting estimator of $\mu^1$ is then
\begin{align*}
\hat \mu^1_{\mbox{\tiny OR}} = \tilde E \{ \hat m^1_{\mbox{\tiny RML}}(X)\} = \frac{1}{n} \sum_{i=1}^n \hat m^1_{\mbox{\tiny RML}}(X_i),
\end{align*}
where $\hat m^1_{\mbox{\tiny RML}}(X) = m_1(X; \hat \alpha^1_{\mbox{\tiny RML}} )$, the fitted outcome regression function.
Various theoretical results have been obtained on Lasso penalized estimation in sparse, high-dimensional regression (e.g., Buhlmann \& van de Geer 2011;
Huang \& Zhang~2012; Negahban et al.~2012).
Such results can be easily adapted to $\hat \alpha^1_{\mbox{\tiny RML}}$, with the data restricted to $\{(Y_i,X_i): T_i=1, i=1,\ldots,n\}$.
If model (\ref{model-OR}) is correctly specified, then it can be shown under suitable conditions that
$\|\hat \alpha^1_{\mbox{\tiny RML}}-\alpha^{\one*} \|_1 = O_p ( \|\alpha^{\one*}\|_0  \{\log(d)/n\}^{1/2} )$ and
$\hat \mu^1_{\mbox{\tiny OR}} = \mu^1 + O_p ( \{\|\alpha^{\one*}\|_0\log(d)/n\}^{1/2} )$,
where $\alpha^{\one*}$ is the true value for model (\ref{model-OR}) such that $m^*_1(X) = m_1(X; \alpha^{\one*})$.

Alternatively, consider a propensity score (PS) model
\begin{align}
P(T=1|X) = \pi(X;\gamma) = \Pi\{\gamma^\T f(X)\}, \label{model-PS}
\end{align}
where $\Pi()$ is an inverse link function, $f(x) = \{1,f_1(x),\ldots, f_p(x)\}^\T$ is a vector of known functions, and $\gamma=(\gamma_0,\gamma_1,\ldots,\gamma_p)^\T$ is a vector of unknown parameters.
For concreteness, assume that model (\ref{model-PS}) is logistic regression with $\pi(X;\gamma)= [1+\exp\{-\gamma^\T  f(X)\}]^{-1}$, and hence the average negative log-likelihood function is
\begin{align}
\ell_{\mbox{\tiny ML}}(\gamma) % & = \tilde E \left[ T_i \log\{ 1+\me^{-\gamma^\T f(X)} \} + (1-T_i) \log\{ 1+\me^{\gamma^\T f(X)}\}  \right] \nonumber \\
&= \tilde E \left[ \log\{ 1+\me^{ \gamma^\T f(X)} \}  -T\, \gamma^\T f(X) \right] .  \label{loss-ML-PS}
\end{align}
To handle high-dimensional data, a Lasso penalized maximum likelihood estimator, $\hat\gamma_{\mbox{\tiny RML}}$, is defined by minimizing the objective function
\begin{align}
\ell_{\mbox{\tiny RML}} (\gamma) &=\ell_{\mbox{\tiny ML}} (\gamma) + \lambda \|\gamma_{1:p}\|_1, \label{reg-ml-loss-PS}
\end{align}
where $\|\gamma_{1:p}\|_1= \sum_{j=1}^p | \gamma_j |$ is the $L_1$ norm of $\gamma_{1:p}= ( \gamma_1, \ldots,\gamma_p)^\T$ excluding $\gamma_0$,
and $\lambda \ge 0$ is a tuning parameter.
The fitted propensity score, $\hat\pi_{\mbox{\tiny RML}}(X) = \pi( X; \hat \gamma_{\mbox{\tiny RML}})$, can be used in various manners to estimate $(\mu^0,\mu^1)$ and ATE
including matching, stratification, and weighting. In particular, a (ratio) inverse probability weighted (IPW) estimator for $\mu^1$ is
\begin{align*}
\hat \mu^1_{\mbox{\tiny rIPW}} (\hat\pi_{\mbox{\tiny RML}}) = \tilde E\left\{ \frac{TY}{\hat\pi_{\mbox{\tiny RML}}(X)} \right\} \Big/\tilde E\left\{ \frac{T}{\hat\pi_{\mbox{\tiny RML}}(X)} \right\} .
\end{align*}
From previous works (Buhlmann \& van de Geer 2011; Huang \& Zhang~2012; Negahban et al.~2012),
if model (\ref{model-PS}) is correctly specified, then it can be shown under suitable conditions  that
$\| \hat\gamma_{\mbox{\tiny RML}} -\gamma^*\|_1 =O_p ( \|\gamma^*\|_0  \{\log(p)/n\}^{1/2} )$
and $\hat \mu^1_{\mbox{\tiny rIPW}} (\hat\pi_{\mbox{\tiny RML}}) = \mu^1 + O_p ( \{ \|\gamma^*\|_0 \log(p)/n\}^{1/2} )$,
where $\gamma^*$ is the true value for model (\ref{model-PS}) such that $\pi^*(X) = \pi(X; \gamma^*)$.

To attain consistency for $\mu^1$, the estimator $\hat \mu^1_{\mbox{\tiny OR}}$ or $\hat \mu^1_{\mbox{\tiny rIPW}} (\hat\pi_{\mbox{\tiny RML}})$
relies on correct specification of OR model (\ref{model-OR}) or PS model (\ref{model-PS}) respectively.
In contrast, there are doubly robust estimators depending on both OR and PS models in the augmented IPW form (Robins et al.~1994)
\begin{align*}
\hat \mu^1(\hat m_1, \hat\pi) = \tilde E \left\{ \varphi(Y,T,X;\hat m_1, \hat\pi) \right\} ,
%\tilde E \left[ \frac{TY}{\hat\pi (X)} - \left\{\frac{T}{\hat\pi (X)}-1 \right\} \hat m_1(X) \right],
\end{align*}
where $\hat m_1(X)$ and $\hat \pi(X)$ are fitted values of $m^*_1(X)$ and $\pi^*(X)$ respectively and
\begin{align}
\varphi(Y,T,X;\hat m_1, \hat\pi) = \frac{TY}{\hat\pi (X)} - \left\{\frac{T}{\hat\pi (X)}-1 \right\} \hat m_1(X) . \label{influence-function}
\end{align}
%For example, the augmented IPW estimator in Robins et al.~(1994) is obtained as
%$\hat \mu^1(\hat m^1_{\mbox{\tiny ML}}, \hat\pi_{\mbox{\tiny ML}})$, where $\hat m^1_{\mbox{\tiny ML}}(X)$ and $\hat\pi_{\mbox{\tiny ML}}(X)$ are
%the fitted values by maximum likelihood.
See Kang \& Schafer (2007) and Tan (2010) for reviews in low-dimensional settings.
Recently, interesting results  in high-dimensional settings have been obtained by Belloni et al.~(2014) and Farrell (2015)
on the estimator $\hat \mu^1(\hat m^1_{\mbox{\tiny RML}}, \hat\pi_{\mbox{\tiny RML}})$, using the fitted values
$\hat m^1_{\mbox{\tiny RML}}(X)$ and $\hat\pi_{\mbox{\tiny RML}}(X)$ from Lasso penalized estimation or similar methods.
These results are mainly of two types.
The first type shows double robustness: $\hat \mu^1(\hat m^1_{\mbox{\tiny RML}}, \hat\pi_{\mbox{\tiny RML}})$
remains consistent if either OR model (\ref{model-OR}) or PS model (\ref{model-PS}) is correctly specified.
The second type establishes valid confidence intervals: $\hat \mu^1(\hat m^1_{\mbox{\tiny RML}}, \hat\pi_{\mbox{\tiny RML}})$ admits the usual influence function,
\begin{align*}
\hat \mu^1(\hat m^1_{\mbox{\tiny RML}}, \hat\pi_{\mbox{\tiny RML}})=
\tilde E\left\{ \varphi(Y,T,X; m_1^*, \pi^*)  \right\} +o_p (n^{-1/2})
%\tilde E\left[ \frac{TY}{\pi^* (X)} - \left\{\frac{T}{\pi^* (X)}-1 \right\}  m_1^*(X) \right] +o_p (n^{-1/2})
\end{align*}
if both OR model (\ref{model-OR}) and  PS model (\ref{model-PS}) are correctly specified.
In general, the latter result requires a stronger sparsity condition than in consistency results only. For example, it is assumed
that $\{ \|\alpha^{\one*}\|_0 + \| \gamma^*\|_0\} \log(p)= o(n^{1/2})$ in Belloni et al.~(2014).

\section{Theory and methods}

\subsection{Overview} \label{sec:overview}

An important limitation of existing works discussed in Section~\ref{sec:setup} is that valid confidence intervals based on
$\hat \mu^1(\hat m^1_{\mbox{\tiny RML}}, \hat\pi_{\mbox{\tiny RML}})$ is obtained only under the assumption that
both OR model (\ref{model-OR}) and  PS model (\ref{model-PS}) are correctly specified, even though the point estimator
$\hat \mu^1(\hat m^1_{\mbox{\tiny RML}}, \hat\pi_{\mbox{\tiny RML}})$ is
doubly robust, that is, remains consistent if either OR model (\ref{model-OR}) or PS model (\ref{model-PS}) is correctly specified.
To fill this gap, we develop new point estimators and confidence intervals for $\mu^1$, depending on a propensity score model and an outcome regression model, such that
properties (G1) and at least one of (G2)--(G3) are attained as described in Section~\ref{sec:intro}.
Obtaining model-assisted or doubly robust confidence intervals presents a considerable improvement over existing theory and methods in Belloni et al.~(2014) and Farrell (2015).

To illustrate main ideas, consider a logistic propensity score model (\ref{model-PS}) and a linear outcome regression model,
\begin{align}
E (Y | T=1, X) = m_1(X; \alpha^1 ) = \alpha^{\one\T} f(X) , \label{lm-OR}
\end{align}
that is, model (\ref{model-OR}) with the identity link and the vector of covariate functions $g^1(X)$ taken to be the same as $f(X)$ in model (\ref{model-PS}).
This condition can be satisfied possibly after enlarging model (\ref{model-OR}) or (\ref{model-PS}) to reach the same dimension.
Our point estimator of $\mu^1$ is
\begin{align}
\hat \mu^1(\hat m^1_{\mbox{\tiny RWL}}, \hat\pi^1_{\mbox{\tiny RCAL}}) = \tilde E \left\{ \varphi(Y,T,X;\hat m^1_{\mbox{\tiny RWL}}, \hat\pi^1_{\mbox{\tiny RCAL}}) \right\} , \label{mu-estimator}
\end{align}
where $\varphi()$ is defined in (\ref{influence-function}), $\hat\pi^1_{\mbox{\tiny RCAL}}(X) = \pi( X; \hat \gamma^1_{\mbox{\tiny RCAL}})$,
$\hat m^1_{\mbox{\tiny RWL}}(X) = m_1(X; \hat \alpha^1_{\mbox{\tiny RWL}} )$,
and $\hat \gamma^1_{\mbox{\tiny RCAL}}$ and $\hat \alpha^1_{\mbox{\tiny RWL}}$ are defined as follows.
The estimator $\hat\gamma^1_{\mbox{\tiny RCAL}}$ is a regularized calibrated estimator of $\gamma$ from Tan (2017),
defined as a minimizer of
the Lasso penalized objective function,
\begin{align}
\ell_{\mbox{\tiny RCAL}} (\gamma) &=\ell_{\mbox{\tiny CAL}} (\gamma) + \lambda \|\gamma_{1:p}\|_1, \label{reg-cal-loss}
\end{align}
where $\ell_{\mbox{\tiny RCAL}} (\gamma)$ is the calibration loss,
\begin{align}
\ell_{\mbox{\tiny CAL}} (\gamma) &= \tilde E \left\{ T \me^{-\gamma^\T f(X)} + (1-T) \gamma^\T f(X) \right\}, \label{loss-CAL}
%&= \tilde E \left[ T \frac{1-\pi(X;\gamma)}{\pi(X;\gamma)} + (1-T) \log \left\{\frac{\pi(X;\gamma)}{1-\pi(X;\gamma)} \right\} \right]
\end{align}
and $\|\gamma_{1:p}\|_1$ is the $L_1$ norm of $\gamma_{1:p}$ and $\lambda \ge 0$ is a tuning parameter.
The estimator $\hat \alpha^1_{\mbox{\tiny RWL}} $ is a regularized weighted least-squares estimator of $\alpha^1$, defined as a minimizer of
\begin{align}
\ell_{\mbox{\tiny RWL}} (\alpha^1; \hat\gamma^1_{\mbox{\tiny RCAL}} ) &=\ell_{\mbox{\tiny WL}} (\alpha^1; \hat\gamma^1_{\mbox{\tiny RCAL}}) + \lambda  \|\alpha^1_{1:p}\|_1, \label{reg-wl-loss}
\end{align}
where $\ell_{\mbox{\tiny WL}} (\alpha^1 ; \hat\gamma^1_{\mbox{\tiny RCAL}}) $ is the weighted least-squares loss,
\begin{align}
\ell_{\mbox{\tiny WL}} (\alpha^1; \hat\gamma^1_{\mbox{\tiny RCAL}} ) = \tilde E \left[ T \frac{1-\hat\pi^1_{\mbox{\tiny RCAL}}(X)}{\hat\pi^1_{\mbox{\tiny RCAL}}(X)} \left\{Y - \alpha^{\one\T} f(X) \right\}^2 \right]/2, \label{loss-WL}
\end{align}
and $\|\alpha^1_{1:p}\|_1$ is the $L_1$ norm of $\alpha^1_{1:p}$ and $\lambda \ge 0$ is a tuning parameter.
That is, the observations in the treated group are weighted by $\{1-\hat\pi^1_{\mbox{\tiny RCAL}}(X_i)\}/\hat\pi^1_{\mbox{\tiny RCAL}}(X_i)$,
which differs slightly from the commonly used inverse propensity score weight $1/\hat\pi^1_{\mbox{\tiny RCAL}}(X_i)$.
%High (or low) weights are associated with observations with low (or high) fitted propensity scores, thereby focusing the fitting of model (\ref{lm-OR}) in
%covariate regions with most missing data.

There are simple and interesting interpretations of the preceding estimators.
By the Karush--Kuhn--Tucker condition for minimizing (\ref{reg-cal-loss}), the fitted propensity score $\hat \pi^1_{\mbox{\tiny RCAL}}(X)$ satisfies
\begin{align}
& \frac{1}{n} \sum_{i=1}^n \frac{T_i}{\hat \pi^1_{\mbox{\tiny RCAL}}(X_i) } = 1, \label{ineq-CAL-1} \\
& \frac{1}{n} \left| \sum_{i=1}^n \frac{T_i f_j(X_i)}{\hat \pi^1_{\mbox{\tiny RCAL}}(X_i) } - \sum_{i=1}^n f_j(X_i) \right| \le \lambda, \quad j=1,\ldots,p , \label{ineq-CAL-2}
\end{align}
where equality holds in (\ref{ineq-CAL-2}) for any $j$ such that the $j$th estimate $(\hat\gamma^1_{\mbox{\tiny RCAL}})_j$ is nonzero.
Eq.~(\ref{ineq-CAL-1}) shows that
the inverse probability weights, $1/\hat \pi^1_{\mbox{\tiny RCAL}}(X_i)$ with $T_i=1$, sum to the sample size $n$ by (\ref{ineq-CAL-1}), whereas
Eq.~(\ref{ineq-CAL-2}) implies that
the weighted average of each covariate $f_j(X_i)$ over the treated group may differ from the overall average of $f_j(X_i)$ by no more than $\lambda$.
In fact, the calibration loss $\ell_{\mbox{\tiny CAL}} (\gamma)$ in (\ref{loss-CAL}) is derived such that its gradient gives the left hand side of (\ref{ineq-CAL-2})
without taking absolute values, as shown in Eq.~(\ref{grad-gamma}). The Lasso penalty is used to induce the box constraints on the gradient of $\ell_{\mbox{\tiny CAL}}(\gamma)$
instead of setting the gradient to 0.

By the Karush--Kuhn--Tucker condition for minimizing (\ref{reg-wl-loss}), the fitted outcome regression function $\hat m^1_{\mbox{\tiny RWL}}(X)$ satisfies
\begin{align}
& \frac{1}{n} \sum_{i=1}^n T_i \frac{1-\hat \pi^1_{\mbox{\tiny RCAL}}(X_i) }{\hat \pi^1_{\mbox{\tiny RCAL}}(X_i) }  \left\{Y_i - \hat m^1_{\mbox{\tiny RWL}}(X_i) \right\}= 0,  \label{ineq-WL-1} \\
& \frac{1}{n}\left| \sum_{i=1}^n T_i \frac{1-\hat \pi^1_{\mbox{\tiny RCAL}}(X_i) }{\hat \pi^1_{\mbox{\tiny RCAL}}(X_i) }  \left\{Y_i - \hat m^1_{\mbox{\tiny RWL}}(X_i) \right\}f_j(X_i) \right| \le \lambda,
\quad j=1,\ldots,p, \label{ineq-WL-2}
\end{align}
where equality holds in (\ref{ineq-WL-2}) for any $j$ such that the $j$th estimate $(\hat\alpha^1_{\mbox{\tiny RWL}})_j$ is nonzero.
Eq.~(\ref{ineq-WL-1}) implies that by simple calculation,
the estimator $\hat \mu^1(\hat m^1_{\mbox{\tiny RWL}}, \hat\pi^1_{\mbox{\tiny RCAL}})$ can be recast as
\begin{align}
\hat \mu^1(\hat m^1_{\mbox{\tiny RWL}}, \hat\pi^1_{\mbox{\tiny RCAL}}) & =
\tilde E \left[\hat m^1_{\mbox{\tiny RWL}}(X) + \frac{T}{\hat \pi^1_{\mbox{\tiny RCAL}}(X)} \left\{ Y - \hat m^1_{\mbox{\tiny RWL}}(X) \right\} \right] \nonumber \\
& = \tilde E \left\{ T Y + (1-T) \hat m^1_{\mbox{\tiny RWL}}(X) \right\}, \label{mu-prediction-form}
\end{align}
which takes the form of linear prediction estimators known in the survey literature (e.g., Sarndal et al.~1992):  %Sarndal \& Wright 1984):
$\tilde E \{ T Y + (1-T) \hat m_1(X)\}$ for some fitted outcome regression function $\hat m_1(X)$.
As a consequence, $\hat \mu^1(\hat m^1_{\mbox{\tiny RWL}}, \hat\pi^1_{\mbox{\tiny RCAL}})$ always falls within the range of
the observed outcomes $\{Y_i: T_i=1, i=1,\ldots,n\}$ and the predicted values $\{ \hat m^1_{\mbox{\tiny RWL}}(X_i): T_i=0, i=1,\ldots,n\}$.
This boundedness property is not satisfied by the estimator $\hat \mu^1(\hat m^1_{\mbox{\tiny RML}}, \hat\pi^1_{\mbox{\tiny RML}})$.

We provide a high-dimensional analysis of the estimator $\hat \mu^1(\hat m^1_{\mbox{\tiny RWL}}, \hat\pi^1_{\mbox{\tiny RCAL}})$  in Section~\ref{sec:linear-OR}, allowing possible model
 misspecification. Our main result shows that under suitable conditions, the estimator $\hat \mu^1(\hat m^1_{\mbox{\tiny RWL}}, \hat\pi^1_{\mbox{\tiny RCAL}})$
admits the asymptotic expansion
\begin{align}
\hat \mu^1(\hat m^1_{\mbox{\tiny RWL}}, \hat\pi^1_{\mbox{\tiny RCAL}}) =
\tilde E \left\{ \varphi(Y,T,X;\bar m^1_{\mbox{\tiny WL}}, \bar \pi^1_{\mbox{\tiny CAL}}) \right\} +o_p (n^{-1/2}), \label{new-expansion}
% = \tilde E\left[ \frac{TY}{\bar \pi^1_{\mbox{\tiny CAL}} (X)} - \left\{\frac{T}{\bar \pi^1_{\mbox{\tiny CAL}} (X)}-1 \right\}  \bar m^1_{\mbox{\tiny WL}}(X) \right]
\end{align}
where  $\bar\pi^1_{\mbox{\tiny CAL}}(X) = \pi( X; \bar \gamma^1_{\mbox{\tiny CAL}})$,
$\bar m^1_{\mbox{\tiny WL}}(X) = m_1(X; \bar \alpha^1_{\mbox{\tiny WL}} )$
and $\bar \gamma^1_{\mbox{\tiny CAL}}$ and $\bar \alpha^1_{\mbox{\tiny WL}}$ are defined as follows.
In the presence of possible model misspecification, the target value $\bar \gamma^1_{\mbox{\tiny CAL}}$ is defined as a minimizer of the expected calibration loss
\begin{align*}
E \left\{ \ell_{\mbox{\tiny CAL}} (\gamma) \right\} &= E \left\{ T \me^{-\gamma^\T f(X)} + (1-T) \gamma^\T f(X) \right\}.
\end{align*}
If model (\ref{model-PS}) is correctly specified, then $\bar \pi^1_{\mbox{\tiny CAL}} (X) =\pi^*(X)$.
Otherwise, $\bar \pi^1_{\mbox{\tiny CAL}} (X)$ may differ from $\pi^*(X)$.
The target value $\bar \alpha^1_{\mbox{\tiny WL}}$ is defined as a minimizer of the expected loss
\begin{align*}
E \left\{ \ell_{\mbox{\tiny WL}} (\alpha^1; \bar\gamma^1_{\mbox{\tiny CAL}} ) \right\}
= E \left[ T \frac{1-\bar\pi^1_{\mbox{\tiny CAL}}(X)}{\bar\pi^1_{\mbox{\tiny CAL}}(X)} \left\{Y - \alpha^{\one\T} f(X) \right\}^2 \right] /2.
\end{align*}
If model (\ref{lm-OR}) is correctly specified, then $\bar m^1_{\mbox{\tiny WL}} (X) =m_1^*(X)$.
But $\bar m^1_{\mbox{\tiny WL}} (X)$ may in general differ from $m^*_1(X)$.
For concreteness, the following result can be deduced from Theorems~\ref{thm-mu} and \ref{thm-var}.
Suppose that the Lasso tuning parameter is specified as $\lambda=A_0^\dag \{\log(p)/n\}^{1/2}$ for $\hat\gamma^1_{\mbox{\tiny RCAL}}$ and
$\lambda = A_1^\dag \{\log(p)/n\}^{1/2}$ for $\hat\alpha^1_{\mbox{\tiny RWL}}$, with some constants $A_0^\dag$ and $A_1^\dag$.
Denote $S_\gamma = \{0\} \cup \{j: \bar\gamma^1_{\mbox{\tiny CAL},j} \not=0, j=1,\ldots,p\}$ and
$S_\alpha = \{0\} \cup \{j: \bar\alpha^1_{\mbox{\tiny WL},j} \not=0, j=1,\ldots,p\}$.

\begin{pro} \label{pro-mu-var}
Suppose that Assumptions~\ref{ass-RCAL} and \ref{ass-RWL} hold as in Section~\ref{sec:linear-OR}, and $( |S_\gamma| + |S_\alpha| ) \log(p) = o(n^{1/2})$.
For sufficiently large constants $A_0^\dag$ and $A_1^\dag$, if either logistic PS model (\ref{model-PS}) or linear OR model (\ref{lm-OR}) is correctly specified, then
the following results hold:
\begin{itemize}\addtolength{\itemsep}{-.05in}
\item[(i)]
$n^{1/2} \{ \hat \mu^1(\hat m^1_{\mbox{\tiny RWL}}, \hat\pi^1_{\mbox{\tiny RCAL}}) - \mu^1\} \to_{\mathcal D} \N(0,V)$,
where $V = \var\{ \varphi(Y,T,X;\bar m^1_{\mbox{\tiny WL}}, \bar \pi^1_{\mbox{\tiny CAL}}) \}$;

\item[(ii)] a consistent estimator of $V$ is
\begin{align*}
\hat V = \tilde E \left[ \left\{\varphi(Y,T,X;\hat m^1_{\mbox{\tiny RWL}}, \hat\pi^1_{\mbox{\tiny RCAL}}) - \hat \mu^1(\hat m^1_{\mbox{\tiny RWL}}, \hat\pi^1_{\mbox{\tiny RCAL}})\right\}^2 \right];
\end{align*}

\item[(iii)] an asymptotic $(1-c)$ confidence interval for $\mu^1$ is $ \hat \mu^1(\hat m^1_{\mbox{\tiny RWL}}, \hat\pi^1_{\mbox{\tiny RCAL}}) \pm z_{c/2} \sqrt{\hat V/n}$,
where $z_{c/2}$ is the $(1-c/2)$ quantile of $\N(0,1)$.
\end{itemize}
That is, a doubly robust confidence interval for $\mu^1$ is obtained.
\end{pro}

We highlight some basic ideas underlying the construction of the estimators $\hat\gamma^1_{\mbox{\tiny RCAL}}$ and $\hat\alpha^1_{\mbox{\tiny RWL}}$ as well as
the proof of the asymptotic expansion (\ref{new-expansion}) for $\hat \mu^1(\hat m^1_{\mbox{\tiny RWL}}, \hat\pi^1_{\mbox{\tiny RCAL}})$.
For an estimator $\hat\gamma$ in model (\ref{model-PS}), suppose that $\hat\gamma$ converges in probability to a limit $\bar\gamma$.
Denote $\hat\pi(X) = \pi(X;\hat\gamma)$ and $\bar\pi(X) = \pi(X;\bar\gamma)$.
Similarly, for an estimator $\hat\alpha^1$ in model (\ref{model-OR}), suppose that $\hat\alpha^1$ converges in probability to a limit $\bar\alpha^1$.
Denote $\hat m_1(X) = \hat\alpha^{\one\T}f(X)$ and $\bar m_1(X) = \bar\alpha^{\one\T}f(X)$.
Consider the following decomposition of $\hat \mu^1(\hat m_1, \hat \pi)$ by direct calculation:
\begin{align}
\hat \mu^1(\hat m_1, \hat \pi) & = \hat \mu^1(\bar m_1, \bar \pi) + \tilde E \left[ \{\hat m_1(X) - \bar m_1(X)\} \left\{ 1 - \frac{T}{\hat \pi(X)} \right\}  \right]  \nonumber \\
& \quad + \tilde E \left[ T\{Y-\bar m_1(X)\} \left\{ \frac{1}{\hat \pi(X)} - \frac{1}{\bar \pi(X)} \right\} \right] . \label{mu-decomp}
% & \quad + \tilde E \left[ \{\hat m_1(X) - \bar m_1(X)\} \left\{ \frac{T}{\bar \pi(X)} - \frac{T}{\hat \pi(X)} \right\} \right].
\end{align}
Eq.~(\ref{mu-decomp}) can also be obtained from a Taylor expansion with $(\hat\alpha^1,\hat\gamma)$ about $(\bar\alpha^1,\bar\gamma)$.
For linear OR model (\ref{lm-OR}), the second term of the decomposition reduces to
\begin{align}
(\hat \alpha^1 - \bar\alpha^1)^\T \times \tilde E \left[ \left\{ 1 - \frac{T}{\hat \pi(X)} \right\}  f(X) \right] . \label{decomp-2}
\end{align}
For logistic PS model (\ref{model-PS}) with $\partial \pi(X;\gamma)/ \partial\gamma = \pi(X;\gamma) \{1-\pi(X;\gamma)\}$, the third term of the decomposition can be approximated via a Taylor expansion by
\begin{align}
-(\hat\gamma-\bar\gamma)^\T \times \tilde E \left[ T \frac{1-\bar\pi(X)}{\bar\pi(X)} \{Y-\bar m_1(X)\} f(X) \right].   \label{decomp-3}
\end{align}
Suppose that $\hat\gamma$ and $\hat\alpha^1$ are Lasso penalized M-estimators such that under suitable conditions,
$ \|\hat\gamma - \bar\gamma\|_1 = O_p ( \{ \log(p)/n \}^{1/2} )$
and $\| \hat \alpha^1 - \bar \alpha^1 \|_1 = O_p(\{ \log(p)/n \}^{1/2} )$, where for simplicity the dependency on the sparsity sizes of $\bar \gamma$ and $\bar \alpha^1$ are suppressed.
The loss functions $\ell_{\mbox{\tiny CAL}}(\gamma) $ and $\ell_{\mbox{\tiny WL}}(\alpha^1; \gamma)$ in (\ref{loss-CAL}) and (\ref{loss-WL}) are constructed such that
\begin{align}
& \frac{\partial \ell_{\mbox{\tiny CAL}}(\gamma)}{\partial \gamma} = \tilde E \left[ \left\{ 1 - \frac{T}{\pi(X;\gamma)} \right\}  f(X) \right] , \label{grad-gamma} \\
& \frac{\partial \ell_{\mbox{\tiny WL}}(\alpha^1; \gamma)}{\partial \alpha^1} =  -\tilde E \left[ T \frac{1-\pi(X;\gamma)}{\pi(X;\gamma)} \{Y- \alpha^{\one\T}f(X)\} f(X) \right]. \label{grad-alpha}
\end{align}
Then the second terms in (\ref{decomp-2}) and (\ref{decomp-3}) can be of order $O_p ( \{ \log(p)/n \}^{1/2} )$ in the supremum norms,
as reflected in conditions (\ref{ineq-CAL-1})--(\ref{ineq-CAL-2})
and (\ref{ineq-WL-1})--(\ref{ineq-WL-2}).
Consequently, the products (\ref{decomp-2}) and (\ref{decomp-3}) can be of order $O_p ( \log(p)/n )$,
which becomes $o_p(n^{-1/2})$ and hence (\ref{new-expansion}) holds provided $\log(p) = o(n^{1/2})$ up to a constant depending on the sparsity sizes  of $\bar \gamma$ and $\bar \alpha^1$.

The estimator $\hat\gamma^1_{\mbox{\tiny RCAL}}$ is called a regularized calibrated estimator of $\gamma$ (Tan 2017), because in the extreme case of $\lambda=0$, Eqs.~(\ref{ineq-CAL-1})--(\ref{ineq-CAL-2})
reduce to calibration equations, which can be traced to Folsom (1991) in the survey literature. Although such equations are intuitively appealing, the preceding discussion
shows that $\hat\gamma^1_{\mbox{\tiny RCAL}}$ can also be derived to reduce the variation associated with estimation of $\alpha^1$ from linear OR model (\ref{lm-OR}) for the estimator $\hat \mu^1(\hat m_1, \hat \pi)$, when
PS model (\ref{model-PS}) may be misspecified. Similarly,
 $\hat\alpha^1_{\mbox{\tiny RWL}}$ is constructed to reduce the variation associated with estimation of $\gamma$ from logistic PS model (\ref{model-PS}) for the estimator $\hat \mu^1(\hat m_1, \hat \pi)$, when
OR model (\ref{lm-OR}) may be misspecified.
By extending the meaning of calibrated estimation, we call
$\hat\alpha^1_{\mbox{\tiny RWL}}$ a regularized calibrated estimator of $\alpha^1$ against model (\ref{model-PS}), as well as
$\hat\gamma^1_{\mbox{\tiny RCAL}}$ a regularized calibrated estimator of $\gamma$ against model (\ref{lm-OR}), when used to define $\hat \mu^1(\hat m_1, \hat \pi)$.

While the preceding discussion outlines our basic reasoning, there are several technical issues we need to address in high-dimensional analysis, including
how to handle the dependency of the estimator $\hat\alpha^1_{\mbox{\tiny RWL}}$ on $\hat\gamma^1_{\mbox{\tiny RCAL}}$,
and what condition is required on the sparsity sizes  of $\bar \gamma$ and $\bar \alpha^1$.
In addition, we develop appropriate methods and theory in the situation where a generalized linear model (\ref{model-OR}), not just linear model (\ref{lm-OR}), is used for outcome regression.

\subsection{Using linear outcome regression} \label{sec:linear-OR}

In this section, we assume that linear outcome model (\ref{lm-OR}) is used together with logistic propensity score model (\ref{model-PS}),
and develop theoretical results for the proposed estimator $\hat \mu^1(\hat m^1_{\mbox{\tiny RWL}}, \hat\pi^1_{\mbox{\tiny RCAL}})$, leading to
Proposition~\ref{pro-mu-var} among others, in high-dimensional settings.

First we describe relevant results from Tan (2017) about the behavior of the regularized calibrated estimator $\hat\gamma^1_{\mbox{\tiny RCAL}}$ in model (\ref{model-PS}).
The tuning parameter $\lambda$ used in (\ref{reg-cal-loss}) for defining $\hat\gamma^1_{\mbox{\tiny RCAL}}$ is specified as
$\lambda = A_0 \lambda_0$, with a constant $A_0>1$ and
\begin{align*}
\lambda_0 = C_1\sqrt{ \log\{(1+p)/\epsilon\} / n},
\end{align*}
where $C_1 >0$ is a constant depending only on $(C_0,B_0)$ from Assumption~\ref{ass-RCAL} below and $0 < \epsilon <1$ is a tail probability for the error bound.
For example, taking $\epsilon=1/(1+p)$ gives $\lambda_0 =C_1\sqrt{ 2 \log(1+p)/n}$, a familiar rate in high-dimensional analysis.

With possible model misspecification, the target value $\bar\gamma^1_{\mbox{\tiny CAL}}$ is defined
as a minimizer of the expected calibration loss $E\{\ell_{\mbox{\tiny CAL}}(\gamma)\}$ as in Section~\ref{sec:overview}.
From a functional perspective, we write $\ell_{\mbox{\tiny CAL}}(\gamma) = \kappa_{\mbox{\tiny CAL}}(\gamma^\T f)$, where for a function $\myeta(x)$,
\begin{align*}
\kappa_{\mbox{\tiny CAL}} (\myeta) = \tilde E \left[ T \me^{-\myeta(X)} + (1-T) \myeta(X) \right].
\end{align*}
As $\kappa_{\mbox{\tiny CAL}} (\myeta)$ is easily shown to be convex in $\myeta$, the Bregman divergence associated with $\kappa_{\mbox{\tiny CAL}}$
is defined such that for two functions $\myeta(x)$ and $\myeta^\prime(x)$,
\begin{align*}
D_{\mbox{\tiny CAL}} (\myeta^\prime, \myeta) = \kappa_{\mbox{\tiny CAL}} (\myeta^\prime) - \kappa_{\mbox{\tiny CAL}} (\myeta) - \langle \nabla\kappa_{\mbox{\tiny CAL}} (\myeta), \myeta^\prime - \myeta\rangle,
\end{align*}
where $\myeta$ is identified as a vector $(\myeta_1,\ldots,\myeta_n)$ with $\myeta_i = \myeta(X_i)$,
and $\nabla\kappa_{\mbox{\tiny CAL}} (\myeta)$ denotes the gradient of $\kappa_{\mbox{\tiny CAL}} (\myeta)$ with respect to $(\myeta_1,\ldots,\myeta_n)$.
The following result (Theorem~\ref{thm-RCAL}) is restated from Tan (2017, Corollary 2), where the convergence of $\hat\gamma^1_{\mbox{\tiny RCAL}}$ to $\bar\gamma^1_{\mbox{\tiny CAL}}$ is
obtained in the $L_1$ norm
$\|\hat\gamma^1_{\mbox{\tiny RCAL}} -\bar\gamma^1_{\mbox{\tiny CAL}} \|_1$ and
the symmetrized Bregman divergence
\begin{align*}
D^\dag_{\mbox{\tiny CAL}} (\hat\myeta^1_{\mbox{\tiny RCAL}}, \bar\myeta^1_{\mbox{\tiny CAL}} ) =
D_{\mbox{\tiny CAL}} (\hat\myeta^1_{\mbox{\tiny RCAL}}, \bar\myeta^1_{\mbox{\tiny CAL}} )  +
D_{\mbox{\tiny CAL}} (\bar\myeta^1_{\mbox{\tiny CAL}}, \hat\myeta^1_{\mbox{\tiny RCAL}} ) ,
\end{align*}
where $\hat\myeta^1_{\mbox{\tiny RCAL}}(X) = \hat\gamma^{\one\T}_{\mbox{\tiny RCAL}} f(X)$ and $\bar\myeta^1_{\mbox{\tiny CAL}}(X) = \bar\gamma^{\one\T}_{\mbox{\tiny CAL}} f(X)$.
See Lemma~\ref{lem-add-ineq} in the Supplementary Material for an explicit expression of $D^\dag_{\mbox{\tiny CAL}}$.

For a matrix $\Sigma$ with row indices $\{0,1,\ldots,k\}$,
a compatibility condition (Buhlmann \& van de Geer 2011) is said to hold with a subset $S \in \{0,1,\ldots,k\}$ and constants $\nu_0 >0$ and $\xi_0>1$ if
$\nu_0^2  (\sum_{j\in S} |b_j|)^2 \le |S| ( b^\T \Sigma b )$
for any vector $b=(b_0,b_1,\ldots,b_k)^\T \in \mathbb R^{1+k} $ satisfying
\begin{align}
\sum_{j\not\in S} |b_j| \le \xi_0 \sum_{j\in S} |b_j| . \label{comp-cond-eq}
\end{align}
Throughout, $|S|$ denotes the size of a set $S$.
By Cauchy--Schwartz inequality, this compatibility condition is implied by (hence weaker than) a restricted eigenvalue condition (Bickel et al.~2009) such that
$\nu_0^2 (\sum_{j\in S} b_j^2) \le b^\T \Sigma b$
for any $b\in \mathbb R^{1+k}$ satisfying (\ref{comp-cond-eq}).

\begin{ass} \label{ass-RCAL}
Suppose that the following conditions are satisfied:
\begin{itemize} \addtolength{\itemsep}{-.05in}
\item[(i)] $\max_{j=0,1,\ldots,p} |f_j(X)| \le C_0$ almost surely for a constant $C_0 \ge 1$;

\item[(ii)] $\bar \myeta^1_{\mbox{\tiny CAL}}(X) \ge B_0$ almost surely for a constant $B_0 \in \mathbb R$, that is,
$\pi(X; \bar\gamma^1_{\mbox{\tiny CAL}})$ is bounded from below by $( 1+\me^{-B_0})^{-1}$;

\item[(iii)] the compatibility condition holds for $\Sigma_\gamma $ with
the subset $S_\gamma = \{0\} \cup \{j: \bar\gamma^1_{\mbox{\tiny CAL},j} \not=0, j=1,\ldots,p\}$ and some constants $\nu_0>0$ and $\xi_0>1$,
where $\Sigma_\gamma= E [ T w(X; \bar \gamma^1_{\mbox{\tiny CAL}}) f(X) f^\T (X)]$ is the Hessian of $E\{\ell_{\mbox{\tiny CAL}}(\gamma)\}$
at $\gamma=\bar\gamma^1_{\mbox{\tiny CAL}}$ and $w(X;\gamma) = \me^{-\gamma^\T f(X)}$;

\item[(iv)] $|S_\gamma| \lambda_0 \le \eta_0$ for a sufficiently small constant $\eta_0 >0$, depending only on $(A_0,C_0,\xi_0,\nu_0)$.
\end{itemize}
\end{ass}

\begin{thm}[Tan 2017] \label{thm-RCAL}
Suppose that Assumption~\ref{ass-RCAL} holds. Then for $A_0 > (\xi_0+1)/(\xi_0-1)$, we have
with probability at least $1-4\epsilon$,
\begin{align}
D^\dag_{\mbox{\tiny CAL}} ( \hat \myeta^1_{\mbox{\tiny RCAL}}, \bar \myeta^1_{\mbox{\tiny CAL}} )+(A_0-1) \lambda_0 \| \hat\gamma^1_{\mbox{\tiny RCAL}} -\bar\gamma^1_{\mbox{\tiny CAL}}  \|_1
\le  M_0 |S_\gamma| \lambda_0^2,  \label{thm-RCAL-eq}
\end{align}
where $M_0>0$ is a constant depending only on $(A_0, C_0, B_0, \xi_0, \nu_0, \eta_0)$.
\end{thm}

\begin{rem} \label{rem-RCAL}
We provide comments about the conditions involved.
First, Assumption~\ref{ass-RCAL}(iii) can be justified from a compatibility condition for the Gram matrix $E \{ f(X) f^\T (X)\}$ in conjunction with
additional conditions such as for some constant $\tau_0 >0$,
\begin{align}
b ^\T E \{ f(X) f^\T (X)\} b  \le (b^\T \Sigma_\gamma b) / \tau_0 ,  \quad \forall\, b \in \mathbb R^{1+p}. \label{compat-add}
\end{align}
For example, (\ref{compat-add}) holds provided that $\pi^*(X)$ is bounded from below by a positive constant and $\pi(X; \bar\gamma^1_{\mbox{\tiny CAL}})$ is bounded away from 1.
%from above by $( 1+\me^{-B_1})^{-1}$ for a constant $B_1 \in \mathbb R$.
But it is also possible that Assumption~\ref{ass-RCAL}(iii) is satisfied even if (\ref{compat-add}) does not hold for any $\tau_0>0$.
Therefore, Assumption 1 requires that $\pi(X; \bar\gamma^1_{\mbox{\tiny CAL}})$ is bounded away from 0, but may not be bounded away from 1.
Second, Assumption~\ref{ass-RCAL}(iv) can be relaxed to only require that $|S_\gamma| \lambda_0^2$ is sufficiently small, albeit under stronger conditions, for example,
the variables $f_1(X),\ldots,f_p(X)$ are jointly (not just marginally) sub-gaussian (Huang \& Zhang 2012; Negahban et al.~2012).
On the other hand, Assumption~\ref{ass-RCAL}(iv) is already weaker than the sparsity condition,
$|S_\gamma| \log(p) = o(n^{1/2})$, which is needed for obtaining valid confidence intervals for $\mu^1$
from existing works (Belloni et al.~2014) and our later results.
\end{rem}

\begin{rem} \label{rem2-RCAL}
For the Hessian $\Sigma_\gamma$, the weight $w(X;\bar \gamma^1_{\mbox{\tiny CAL}})$ with $\bar \gamma^1_{\mbox{\tiny CAL}}$ replaced by
$\hat\gamma^1_{\mbox{\tiny RCAL}}$ is identical to that used in the weighted least-square loss (\ref{loss-WL}) to define $\hat\alpha^1_{\mbox{\tiny RWL}}$,
that is, $w(X; \hat\gamma^1_{\mbox{\tiny RCAL}}) =\{1-\hat\pi^1_{\mbox{\tiny RCAL}}(X)\}/\hat\pi^1_{\mbox{\tiny RCAL}}(X)$.
The Hessian of
$\ell_{\mbox{\tiny CAL}}( \gamma)$ at $\bar \gamma^1_{\mbox{\tiny CAL}}$ is also the same as the Hessian of $\ell_{\mbox{\tiny WL}}(\alpha^1; \bar \gamma^1_{\mbox{\tiny CAL}})$ in $\alpha^1$.
As later discussed in Section~\ref{sec:discussion}, this coincidence is a consequence of the construction of the loss functions $\ell_{\mbox{\tiny CAL}}(\gamma)$ and $\ell_{\mbox{\tiny WL}}(\alpha^1; \gamma)$
in (\ref{loss-CAL}) and (\ref{loss-WL}).
\end{rem}

Now we turn to the regularized weighted least-squares estimator $\hat\alpha^1_{\mbox{\tiny RWL}}$.
We develop a new strategy of inverting a quadratic inequality
to address the dependency of $\hat\alpha^1_{\mbox{\tiny RWL}}$ on $\hat\gamma^1_{\mbox{\tiny RCAL}}$
and establish convergence of  $\hat\alpha^1_{\mbox{\tiny RWL}}$ under similar conditions as needed for Lasso penalized unweighted least-squares estimators in high-dimensional settings.
The error bound obtained, however, depends on the sparsity size $|S_\gamma|$ and various constants in Assumption~\ref{ass-RCAL}.

For theoretical analysis, the tuning parameter $\lambda$ used in (\ref{reg-wl-loss}) for defining $\hat\alpha^1_{\mbox{\tiny RWL}}$ is specified as
$\lambda = A_1 \lambda_1$, with a constant $A_1>1$ and
\begin{align*}
\lambda_1 = \max\left\{\lambda_0, \,  \me^{-B_0} C_0 \textstyle{\sqrt{8(D_0^2+D_1^2)}\sqrt{ \log\{(1+p)/\epsilon\} / n}} \right\} ,
\end{align*}
where $0 < \epsilon <1$ is a tail probability for the error bound, $(C_0,B_0)$ are from Assumption~\ref{ass-RCAL}, and
$(D_0,D_1)$ are from Assumption~\ref{ass-RWL} below.
With possible model misspecification, the target value $\bar\alpha^1_{\mbox{\tiny WL}}$ is defined
as a minimizer of the expected loss $E \{ \ell_{\mbox{\tiny WL}} (\alpha^1; \bar\gamma^1_{\mbox{\tiny CAL}} ) \}$ as in Section~\ref{sec:overview}.
The following result gives the convergence of $\hat\alpha^1_{\mbox{\tiny RWL}}$ to $\bar\alpha^1_{\mbox{\tiny WL}}$
in the $L_1$ norm $\| \hat\alpha^1_{\mbox{\tiny RWL}} -\bar\alpha^1_{\mbox{\tiny WL}}\|_1$ and the weighted (in-sample) prediction error defined as
\begin{align}
Q_{\mbox{\tiny WL}} (  \hat m^1_{\mbox{\tiny RWL}} , \bar m^1_{\mbox{\tiny WL}};  \bar \gamma^1_{\mbox{\tiny CAL}} ) =
\tilde E \left[T w(X; \bar\gamma^1_{\mbox{\tiny CAL}}) \{ \hat m^1_{\mbox{\tiny RWL}} (X)- \bar m^1_{\mbox{\tiny WL}} (X) \}^2 \right], \label{Q-def}
\end{align}
where $\hat m^1_{\mbox{\tiny RWL}}(X) = \hat \alpha^{\one\T}_{\mbox{\tiny RWL}} f(X)$ and $\bar m^1_{\mbox{\tiny WL}}(X) = \bar \alpha^{\one\T}_{\mbox{\tiny WL}} f(X)$.
In fact, $Q_{\mbox{\tiny WL}} (  \hat m^1_{\mbox{\tiny RWL}} , \bar m^1_{\mbox{\tiny WL}} ;  \bar \gamma^1_{\mbox{\tiny CAL}})$
is the symmetrized Bregman divergence between $\hat m^1_{\mbox{\tiny RWL}}(X)$ and $\bar m^1_{\mbox{\tiny WL}}(X)$
associated with the loss $\kappa_{\mbox{\tiny WL}} (h; \bar\gamma^1_{\mbox{\tiny CAL}} ) = \tilde E [ Tw(X; \bar\gamma^1_{\mbox{\tiny CAL}})   \{Y - h(X) \}^2]/2$.
See Section~\ref{sec:glm-OR} for further discussion.

\begin{ass} \label{ass-RWL}
Suppose that the following conditions are satisfied:
\begin{itemize}\addtolength{\itemsep}{-.05in}
\item[(i)] $Y^1 - \bar m^1_{\mbox{\tiny WL}}(X)$ is uniformly sub-gaussian given $X$: $ D_0^2 E( \exp[ \{ Y^1 - \bar m^1_{\mbox{\tiny WL}}(X)\}^2  /D_0^2] -1 |X) \le D_1^2$
for some positive constants $(D_0,D_1)$;

\item[(ii)] the compatibility condition holds for $\Sigma_\gamma$ with
the subset $S_\alpha= \{0\} \cup \{j: \bar\alpha^1_{\mbox{\tiny WL},j} \not=0, j=1,\ldots,p\}$ and some constants $\nu_1>0$ and $\xi_1>1$;

\item[(iii)] $(1+\xi_1)^2 \nu_1^{-2} |S_\alpha| \lambda_1 \le \eta_1$ for a constant $0<\eta_1 <1$.
\end{itemize}
\end{ass}

\begin{thm} \label{thm-RWL}
Suppose that linear outcome model (\ref{lm-OR}) is used, $A_0 > (\xi_0+1)/(\xi_0-1)$, $A_1 > (\xi_1+1)/(\xi_1-1)$, and
Assumptions~\ref{ass-RCAL} and \ref{ass-RWL} hold.
If $\log\{(1+p)/\epsilon\} / n\le 1$, then we have
with probability at least $1-8\epsilon$,
\begin{align}
&  Q_{\mbox{\tiny WL}}(  \hat m^1_{\mbox{\tiny RWL}} , \bar m^1_{\mbox{\tiny WL}};  \bar \gamma^1_{\mbox{\tiny CAL}} )
+ \me^{\eta_{01}} (A_1-1) \lambda_1 \| \hat \alpha^1_{\mbox{\tiny RWL}} - \bar \alpha^1_{\mbox{\tiny WL}} \|_1 \nonumber \\
& \le  \me^{4\eta_{01}} \xi_2^{-2} \left( M_{01} |S_\gamma| \lambda_0^2\right) +
 \me^{2\eta_{01}} \xi_3^2 \left( \nu_2^{-2} |S_\alpha| \lambda_1^2  \right) ,  \label{thm-RWL-eq}
\end{align}
where $\xi_2 =1-2A_1/\{(\xi_1+1)(A_1-1)\} \in (0,1]$, $\xi_3 = (\xi_1+1)(A_1-1)$, and
$\nu_2 = \nu_1 (1-\eta_1)^{1/2}$, depending only on $(A_1,\xi_1,\nu_1, \eta_1)$, and
$M_{01} =(D_0^2 + D_1^2) (\me^{\eta_{01}} M_0 + \eta_{02}) + (D_0^2 + D_0 D_1) \eta_{02}$, $\eta_{01} = (A_0-1)^{-1} M_0  \eta_0 C_0$,
and $\eta_{02} = (A_0-1)^{-2} M_0^2  \eta_0 $, depending only on $(A_0, C_0, B_0, \xi_0, \nu_0, \eta_0)$ in Theorem~\ref{thm-RCAL} and $(D_0,D_1)$.
\end{thm}

\begin{rem} \label{rem-RWL}
Assumption~\ref{ass-RWL}(ii) is concerned about the same matrix $\Sigma_\gamma$ as in Assumption~\ref{ass-RCAL}(iii), but
with the sparsity subset $S_\alpha$ from $\bar \alpha^1_{\mbox{\tiny WL}}$ instead of $S_\gamma$ from $\bar \gamma^1_{\mbox{\tiny CAL}}$.
The matrix $\Sigma_\gamma$ is also the Hessian of the expected loss $E\{ \ell_{\mbox{\tiny WL}}(\alpha^1; \bar \gamma^1_{\mbox{\tiny CAL}})\}$ at
$\alpha^1 = \bar \alpha^1_{\mbox{\tiny WL}}$, for reasons mentioned in Remark~\ref{rem2-RCAL}.
Assumptions~\ref{ass-RWL}(ii)--(iii) are combined to derive a compatibility condition for the sample matrix
$\tilde \Sigma_\gamma =\tilde E [ T w(X; \bar \gamma^1_{\mbox{\tiny CAL}}) f(X) f^\T (X)]$.
Assumption~\ref{ass-RWL}(iii) can be relaxed such that $|S_\alpha| \lambda_1^2$ is sufficiently small under further side conditions,
but it is already weaker than the sparsity condition, $|S_\alpha| \log(p) = o(n^{1/2})$,
later needed for valid confidence intervals for $\mu^1$. Essentially, the conditions in Assumption~\ref{ass-RWL} are comparable to those for high-dimensional analysis of standard Lasso estimators
(Bickel et al.~2009; Buhlmann \& van de Geer 2011).
\end{rem}

\begin{rem} \label{rem2-RWL}
One of the key steps in our proof is to upper-bound the product
\begin{align}
(\hat\alpha^1_{\mbox{\tiny RWL}} - \bar \alpha^1_{\mbox{\tiny WL}})^\T \tilde E \left[ T w(X; \hat\gamma^1_{\mbox{\tiny RCAL}}) \{ Y - \bar m^1_{\mbox{\tiny WL}} (X)\} f(X) \right] . \label{key-prod}
\end{align}
If $\hat\gamma^1_{\mbox{\tiny RCAL}}$ were replaced by $\bar\gamma^1_{\mbox{\tiny CAL}}$, then it is standard to use the following bound,
\begin{align}
& (\hat\alpha^1_{\mbox{\tiny RWL}} - \bar \alpha^1_{\mbox{\tiny WL}})^\T \tilde E \left[ T w(X; \bar\gamma^1_{\mbox{\tiny CAL}}) \{ Y - \bar m^1_{\mbox{\tiny WL}} (X)\} f(X) \right] \label{key-prod2} \\
& \le \| \hat\alpha^1_{\mbox{\tiny RWL}} - \bar \alpha^1_{\mbox{\tiny WL}} \|_1 \times \| \tilde E \left[ T w(X; \bar\gamma^1_{\mbox{\tiny CAL}}) \{ Y - \bar m^1_{\mbox{\tiny WL}} (X)\} f(X) \right] \|_\infty. \nonumber
\end{align}
To handle the dependency on $\hat\gamma^1_{\mbox{\tiny RCAL}}$, our strategy is to derive an upper bound of the difference between (\ref{key-prod}) and (\ref{key-prod2}),
depending on $Q_{\mbox{\tiny WL}}(  \hat m^1_{\mbox{\tiny RWL}} , \bar m^1_{\mbox{\tiny WL}} ;  \bar \gamma^1_{\mbox{\tiny CAL}})$,
which we seek to control. Carrying this bound leads to a quadratic inequality in
$Q_{\mbox{\tiny WL}}(  \hat m^1_{\mbox{\tiny RWL}} , \bar m^1_{\mbox{\tiny WL}} ;  \bar \gamma^1_{\mbox{\tiny CAL}})$, which can be inverted to obtain an explicit bound
on $Q_{\mbox{\tiny WL}}(  \hat m^1_{\mbox{\tiny RWL}} , \bar m^1_{\mbox{\tiny WL}} ;  \bar \gamma^1_{\mbox{\tiny CAL}})$.
The resulting error bound (\ref{thm-RWL-eq}) is of order $( |S_\gamma| + |S_\alpha|) \log(p)/n$, much sharper
than what we could obtain using other approaches, for example, directly bounding $\| \tilde E [ T w(X; \hat\gamma^1_{\mbox{\tiny RCAL}}) \{ Y - \bar m^1_{\mbox{\tiny WL}} (X)\} f(X) ] \|_\infty$.
\end{rem}

Finally, we study the proposed estimator $\hat \mu^1(\hat m^1_{\mbox{\tiny RWL}}, \hat\pi^1_{\mbox{\tiny RCAL}}) $ for $\mu^1$, depending on the regularized estimators
$\hat \gamma^1_{\mbox{\tiny RCAL}}$ and $\hat \alpha^1_{\mbox{\tiny RWL}}$ from logistic propensity score model (\ref{model-PS})
and linear outcome regression model (\ref{lm-OR}).
The following result gives an error bound for $\hat \mu^1(\hat m^1_{\mbox{\tiny RWL}}, \hat\pi^1_{\mbox{\tiny RCAL}}) $, allowing that
both models (\ref{model-PS}) and (\ref{lm-OR}) may be misspecified.

\begin{thm} \label{thm-mu}
Under the conditions of Theorem~\ref{thm-RWL}, if $\log\{(1+p)/\epsilon\} / n\le 1$, then we have
with probability at least $1-10\epsilon$,
\begin{align}
& \left| \hat \mu^1(\hat m^1_{\mbox{\tiny RWL}}, \hat\pi^1_{\mbox{\tiny RCAL}}) - \hat \mu^1(\bar m^1_{\mbox{\tiny WL}}, \bar\pi^1_{\mbox{\tiny CAL}})  \right| \nonumber \\
& \le M_{11} |S_\gamma| \lambda_0^2 + M_{12} |S_\gamma| \lambda_0 \lambda_1 + M_{13} |S_\alpha| \lambda_0 \lambda_1 , \label{thm-mu-eq}
\end{align}
where $M_{11} = M_{13} + \sqrt{D_0^2 + D_1^2} \me^{\eta_{01}} ( \me^{\eta_{01}}  M_0/2 + \eta_{02}) $,
$M_{12} = (A_0-1)^{-1} M_0$, $M_{13} = A_0 (A_1-1)^{-1} M_1$,
and $M_1$ is a constant such that the right hand side of (\ref{thm-RWL-eq}) in Theorem~\ref{thm-RWL} is upper-bounded by  $\me^{\eta_{01}} M_1 (|S_\gamma| \lambda_0\lambda_1 + |S_\alpha|\lambda_1^2 )$.
\end{thm}

Theorem~\ref{thm-mu} shows that $\hat \mu^1(\hat m^1_{\mbox{\tiny RWL}}, \hat\pi^1_{\mbox{\tiny RCAL}})$ is doubly robust for $\mu^1$
provided $ ( |S_\gamma| + |S_\alpha|) \lambda_1^2 = o(1)$, that is, $( |S_\gamma| + |S_\alpha|) \log(p) = o(n)$.
In addition, Theorem~\ref{thm-mu} gives the $n^{-1/2}$ asymptotic expansion (\ref{new-expansion}) provided $n^{1/2} ( |S_\gamma| + |S_\alpha|) \lambda_1^2 = o(1)$,
that is, $( |S_\gamma| + |S_\alpha|) \log(p) = o(n^{1/2})$.
To obtain valid confidence intervals for $\mu^1$ via the Slutsky theorem,
the following result gives the convergence of the variance estimator $\hat V$ to $V$, as defined in Proposition~\ref{pro-mu-var}, allowing that
both models (\ref{model-PS}) and (\ref{lm-OR}) may be misspecified. For notational simplicity,
denote $\hat\varphi = \varphi(T,Y,X; \hat m^1_{\mbox{\tiny RWL}}, \hat \pi^1_{\mbox{\tiny RCAL}})$ and
$\hat \varphi_c = \hat\varphi - \hat\mu^1(\hat m^1_{\mbox{\tiny RWL}}, \hat \pi^1_{\mbox{\tiny RCAL}})$
such that $\hat V = \tilde E ( \hat \varphi_c^2)$.
Similarly, denote
$\bar\varphi = \varphi(T,Y,X; \bar m^1_{\mbox{\tiny WL}}, \bar \pi^1_{\mbox{\tiny CAL}})$ and
$\bar \varphi_c = \bar\varphi - \hat\mu^1(\bar m^1_{\mbox{\tiny WL}}, \bar \pi^1_{\mbox{\tiny CAL}})$
such that $ V = E ( \bar \varphi_c^2)$.

\begin{thm} \label{thm-var}
Under the conditions of Theorem~\ref{thm-RWL}, if $\log\{(1+p)/\epsilon\} / n\le 1$, then we have
with probability at least $1-10\epsilon$,
\begin{align}
\left| \tilde E \left(\hat \varphi_c^2 - \bar \varphi_c^2 \right) \right|
\le 2 M_{14} \{\tilde E(\bar \varphi_c^2 )\}^{1/2}  (|S_\gamma|\lambda_0 + |S_\alpha|\lambda_1)  + M_{14} (|S_\gamma|\lambda_0 + |S_\alpha|\lambda_1)^2 , \label{thm-var-eq1}
\end{align}
where $M_{14}$ is a positive constant depending only on $(A_0, C_0, B_0, \xi_0, \nu_0, \eta_0)$ in Theorem~\ref{thm-RCAL}
and $(A_1, D_0, D_1, \xi_1, \nu_1, \eta_1)$ in Thorem~\ref{thm-RWL}.
If, in addition, condition (\ref{compat-add}) holds, then we have with probability at least $1-12\epsilon$,
\begin{align}
\left| \tilde E \left(\hat \varphi_c^2 - \bar \varphi_c^2 \right) \right|
\le 2 M_{15} \{\tilde E(\bar \varphi_c^2 )\}^{1/2}  (|S_\gamma|\lambda_0 \lambda_1 + |S_\alpha|\lambda_1^2)^{1/2}  + M_{15} (|S_\gamma|\lambda_0 \lambda_1 + |S_\alpha|\lambda_1^2) , \label{thm-var-eq2}
\end{align}
where $M_{15}$ is a positive constant, depending on $\tau_0$ from (\ref{compat-add}) as well as $(A_0, C_0, B_0, $ $\xi_0, \nu_0, \eta_0)$ and
$(A_1, D_0, D_1, \xi_1, \nu_1, \eta_1)$.
\end{thm}

\begin{rem} \label{rem-var}
Theorem~\ref{thm-var} provides two rates of convergence for $\hat V$ under different conditions.
Inequality (\ref{thm-var-eq1}) shows that $\hat V$ is a consistent estimator of $V$, that is, $\hat V - V =o_p(1)$, provided $( |S_\gamma| + |S_\alpha|) (\log p)^{1/2} = o(n^{1/2})$.
Technically, consistency of $\hat V$ is sufficient for applying Slutsky theorem to establish confidence intervals for $\mu^1$ in Proposition~\ref{pro-mu-var}(iii).
With additional condition (\ref{compat-add}), inequality (\ref{thm-var-eq2}) shows that $\hat V$ achieves the parametric rate of convergence,
$\hat V - V =o_p(n^{-1/2})$, provided $( |S_\gamma| + |S_\alpha|) \log(p) = o(n^{1/2})$.
\end{rem}

\begin{rem} \label{rem2-var}
Combining Theorems~\ref{thm-mu}--\ref{thm-var} directly leads to Proposition~\ref{pro-mu-var}, which gives doubly robust confidence intervals of $\mu^1$.
In addition, a broader interpretation can also be accommodated. All the results, Theorems~\ref{thm-RCAL}--\ref{thm-var}, are developed to remain valid
in the presence of misspecification of models (\ref{model-PS}) and (\ref{lm-OR}),
similarly as in classical theory of estimation with misspecified models (e.g., White 1982; Manski 1988).
If both models (\ref{model-PS}) and (\ref{lm-OR}) may be misspecified, then $ \hat \mu^1(\hat m^1_{\mbox{\tiny RWL}}, \hat\pi^1_{\mbox{\tiny RCAL}}) \pm z_{c/2} \sqrt{\hat V/n}$
is an asymptotic $(1-c)$ confidence interval for the target value $\bar \mu^1= E ( \bar \varphi)$, which in general differs from the true value $\mu^1$.
By comparison, the standard estimator $\hat \mu^1(\hat m^1_{\mbox{\tiny RML}}, \hat\pi^1_{\mbox{\tiny RML}})$ can be shown to converge to
a target value, different from $\mu^1$ as well as $\bar \mu^1$ in the presence of model misspecification.
But it seems difficult to obtain valid confidence intervals based on $\hat \mu^1(\hat m^1_{\mbox{\tiny RML}}, \hat\pi^1_{\mbox{\tiny RML}})$
under similar conditions as in our results, because (\ref{decomp-2}) and (\ref{decomp-3}) are then $O_p(\{\log (p)/n\}^{1/2} )$
if either model (\ref{model-PS}) or (\ref{lm-OR}) is misspecified.
\end{rem}

\subsection{Using generalized linear outcome models} \label{sec:glm-OR}

In this section, we turn to the situation where a generalized linear model is used for
outcome regression together with a logistic propensity score model,
and develop appropriate methods and theory for obtaining confidence intervals for $\mu^1$ in high-dimensional settings.

A technical complication compared with the situation of a linear outcome model in Section~\ref{sec:linear-OR} is that
the reasoning outlined through (\ref{mu-decomp})--(\ref{grad-alpha}) for deriving doubly robust confidence intervals for $\mu^1$ does not directly
hold with a non-linear outcome model, where the second term of (\ref{mu-decomp}) does not in general reduce to the simple product in (\ref{decomp-2}).
There are, however, different approaches that can be used to derive model-assisted confidence intervals, that is, satisfying either property (G2) or (G3) described in Section~\ref{sec:overview}.
For concreteness, we focus on a PS based, OR assisted approach to obtain confidence intervals with property (G2), that is, being
valid if the propensity score model used is correctly specified but the outcome regression model may be misspecified.
See Section~\ref{sec:discussion} for further discussion of related issues.

Consider a logistic propensity score model (\ref{model-PS}) and a generalized linear outcome model,
\begin{align}
E (Y | T=1, X) = m_1(X; \alpha^1 ) = \psi\{ \alpha^{\one\T} f(X) \}, \label{glm-OR}
\end{align}
that is, model (\ref{model-OR}) with the vector of covariate functions $g^1(X)$ taken to be the same as $f(X)$ in model (\ref{model-PS}).
This choice of covariate functions can be more justified than in the setting of Section~\ref{sec:linear-OR}, because OR model (\ref{glm-OR})
plays an assisting role when confidence intervals for $\mu^1$ are concerned.
Our point estimator of $\mu^1$ is $\hat \mu^1(\hat m^1_{\mbox{\tiny RWL}}, \hat\pi^1_{\mbox{\tiny RCAL}})$ as defined in (\ref{mu-estimator}),
where $\hat\pi^1_{\mbox{\tiny RCAL}}(X) = \pi( X; \hat \gamma^1_{\mbox{\tiny RCAL}})$ and
$\hat m^1_{\mbox{\tiny RWL}}(X) = m_1(X; \hat \alpha^1_{\mbox{\tiny RWL}} )$.
The estimator $\hat\gamma^1_{\mbox{\tiny RCAL}}$ is a regularized calibrated estimator of $\gamma$ from Tan (2017) as in Section~\ref{sec:linear-OR}. But
$\hat \alpha^1_{\mbox{\tiny RWL}} $ is a regularized weighted likelihood estimator of $\alpha^1$, defined as a minimizer of
\begin{align}
\ell_{\mbox{\tiny RWL}} (\alpha^1; \hat\gamma^1_{\mbox{\tiny RCAL}} ) &=\ell_{\mbox{\tiny WL}} (\alpha^1; \hat\gamma^1_{\mbox{\tiny RCAL}}) + \lambda  \|\alpha^1_{1:p}\|_1, \label{reg-wl2-loss}
\end{align}
where $\ell_{\mbox{\tiny WL}} (\alpha^1 ; \hat\gamma^1_{\mbox{\tiny RCAL}}) $ is the weighted likelihood loss as follows, with
$w(X;\gamma) = \{1-\pi(X;\gamma)\}/ $ $\pi(X; \gamma)$ $ = \me^{-\gamma^\T f(X)}$ for logistic model (\ref{model-PS}),
\begin{align}
\ell_{\mbox{\tiny WL}} (\alpha^1; \hat\gamma^1_{\mbox{\tiny RCAL}}) = \tilde E \Big( T w(X;\hat\gamma^1_{\mbox{\tiny RCAL}})
\left[ -Y \alpha^{\one\T} f(X) + \Psi\{ \alpha^{\one\T} f(X)\} \right] \Big),  \label{loss-WL2}
\end{align}
and $\|\alpha^1_{1:p}\|_1$ is the $L_1$ norm of $\alpha^1_{1:p}$ and $\lambda \ge 0$ is a tuning parameter.
The regularized weighted least-squares estimator $\hat \alpha^1_{\mbox{\tiny RWL}}$ used in Section~\ref{sec:linear-OR}
is recovered in the special case of the identity link, $\psi(u)=u$ and $\Psi(u) = u^2/2$.
In addition, the Kuhn--Tucker--Karush condition for minimizing (\ref{reg-wl2-loss}) remains the same as in (\ref{ineq-WL-1})--(\ref{ineq-WL-2}),
and hence the estimator $\hat \mu^1(\hat m^1_{\mbox{\tiny RWL}}, \hat\pi^1_{\mbox{\tiny RCAL}})$ can be put in the prediction form (\ref{mu-prediction-form}),
which ensures the boundedness property that $\hat \mu^1(\hat m^1_{\mbox{\tiny RWL}}, \hat\pi^1_{\mbox{\tiny RCAL}})$ always falls in the range of the observed outcomes $Y_i$
in the treated group ($T_i=1$) and the predicted values $\hat m^1_{\mbox{\tiny RWL}}(X_i)$ in the untreated group ($T_i=0$).

With possible model misspecification, the target value $\bar \alpha^1_{\mbox{\tiny WL}}$ is defined as a minimizer of the
expected loss $E\{ \ell_{\mbox{\tiny WL}}( \alpha^1; \bar \gamma^1_{\mbox{\tiny CAL}})\}$. From a functional perspective, we write
$ \ell_{\mbox{\tiny WL}}( \alpha^1;  \gamma) =  \kappa_{\mbox{\tiny WL}} ( \alpha^{\one\T} f; \gamma)$, where
for a function $\myeta(x)$ which may not be in the form $\alpha^{\one\T} f$,
\begin{align*}
\kappa _{\mbox{\tiny WL}} ( h; \gamma) =\tilde E \left( T w (X; \gamma)
\left[ -Y \myeta (X) + \Psi\{ \myeta (X)\} \right] \right) .
\end{align*}
As $\kappa _{\mbox{\tiny WL}} ( h; \gamma)$ is convex in $\myeta$ by the convexity of $\Psi()$,
the Bregman divergence associated with
$\kappa _{\mbox{\tiny WL}} ( h; \gamma)$ is defined as
\begin{align*}
D_{\mbox{\tiny WL}} (\myeta^\prime, \myeta; \gamma) = \kappa_{\mbox{\tiny WL}}(\myeta^\prime; \gamma) -  \kappa_{\mbox{\tiny WL}}(\myeta; \gamma)
- \langle \nabla\kappa_{\mbox{\tiny WL}} (\myeta; \gamma), \myeta^\prime - \myeta \rangle,
\end{align*}
where $\nabla\kappa_{\mbox{\tiny WL}} (\myeta; \gamma)$ denotes the gradient of $\kappa_{\mbox{\tiny WL}} (\myeta; \gamma)$ with respect to $(\myeta_1,\ldots,\myeta_n)$ with $\myeta_i = \myeta(X_i)$.
The symmetrized Bregman divergence is easily shown to be
\begin{align}
D^\dag_{\mbox{\tiny WL}} (\myeta^\prime, \myeta; \gamma) & =  D_{\mbox{\tiny WL}} (\myeta^\prime, \myeta; \gamma)  +
D_{\mbox{\tiny WL}} (\myeta, \myeta^\prime; \gamma) \nonumber \\
& = \tilde E \left( T w(X;\gamma) \left[ \psi\{\myeta^\prime(X)\} - \psi\{\myeta(X)\} \right] \{ \myeta^\prime(X)- \myeta(X)\} \right) . \label{symD-WL}
\end{align}
The following result establishes the convergence of $\hat\alpha^1_{\mbox{\tiny RWL}}$ to $\bar\alpha^1_{\mbox{\tiny WL}}$
in the $L_1$ norm $\| \hat\alpha^1_{\mbox{\tiny RWL}} -\bar\alpha^1_{\mbox{\tiny WL}}\|_1$ and the symmetrized Bregman divergence
$D^\dag_{\mbox{\tiny WL}} (  \hat \myeta^1_{\mbox{\tiny RWL}} , \bar \myeta^1_{\mbox{\tiny WL}};  \bar \gamma^1_{\mbox{\tiny CAL}} ) $,
where $\hat \myeta^1_{\mbox{\tiny RWL}}(X) = \hat \alpha^{\one\T}_{\mbox{\tiny RWL}} f(X)$ and $\bar \myeta^1_{\mbox{\tiny WL}}(X) = \bar \alpha^{\one\T}_{\mbox{\tiny WL}} f(X)$.
In the case of the identity link, $\psi(u)=u $, the symmetrized Bregman divergence
$D^\dag_{\mbox{\tiny WL}} (  \hat \myeta^1_{\mbox{\tiny RWL}} , \bar \myeta^1_{\mbox{\tiny WL}};  \bar \gamma^1_{\mbox{\tiny CAL}} ) $
becomes
$Q_{\mbox{\tiny WL}} (  \hat m^1_{\mbox{\tiny RWL}} , \bar m^1_{\mbox{\tiny WL}};  \bar \gamma^1_{\mbox{\tiny CAL}} ) $ in (\ref{Q-def}).
%$D^\dag_{\mbox{\tiny WL}} (\myeta^\prime, \myeta; \gamma) =\tilde E [ ( T w(X;\gamma) \{ \myeta^\prime(X) - \myeta(X)\}^2 ] $
Inequality (\ref{thm-RWL2-eq}) also reduces to (\ref{thm-RWL-eq}) in Theorem~\ref{thm-RWL} with the choices $C_2=1$ and $C_3=\eta_2=\eta_3=0$.

\begin{ass} \label{ass-RWL2}
Assume that $\psi()$ is differentiable and denote $\psi_2(u) = \dif \psi(u)/\dif u$.
Suppose that the following conditions are satisfied:
\begin{itemize} \addtolength{\itemsep}{-.05in}
\item[(i)] $\psi_2\{ \bar \myeta^1_{\mbox{\tiny WL}}(X) \}  \le C_1$ almost surely for a constant $C_1 >0$;

\item[(ii)] $\psi_2\{ \bar \myeta^1_{\mbox{\tiny WL}}(X) \}  \ge C_2$ almost surely for a constant $C_2 >0$;

\item[(iii)] $\psi_2(u) \le \psi_2(u^\prime) \me^{C_3 |u - u^\prime|}$ for any $(u, u^\prime)$, where $C_3 \ge 0$ is a constant.

\item[(iv)] $C_0 C_3 (A_1-1)^{-1} \xi_3^2 \nu_2^{-2} C_2^{-1} |S_\alpha| \lambda_1 \le \eta_2$ for a constant $0 \le \eta_2 <1$
and $C_0 C_3 \me^{3 \eta_{01}}$ $ (A_1-1)^{-1} \xi_2^{-2} C_2^{-1} (M_{01} |S_\gamma| \lambda_0) \le \eta_3$ for a constant $0 \le \eta_3 <1$,
where $(\eta_{01},\nu_2,\xi_2,\xi_3, M_{01} )$ are as in Theorem~\ref{thm-RWL}.
\end{itemize}
\end{ass}

\begin{thm} \label{thm-RWL2}
Suppose that Assumptions~\ref{ass-RCAL}, \ref{ass-RWL}, and \ref{ass-RWL2}(ii)--(iv) hold. If $\log\{(1+p)/\epsilon\} / n\le 1$, then for $A_0 > (\xi_0+1)/(\xi_0-1)$ and $A_1 > (\xi_1+1)/(\xi_1-1)$, we have
with probability at least $1-8\epsilon$,
\begin{align}
&  D^\dag_{\mbox{\tiny WL}}(  \hat m^1_{\mbox{\tiny RWL}} , \bar m^1_{\mbox{\tiny WL}} ) + \me^{\eta_{01}} (A_1-1) \lambda_1 \| \hat \alpha^1_{\mbox{\tiny RWL}} - \bar \alpha^1_{\mbox{\tiny WL}} \|_1 \nonumber \\
& \le  \me^{4\eta_{01}} \xi_4^{-2} \left( M_{01} |S_\gamma| \lambda_0^2  \right) +
 \me^{2\eta_{01}} \xi_3^2 \left(\nu_3^{-2} |S_\alpha| \lambda_1^2 \right)  ,  \label{thm-RWL2-eq}
\end{align}
where
$\xi_4 = \xi_2 (1-\eta_3)^{1/2}C_2^{1/2}$, $\nu_3 = \nu_2^{1/2} (1-\eta_2)^{1/2} C_2^{1/2}$,
and $(\eta_{01}, \nu_2,\xi_2,\xi_3, M_{01} )$ are as in Theorem~\ref{thm-RWL}.
\end{thm}

\begin{rem} \label{rem-RWL2}
We discuss the conditions involved in Theorem~\ref{thm-RWL2}. Assumption~\ref{ass-RWL2}(i) is not needed, but will be used in later results.
Assumption~\ref{ass-RWL2}(iii), adapted from Huang \& Zhang (2012), is used along with
Assumption~\ref{ass-RCAL}(i) to bound the curvature of $D^\dag_{\mbox{\tiny WL}} (\myeta^\prime, \myeta; \bar \gamma^1_{\mbox{\tiny CAL}})$
and then with Assumption~\ref{ass-RWL2}(iv) to achieve a localized analysis when handling a non-quadratic loss function.
Assumption~\ref{ass-RWL2}(ii) is used for two distinct purposes.
First, it is combined with Assumptions~\ref{ass-RWL}(ii)--(iii) to yield a compatibility condition for
$\tilde \Sigma_{\alpha} = \tilde E [ T w(X; \bar\gamma^1_{\mbox{\tiny CAL}}) \psi_2\{\bar \myeta^1_{\mbox{\tiny WL}}(X)\} f(X) f^\T (X)] $,
which is the sample version of
the Hessian of the expected loss $E \{ \ell_{\mbox{\tiny WL}} (\alpha^1;  \bar \gamma^1_{\mbox{\tiny CAL}})\}$ at $\alpha^1 =\bar\alpha^1_{\mbox{\tiny WL}}$,
that is,
$\Sigma_{\alpha} = E [ T w(X; \bar\gamma^1_{\mbox{\tiny CAL}}) \psi_2\{\bar \myeta^1_{\mbox{\tiny WL}}(X)\} f(X) f^\T (X)] $.
Second, Assumption~\ref{ass-RWL2}(ii) is also used in deriving a quadratic inequality to be inverted in our strategy to deal with the dependency of $\hat\alpha^1_{\mbox{\tiny RWL}}$
on $\hat\gamma^1_{\mbox{\tiny RCAL}}$ as mentioned in Remark~\ref{rem2-RWL}.
As seen from the proofs in Supplementary Material,  similar results as in Theorem~\ref{thm-RWL2} can be obtained with Assumption~\ref{ass-RWL2}(ii) replaced by the weaker condition
that for some constant $\tau_1 >0$,
\begin{align*}
b^\T \Sigma_\gamma b \le (b^\T \Sigma_\alpha b ) / \tau_1, \quad \forall\, b \in \mathbb R^{1+p}, %\label{psi2-lower-bound}
\end{align*}
provided that the condition on $A_1$ and Assumption~\ref{ass-RWL2}(iv) are modified accordingly, depending on $\tau_1$. This extension is not pursued here for simplicity.
\end{rem}

Now we study the proposed estimator $\hat \mu^1(\hat m^1_{\mbox{\tiny RWL}}, \hat\pi^1_{\mbox{\tiny RCAL}}) $ for $\mu^1$,
with the regularized estimators
$\hat \gamma^1_{\mbox{\tiny RCAL}}$ and $\hat \alpha^1_{\mbox{\tiny RWL}}$
obtained using logistic propensity score model (\ref{model-PS}) and genealized linear outcome model (\ref{glm-OR}).
Theorem~\ref{thm-mu2} gives an error bound for $\hat \mu^1(\hat m^1_{\mbox{\tiny RWL}}, \hat\pi^1_{\mbox{\tiny RCAL}}) $, allowing that
both models (\ref{model-PS}) and (\ref{glm-OR}) may be misspecified, but depending on
additional terms in the presence of misspecification of model (\ref{model-PS}).
Denote $h(X; \alpha^1) = \alpha^{\one\T} f(X)$ and for $r \ge 0$,
\begin{align*}
& \Lambda_0( r) = \sup_{j=0,1,\ldots,p, \,\|\alpha^1 - \bar\alpha^1_{\mbox{\tiny WL}} \|_1 \le r}
\left| E \left[ \psi_2\{\myeta (X; \alpha^1)\} f_j(X) \left\{ \frac{T}{\bar \pi^1_{\mbox{\tiny CAL}}(X)} -1 \right\}  \right] \right| .
\end{align*}
As a special case, the quantity $\Lambda_0(0)$ is defined as
\begin{align*}
& \Lambda_1 = \sup_{j=0,1,\ldots,p}
\left| E \left[ \psi_2\{\bar \myeta^1_{\mbox{\tiny WL}}(X)\} f_j(X) \left\{ \frac{T}{\bar \pi^1_{\mbox{\tiny CAL}}(X)} -1 \right\}  \right] \right| .
\end{align*}
By the definition of $\bar \gamma^1_{\mbox{\tiny CAL}}$, it holds that $E [ \{T/\bar \pi^1_{\mbox{\tiny CAL}}(X) -1\} f_j(X)]= 0$ for $j=0,1,\ldots,p$ whether or not
model (\ref{model-PS}) is correctly specified. But $\Lambda_0(r)$ is in general either zero or positive respectively if model (\ref{model-PS}) is correctly specified or misspecified,
except in the case of linear outcome model (\ref{lm-OR}) where $\Lambda_0(r)$ is automatically zero because $\psi_2()$ is constant.

\begin{thm} \label{thm-mu2}
Suppose that Assumptions~\ref{ass-RCAL}, \ref{ass-RWL}, and \ref{ass-RWL2} hold.
If $\log\{(1+p)/\epsilon\} / n\le 1$, then for $A_0 > (\xi_0+1)/(\xi_0-1)$ and $A_1 > (\xi_1+1)/(\xi_1-1)$,
we have with probability at least $1-12\epsilon$,
\begin{align}
& \left| \hat \mu^1(\hat m^1_{\mbox{\tiny RWL}}, \hat\pi^1_{\mbox{\tiny RCAL}}) - \hat \mu^1(\bar m^1_{\mbox{\tiny WL}}, \bar\pi^1_{\mbox{\tiny CAL}})  \right| \nonumber \\
& \le M_{21} |S_\gamma| \lambda_0^2 + M_{22} |S_\gamma| \lambda_0 \lambda_1 + M_{23} |S_\alpha| \lambda_0 \lambda_1 + \eta_{11} \Lambda_0 (\eta_{11}) , \label{thm-mu2-eq1}
\end{align}
where $M_{21}$, $M_{22}$, and $M_{23}$ are positive constants, depending only on $(A_0, C_0, B_0, \xi_0, \nu_0, \eta_0)$,
$(A_1, D_0, D_1, \xi_1, \nu_1, \eta_1)$, and $(C_1,C_2,C_3, \eta_2,\eta_3)$,
$\eta_{11} = (A_1-1)^{-1} M_2 (|S_\gamma| \lambda_0 + |S_\alpha| \lambda_1)$, and
$M_2 $ is a constant such that the right hand side of (\ref{thm-RWL2-eq}) is upper-bounded by  $\me^{\eta_{01}} M_2 (|S_\gamma| \lambda_0\lambda_1 + |S_\alpha|\lambda_1^2 )$.
If, in addition, condition (\ref{compat-add}) holds, then we have
with probability at least $1-14\epsilon$,
\begin{align}
& \left| \hat \mu^1(\hat m^1_{\mbox{\tiny RWL}}, \hat\pi^1_{\mbox{\tiny RCAL}}) - \hat \mu^1(\bar m^1_{\mbox{\tiny WL}}, \bar\pi^1_{\mbox{\tiny CAL}})  \right| \nonumber \\
& \le M_{24} |S_\gamma| \lambda_0^2 + M_{25} |S_\gamma| \lambda_0 \lambda_1 + M_{26} |S_\alpha| \lambda_0 \lambda_1  + \eta_{11} \Lambda_1 , \label{thm-mu2-eq2}
\end{align}
where $M_{24}$, $M_{25}$, and $M_{26}$ are positive constants, also depending on $\tau_0$ from (\ref{compat-add}).
\end{thm}

\begin{rem} \label{rem-mu2}
Two different error bounds are obtained in Theorem~\ref{thm-mu2}. Because $\Lambda_0(\eta_{11}) \ge \Lambda_1$, the error bound (\ref{thm-mu2-eq2}) is tighter than (\ref{thm-mu2-eq1}),
but with the additional condition (\ref{compat-add}), which requires that the generalized eigenvalues of $\Sigma_\gamma$ relative to the gram matrix $E \{ f(X) f^\T(X)\}$ is bounded away from 0.
In either case, the result shows that $\hat \mu^1(\hat m^1_{\mbox{\tiny RWL}}, \hat\pi^1_{\mbox{\tiny RCAL}})$ is doubly robust for $\mu^1$ provided
$ (|S_\gamma| + |S_\alpha|) \lambda_1 = o(1)$, that is, $( |S_\gamma| + |S_\alpha|) (\log p)^{1/2} = o(n^{1/2})$.
In addition, the error bounds imply that $\hat \mu^1(\hat m^1_{\mbox{\tiny RWL}}, \hat\pi^1_{\mbox{\tiny RCAL}})$ admits the $n^{-1/2}$ asymptotic expansion (\ref{new-expansion})
provided $( |S_\gamma| + |S_\alpha|) \log (p) = o(n^{1/2})$,
when PS model (\ref{model-PS}) is correctly specified but OR model (\ref{glm-OR}) may be misspecified, because the term involving $\Lambda_0(\eta_{11})$ or
$\Lambda_1$ vanishes as discussed above. Unfortunately, expansion (\ref{new-expansion}) may fail when PS model (\ref{model-PS}) is misspecified.
\end{rem}

Similarly as Theorem~\ref{thm-var}, the following result establishes the convergence of $\hat V$ to $V$ as defined in Proposition~\ref{pro-mu-var}, allowing that
both models (\ref{model-PS}) and (\ref{glm-OR}) may be misspecified.

\begin{thm} \label{thm-var2}
Under the conditions of Theorem~\ref{thm-mu2}, if $\log\{(1+p)/\epsilon\} / n\le 1$, then we have with probability at least $1-12\epsilon$,
\begin{align}
\left| \tilde E \left(\hat \varphi_c^2 - \bar \varphi_c^2 \right) \right|
& \le 2 M_{27} \{\tilde E(\bar \varphi_c^2 )\}^{1/2} \{1+\Lambda_0(\eta_{11})\} (|S_\gamma|\lambda_0 + |S_\alpha|\lambda_1) \nonumber \\
& \quad + M_{27} \{1+\Lambda_0^2(\eta_{11})\} (|S_\gamma|\lambda_0 + |S_\alpha|\lambda_1)^2 , \label{thm-var2-eq1}
\end{align}
where $M_{27}$ is a positive constant depending only on $(A_0, C_0, B_0, $ $ \xi_0, \nu_0, \eta_0)$,
$(A_1, D_0, D_1, \xi_1, \nu_1, \eta_1)$, and $(C_1,C_2,C_3, \eta_2,\eta_3)$.
If, in addition, condition (\ref{compat-add}) holds, then we have with probability at least $1-14\epsilon$,
\begin{align}
\left| \tilde E \left(\hat \varphi_c^2 - \bar \varphi_c^2 \right) \right|
& \le 2 M_{28} \{\tilde E(\bar \varphi_c^2 )\}^{1/2}  \left\{ (|S_\gamma|\lambda_0 \lambda_1 + |S_\alpha|\lambda_1^2)^{1/2} + \Lambda_1 (|S_\gamma|\lambda_0 + |S_\alpha|\lambda_1) \right\} \nonumber \\
& \quad + M_{28} \left\{ (|S_\gamma|\lambda_0 \lambda_1 + |S_\alpha|\lambda_1^2) + \Lambda_1^2 (|S_\gamma|\lambda_0 + |S_\alpha|\lambda_1)^2 \right\}, \label{thm-var2-eq2}
\end{align}
where $M_{28}$ is a positive constant, similar to $M_{27}$ but also depending on $\tau_0$ from (\ref{compat-add}).
\end{thm}

\begin{rem} \label{rem-var2}
Two different rates of convergence are obtained for $\hat V$ in Theorem~\ref{thm-var2}. Similarly as discussed in Remark~\ref{rem-var},
if $( |S_\gamma| + |S_\alpha|) (\log p)^{1/2} = o(n^{1/2})$,
then inequality (\ref{thm-var2-eq1}) implies the consistency of $\hat V$ for $V$, which is sufficient for
applying Slutsky Theorem to establish confidence intervals for $\mu^1$.
With additional condition (\ref{compat-add}), inequality (\ref{thm-var2-eq2}) gives a faster rate of convergence for $\hat V$, which is of order $n^{-1/2}$ provided $( |S_\gamma| + |S_\alpha|) \log (p) = o(n^{1/2})$.
\end{rem}

Combining Theorems~\ref{thm-mu2}--\ref{thm-var2} leads to the following result.

\begin{pro} \label{pro-mu-var2}
Suppose that Assumptions~\ref{ass-RCAL}, \ref{ass-RWL}, and \ref{ass-RWL2} hold, and $( |S_\gamma| + |S_\alpha| ) \log(p) = o(n^{1/2})$.
For sufficiently large constants $A_0$ and $A_1$, if logistic PS model (\ref{model-PS}) is correctly specified but
OR model (\ref{glm-OR}) may be misspecified, then (i)--(iii) in Proposition~\ref{pro-mu-var} hold.
That is, a PS based, OR assisted confidence interval for $\mu^1$ is obtained.
\end{pro}

\begin{rem} \label{rem2-mu-var2}
The conclusion of Proposition~\ref{pro-mu-var2} remains valid if PS model (\ref{model-PS}) is misspecified but only locally such that
$\Lambda_0 (\eta_{11} ) = O( \{\log(p)/n\}^{1/2}) $ or $\Lambda_1 = O(\{\log(p)/n\}^{1/2} )$,
in the case of the error bound (\ref{thm-mu2-eq1}) or (\ref{thm-mu2-eq2}).
Therefore, $ \hat \mu^1(\hat m^1_{\mbox{\tiny RWL}}, \hat\pi^1_{\mbox{\tiny RCAL}}) \pm z_{c/2} \sqrt{\hat V/n}$
can be interpreted as an asymptotic $(1-c)$ confidence interval for the target value $\bar \mu^1= E ( \bar \varphi)$
if model (\ref{model-PS}) is at most locally misspecified but model (\ref{glm-OR}) may be arbitrarily misspecified.
It is an interesting open problem to find broadly valid confidence intervals in the presence of model misspecification
similarly as discussed in Remark~\ref{rem2-var}
when a linear outcome model is used.
\end{rem}

\subsection{Further discussion}  \label{sec:discussion}

\noindent{\bf Estimation of ATE.}
Our theory and methods are presented mainly on estimation of $\mu^1$, but
they can be directly extended for estimating $\mu^0$ and hence ATE, that is, $\mu^1 - \mu^0$.
Consider a logistic propensity score model (\ref{model-PS}) and a generalized linear outcome model,
\begin{align}
E (Y | T=0, X) = m_0(X; \alpha^0 ) = \psi\{ \alpha^{\zero\T} f(X) \}, \label{glm-OR0}
\end{align}
where $f(X)$ is the same vector of covariate functions as in the model (\ref{model-PS}) and
$\alpha^0$ is a vector of unknown parameters.
Our point estimator of ATE is $\hat \mu^1(\hat m^1_{\mbox{\tiny RWL}}, \hat\pi^1_{\mbox{\tiny RCAL}}) - \hat \mu^0(\hat m^0_{\mbox{\tiny RWL}}, \hat\pi^0_{\mbox{\tiny RCAL}})$,
and that of $\mu^0$ is
\begin{align*}
\hat \mu^0(\hat m^0_{\mbox{\tiny RWL}}, \hat\pi^0_{\mbox{\tiny RCAL}}) = \tilde E \left\{ \varphi(Y,1-T,X;\hat m^0_{\mbox{\tiny RWL}}, 1-\hat\pi^0_{\mbox{\tiny RCAL}}) \right\} , %\label{mu0-estimator}
\end{align*}
where $\varphi()$ is defined in (\ref{influence-function}), $\hat\pi^0_{\mbox{\tiny RCAL}}(X) = \pi( X; \hat \gamma^0_{\mbox{\tiny RCAL}})$,
$\hat m^0_{\mbox{\tiny RWL}}(X) = m_0(X; \hat \alpha^0_{\mbox{\tiny RWL}} )$,
and $\hat \gamma^0_{\mbox{\tiny RCAL}}$ and $\hat \alpha^0_{\mbox{\tiny RWL}}$ are defined as follows.
The estimator $\hat \gamma^0_{\mbox{\tiny RCAL}}$ is defined similarly as $\hat \gamma^1_{\mbox{\tiny RCAL}}$, but with the loss function
$\ell_{\mbox{\tiny CAL}} (\gamma)$ in (\ref{loss-CAL}) replaced by
\begin{align*}
\ell^0_{\mbox{\tiny CAL}} (\gamma) &= \tilde E \left\{ (1-T) \me^{\gamma^\T f(X)} - T \gamma^\T f(X) \right\}, %\label{loss-CAL0}
\end{align*}
that is, $T$ and $\gamma$ in $\ell_{\mbox{\tiny CAL}} (\gamma)$ are replaced by $1-T$ and $-\gamma$.
The estimator $\hat \alpha^0_{\mbox{\tiny RWL}}$ is defined similarly as $\hat \alpha^1_{\mbox{\tiny RWL}}$, but with the loss function
$ \ell_{\mbox{\tiny WL}} ( \cdot; \hat \gamma^1_{\mbox{\tiny RCAL}})$ in (\ref{loss-WL2}) replaced by
\begin{align*}
\ell^0_{\mbox{\tiny WL}} (\alpha^0; \hat \gamma^0_{\mbox{\tiny RCAL}}) = \tilde E \Big( (1-T) w^0(X; \hat \gamma^0_{\mbox{\tiny RCAL}})
\left[ -Y \alpha^{\zero\T} g^0(X) + \Psi\{ \alpha^{\zero\T} g^0(X)\} \right] \Big),  %\label{loss-WL0-OR}
\end{align*}
where $w^0( X; \gamma) = \pi(X;\gamma) / \{1-\pi(X;\gamma)\} = \me^{\gamma^\T f(X)}$.
Under similar conditions as in Propositions~\ref{pro-mu-var} and \ref{pro-mu-var2}, the estimator $\hat \mu^0(\hat m^0_{\mbox{\tiny RWL}}, \hat\pi^0_{\mbox{\tiny RCAL}})$
admits the asymptotic expansion
\begin{align}
\hat \mu^0(\hat m^0_{\mbox{\tiny RWL}}, \hat\pi^0_{\mbox{\tiny RCAL}}) =
\tilde E \left\{ \varphi(Y,1-T,X;\bar m^0_{\mbox{\tiny WL}}, 1-\bar \pi^0_{\mbox{\tiny CAL}}) \right\} +o_p (n^{-1/2}), \label{new-expansion0}
\end{align}
where $\bar\pi^0_{\mbox{\tiny RCAL}}(X) = \pi( X; \bar \gamma^0_{\mbox{\tiny RCAL}})$,
$\bar m^0_{\mbox{\tiny RWL}}(X) = m_0(X; \bar \alpha^0_{\mbox{\tiny RWL}} )$,
and $\bar \gamma^0_{\mbox{\tiny RCAL}}$ and $\bar \alpha^0_{\mbox{\tiny RWL}}$ are the target values defined similarly as
 $\bar \gamma^1_{\mbox{\tiny RCAL}}$ and $\bar \alpha^1_{\mbox{\tiny RWL}}$.
Then Wald confidence intervals for $\mu^0$ and ATE cane be derived from  (\ref{new-expansion}) and (\ref{new-expansion0}) similarly as in Propositions~\ref{pro-mu-var} and \ref{pro-mu-var2}
and shown to be either doubly robust in the case of linear outcome models,
or valid if PS model (\ref{model-PS}) is correctly specified but OR models (\ref{glm-OR}) and (\ref{glm-OR0}) may be misspecified for
nonlinear $\psi()$.

An unusual aspect of our approach is that two different estimators of the propensity score are used when
estimating $\mu^0$ and $\mu^1$.
On one hand, the estimators  $\hat\gamma^0_{\mbox{\tiny RCAL}}$ and $\hat\gamma^1_{\mbox{\tiny RCAL}}$
are both consistent, and hence there is no self-contradiction at least asymptotically, when PS model (\ref{model-PS}) is correctly specified.
On the other hand, if model (\ref{model-PS}) is misspecified,
the two estimators may in general have different asymptotic limits,
which can be an advantage from the following perspective.
By definition, the augmented IPW estimators of $\mu^1$ and $\mu^0$ are obtained, depending on
fitted propensity scores within the treated group and untreated groups separately, that is,
$\{\pi(X_i; \gamma^1): T_i=1\}$ and $\{\pi(X_i; \gamma^0): T_i=0\}$.
In the presence of model misspecification, allowing different $\gamma^1$ and $ \gamma^0$
can be helpful in finding suitable approximations of the two sets of propensity scores, without
being constrained by the then-false assumption that they are determined by the same coefficient vector $\gamma^1=\gamma^0$.

\vspace{.1in}
\noindent{\bf Estimation of ATT.}
There is a simple extension of our approach to estimation of ATT,
that is, $\nu^1-\nu^0$ as defined in Section~\ref{sec:setup}.
The parameter $\nu^1 = E(Y^1 | T=1)$ can be directly estimated by $\tilde E (TY) / \tilde E(T)$.
For $\nu^0 = E(Y^1 | T=1)$, our point estimator is
\begin{align*}
\hat \nu^0(\hat m^0_{\mbox{\tiny RWL}}, \hat\pi^0_{\mbox{\tiny RCAL}}) = \tilde E \left\{ \varphi_{\nu^0}(Y, T,X;\hat m^0_{\mbox{\tiny RWL}}, \hat\pi^0_{\mbox{\tiny RCAL}}) \right\} / \tilde E(T),
\end{align*}
where $\hat \pi^0_{\mbox{\tiny RCAL}}(X)$ and $\hat m^0_{\mbox{\tiny RWL}}(X)$ are the same fitted values as used in the estimator
$\hat \mu^0(\hat m^0_{\mbox{\tiny RWL}}, \hat\pi^0_{\mbox{\tiny RCAL}})$ for $\mu^0$,
and $\varphi_{\nu^0} (\cdot ; \hat m_0, \hat \pi)$ is defined as
\begin{align*}
\varphi_{\nu^0} (Y,T,X;\hat m_0, \hat\pi) = \frac{(1-T)\hat\pi(X)}{1-\hat\pi (X)} Y - \left\{\frac{1-T}{1-\hat\pi (X)}-1 \right\} \hat m_0(X) .
\end{align*}
The function $\varphi_{\nu^0} (\cdot ; \hat m_0, \hat \pi)$ can be derived, by substituting fitted values $(\hat m_0, \hat\pi)$ for the true values $(m^*_0, \pi^*)$
in the efficient influence function of $\mu^0$ under a nonparametric model (Hahn 1998). In addition,
the estimator $\tilde E \{ \varphi_{\nu^0}(Y, T,X;\hat m^0 , \hat\pi) \}  $ is also doubly robust:
it remains consistent for $E(TY^0)$ if either $\hat m^0 = m^*_0$ or $\hat \pi= \pi^*$.
In fact, by straightforward calculation, the function $\varphi_{\nu^0} ( )$ is related to $\varphi ( )$ in (\ref{influence-function}) through the simple identify:
\begin{align}
\varphi_{\nu^0} (Y,T,X;\hat m_0, \hat\pi) = \varphi (Y,1-T,X;\hat m_0, 1-\hat\pi) - (1-T)Y. \label{ATT-identity}
\end{align}
As a result, $\hat \nu^0(\hat m^0_{\mbox{\tiny RWL}}, \hat\pi^0_{\mbox{\tiny RCAL}})$ can be equivalently obtained as
\begin{align*}
\hat \nu^0(\hat m^0_{\mbox{\tiny RWL}}, \hat\pi^0_{\mbox{\tiny RCAL}}) & = \left[ \hat \mu^0(\hat m^0_{\mbox{\tiny RWL}}, \hat\pi^0_{\mbox{\tiny RCAL}}) - \tilde E\{ (1-T)Y \} \right]/ \tilde E(T) \\
& = \tilde E \left\{ T \hat m^0_{\mbox{\tiny RWL}}(X) \right\} / \tilde E(T) ,
\end{align*}
where the second step follows from a similar equation for $\hat \mu^0(\hat m^0_{\mbox{\tiny RWL}}, \hat\pi^0_{\mbox{\tiny RCAL}})$ as (\ref{mu-prediction-form}).
Moreover, it can be shown using Eq.~(\ref{ATT-identity}) that
under similar conditions as in Propositions~\ref{pro-mu-var} and \ref{pro-mu-var2}, the estimator $\hat \nu^0(\hat m^0_{\mbox{\tiny RWL}}, \hat\pi^0_{\mbox{\tiny RCAL}})$
admits the asymptotic expansion
\begin{align*}
\hat \nu^0(\hat m^0_{\mbox{\tiny RWL}}, \hat\pi^0_{\mbox{\tiny RCAL}}) - \nu^0 =
\tilde E \left\{ \varphi_{\nu^0}(Y, T,X;\bar m^0_{\mbox{\tiny WL}}, \bar \pi^0_{\mbox{\tiny CAL}}) - T \nu^0 \right\} / \tilde E(T)+o_p (n^{-1/2}),
\end{align*}
similarly as (\ref{new-expansion0}) for $\hat \mu^0(\hat m^0_{\mbox{\tiny RWL}}, \hat\pi^0_{\mbox{\tiny RCAL}})$.
From this expansion, Wald confidence intervals for $\nu^0$ and ATT can be derived and shown to be either doubly robust with linear OR model (\ref{glm-OR0})
or valid at least when PS model (\ref{model-PS}) is correctly specified.

\vspace{.1in}
\noindent{\bf Construction of loss functions.}
We provide additional comments about the construction of loss functions for $\gamma$ and $\alpha^1$ and
alternative approaches when using nonlinear outcome models.
For a linear outcome model (\ref{lm-OR}) as in Section~\ref{sec:overview},
the loss functions $\ell_{\mbox{\tiny CAL}}(\gamma) $ and $\ell_{\mbox{\tiny WL}}(\alpha^1; \gamma)$ are derived such that
their gradients satisfy (\ref{grad-gamma})--(\ref{grad-alpha}), which are in turn obtained as
the coefficients for $\hat\alpha^1 - \bar\alpha^1$ and $\hat\gamma-\bar\gamma$ in the first-order terms (\ref{decomp-2})--(\ref{decomp-3})
from the Taylor expansion~(\ref{mu-decomp}) of $\hat\mu^1(\hat m^1, \hat \pi)$.
Combining the two steps, Eqs.~(\ref{grad-gamma})--(\ref{grad-alpha}) amount to choosing
\begin{align}
& \frac{\partial \ell_{\mbox{\tiny CAL}}(\gamma)}{\partial \gamma} = \frac{\partial}{\partial \alpha^1}
\tilde E \left[ \varphi\{ Y,T,X; m_1(\cdot;\alpha^1), \pi(\cdot; \gamma)\} \right] , \label{grad-gamma2} \\
& \frac{\partial \ell_{\mbox{\tiny WL}}(\alpha^1; \gamma)}{\partial \alpha^1} =  \frac{\partial }{\partial \gamma}
\tilde E \left[ \varphi\{Y,T,X; m_1(\cdot;\alpha^1), \pi(\cdot; \gamma)\} \right]. \label{grad-alpha2}
\end{align}
We say that
the loss function $\ell_{\mbox{\tiny CAL}}(\gamma)$ for $\gamma$ in model (\ref{model-PS}) is calibrated against model (\ref{lm-OR}),
whereas  $\ell_{\mbox{\tiny WL}}(\alpha^1; \gamma)$ for $\alpha^1$ in model (\ref{lm-OR}) is calibrated against model (\ref{model-PS}).
The estimators $\hat \gamma^1_{\mbox{\tiny RCAL}}$ and $\hat \alpha^1_{\mbox{\tiny RWL}}$ are called regularized calibrated estimators of $\gamma$ and $\alpha^1$ respectively.
The pair of equations (\ref{grad-gamma2})--(\ref{grad-alpha2}) also underlie the coincidence of the Hessian of
$\ell_{\mbox{\tiny CAL}}( \gamma)$ at $\bar \gamma^1_{\mbox{\tiny CAL}}$ and that of $\ell_{\mbox{\tiny WL}}(\alpha^1; \bar \gamma^1_{\mbox{\tiny CAL}})$ in $\alpha^1$
with a linear outcome model, as mentioned in Remark~\ref{rem2-RCAL}.

Previously, an augmented IPW estimator $\hat\mu^1(\hat m^1, \hat \pi)$ for $\mu^1$ was proposed in low-dimensional settings
by Kim \& Haziza (2014) and Vermeulen \& Vansteelandt (2015), where $(\hat\alpha^1, \hat\gamma)$ are non-penalized, defined by directly setting
the right-hand sides of (\ref{grad-gamma2})--(\ref{grad-alpha2}) to zero. One of their motivations is to
enable simple calculation of confidence intervals without the need of correcting for estimation of $(\alpha^1,\gamma)$.
Our work generalizes these previous estimators to high-dimensional settings,
where the motivation for using $(\hat \alpha^1_{\mbox{\tiny RWL}}, \hat \gamma^1_{\mbox{\tiny RCAL}})$, instead of
$(\hat \alpha^1_{\mbox{\tiny RML}}, \hat \gamma_{\mbox{\tiny RML}})$
is mainly statistical:
to reduce the variation caused by estimation of $(\alpha^1,\gamma)$ from $O_p(\{\log (p)/n\}^{1/2} )$ to $o_p(n^{-1/2})$ for the estimator
$\hat\mu^1(\hat m^1_{\mbox{\tiny RWL}}, \hat \pi^1_{\mbox{\tiny RCAL}})$, so that
valid confidence intervals for $\mu^1$ can be obtained
even in the presence of model misspecification.

For a possibly nonlinear outcome model (\ref{glm-OR}), the augmented IPW estimator of $\mu^1$
in Kim \& Haziza (2014) and Vermeulen \& Vansteelandt (2015) is also defined as described above.
However, the gradients from the right-hand sides of (\ref{grad-gamma2})--(\ref{grad-alpha2}) become
\begin{align}
\frac{\partial}{\partial \alpha^1}
\tilde E \left[ \varphi\{Y,T,X; m_1(\cdot;\alpha^1), \pi(\cdot; \gamma) \}\right] &= \tilde E \left[ \left\{ 1 - \frac{T}{\pi(X;\gamma)} \right\} \psi_2\{\alpha^{\one\T}f(X)\} f(X) \right] , \label{grad-gamma3} \\
\frac{\partial }{\partial \gamma}
\tilde E \left[ \varphi\{Y,T,X; m_1(\cdot;\alpha^1), \pi(\cdot; \gamma)\} \right] &=  -\tilde E \left[ T \frac{1-\pi(X;\gamma)}{\pi(X;\gamma)} \{Y- m_1(X; \alpha^1)\} f(X) \right], \label{grad-alpha3}
\end{align}
where $\psi_2()$ denotes the derivative of $\psi()$.
The pair of equations obtained by setting (\ref{grad-gamma3})--(\ref{grad-alpha3}) to zero are intrinsically coupled, unless outcome model (\ref{glm-OR}) is linear and
hence the dependency of (\ref{grad-gamma3}) on $\alpha^1$ vanishes. This complication, although mainly computational in low-dimensional settings,
presents a statistical as well as computational obstacle to developing doubly robust confidence intervals with regularized estimation in high-dimensional settings.

The development in Section~\ref{sec:glm-OR} involves using (\ref{grad-gamma}) instead of (\ref{grad-gamma3}) but retaining (\ref{grad-alpha3}),
which lead to the loss functions $\ell_{\mbox{\tiny CAL}} (\gamma)$ in (\ref{loss-CAL}) and $\ell_{\mbox{\tiny WL}} (\alpha^1; \gamma)$ in (\ref{loss-WL2}).
The resulting confidence intervals are PS based, OR assisted, that is, being valid if PS model (\ref{model-PS}) is correctly specified but OR model (\ref{glm-OR}) may be misspecified.
Alternatively, it is possible to develop an OR based, PS assisted approach using the regularized maximum likelihood estimator $\hat\alpha^1_{\mbox{\tiny RML}}$
in conjunction with a regularized estimator of $\gamma$ based on a weighted calibration loss,
\begin{align}
\ell_{\mbox{\tiny WL}} (\gamma; \hat\alpha^1_{\mbox{\tiny RML}}) = \tilde E \left[ \psi_2 \left\{ \hat\alpha^{\one\T}_{\mbox{\tiny RML}} f(X) \right\} \left\{ T \me^{-\gamma^\T f(X)} + (1-T) \gamma^\T f(X) \right\} \right].
\label{loss-WCAL}
\end{align}
The gradient of (\ref{loss-WCAL}) in $\gamma$ is (\ref{grad-gamma3}), with $\alpha^1=\hat\alpha^1_{\mbox{\tiny RML}}$.
Similar results can be established as in Section~\ref{sec:glm-OR}, to provide valid confidence intervals for $\mu^1$
if OR model (\ref{glm-OR}) is correctly specified but PS model (\ref{model-PS}) may be misspecified.
This work can be pursued elsewhere.

\section{Simulation study}

We present a simulation study with the design of Kang \& Schafer (2007) modified and extended to high-dimensional, sparse settings.
It is of interest to empirically compare $\hat \mu^1(\hat m^1_{\mbox{\tiny RML}}, \hat\pi^1_{\mbox{\tiny RML}})$ and
$\hat \mu^1(\hat m^1_{\mbox{\tiny RWL}}, \hat\pi^1_{\mbox{\tiny RCAL}})$ and their associated confidence intervals.

In our implementation,
the penalized loss function (\ref{reg-ml-loss-OR}) or (\ref{reg-ml-loss-PS}) for computing
$\hat \alpha^1_{\mbox{\tiny RML}}$ or $\hat\gamma^1_{\mbox{\tiny RML}}$
or (\ref{reg-cal-loss}), (\ref{reg-wl-loss}), or (\ref{reg-wl2-loss})
for computing
$\hat \alpha^1_{\mbox{\tiny RWL}}$ or $\hat\gamma^1_{\mbox{\tiny RCAL}}$
is minimized for a fixed tuning parameter $\lambda$, using algorithms similar to those in Friedman et al.~(2010), but with the coordinate descent method replaced by an active set method
as in Osborne et al.~(2000) for solving each Lasso penalized least squares problem.
In addition, the penalized loss (\ref{reg-cal-loss}) for computing $\hat\gamma^1_{\mbox{\tiny RCAL}}$ is minimized using the algorithm in Tan (2017),
where a nontrivial Fisher scoring step is involved for quadratic approximation.
The tuning parameter $\lambda$ is determined using 5-fold cross validation
based on the corresponding loss function as follows.

For $k=1,\ldots,5$, let $\mathcal I_k$ be a random subsample of size $n/5$ from $\{1,2,\ldots,n\}$.
For a loss function $\ell (\gamma)$, either $\ell_{\mbox{\tiny ML}}(\gamma)$ in (\ref{loss-ML-PS}) or  $\ell_{\mbox{\tiny CAL}}(\gamma)$ in (\ref{loss-CAL}), denote by $\ell (\gamma; \mathcal I)$ the loss function
obtained when the sample average $\tilde E()$ is computed over only the subsample $\mathcal I$.
The 5-fold cross-validation criterion is defined as
\begin{align*}
\mbox{CV}_5 (\lambda) = \frac{1}{5}\sum_{k=1}^5 \ell ( \hat\gamma_\lambda^{(k)}; \mathcal I_k ),
\end{align*}
where $\hat\gamma^{(k)}_\lambda$ is a minimizer of the penalized loss $\ell(\gamma; \mathcal I^c_k) + \lambda \| \gamma_{1:p} \|_1$ over the subsample $\mathcal I^c_k$ of size $4n/5$,
i.e., the complement to $\mathcal I_k$.
Then $\lambda$ is selected by minimizing $\mbox{CV}_5(\lambda)$ over the discrete set $\{ \lambda^* / 2^j: j=0,1,\ldots,10\}$,
where for $\hat\pi_0=\tilde E(T)$, the value $\lambda^*$ is computed as either
\begin{align*}
\lambda^* = \max_{j=1,\ldots,p} \left| \tilde E \{ (T -\hat\pi_0) f_j(X) \} \right|
\end{align*}
when the likelihood loss (\ref{loss-ML-PS}) is used, or
\begin{align*}
\lambda^* = \max_{j=1,\ldots,p} \left| \tilde E \{ (T/\hat\pi_0 -1) f_j(X) \} \right|
\end{align*}
when the calibration loss (\ref{loss-CAL}) is used. It can be shown that in either case, the penalized loss
$\ell(\gamma) + \lambda \| \gamma_{1:p} \|_1$ over the original sample has a minimum at $\gamma_{1:p} =0$ for all $\lambda \ge \lambda^*$.

For computing $\hat \alpha^1_{\mbox{\tiny RML}}$ or $\hat \alpha^1_{\mbox{\tiny RWL}}$, cross validation is conducted similarly as above using the loss function $\ell_{\mbox{\tiny ML}}(\alpha^1)$ in (\ref{loss-ML-OR}) or
$\ell_{\mbox{\tiny WL}}(\alpha^1; \hat\gamma^1_{\mbox{\tiny RCAL}} )$ in (\ref{loss-WL2}).
In the latter case, $\hat\gamma^1_{\mbox{\tiny RCAL}}$ is determined separately and then fixed during cross validation for computing $\hat \alpha^1_{\mbox{\tiny RWL}}$.

%$Z_1^\prime=\exp(0.5 Z_1)$,
%$Z_2^\prime=10 + \{ 1+ \exp(Z_1) \}^{-1} Z_2$, $Z_3^\prime=(0.04 Z_1 Z_3 + 0.6)^3$, and $Z_4^\prime=(Z_2 + Z_4 + 20)^2$

\subsection{Linear outcome models} \label{sec:sim-linearOR}

Let $X=(X_1,\ldots,X_p)$ be independent variables, where each $X_j$ is $\N(0,1)$ truncated to the interval $(-2.5, 2.5)$ and then standardized to have mean 0 and variance 1.
In addition, let $X^\dag= (X^\dag_1,\ldots, X^\dag_p)$, where $X^\dag_j = X_j$ for $j=5,\ldots,p$, and $X^\dag_1$, $X^\dag_2$, $X^\dag_3$, and $X^\dag_4$ are standardized versions of $\exp(0.5 X_1)$,
$10 + \{ 1+ \exp(X_1) \}^{-1} X_2$, $(0.04 X_1 X_3 + 0.6)^3$, and $(X_2 + X_4 + 20)^2$
to have means 0 and variances 1. The truncation of $X_j$ prevents propensity scores arbitrarily close to 0, and
ensures that the mapping between $X$ and $X^\dag$ are strictly one-to-one.
See the Supplementary Material for calculation to perform the standardization and for scatter plots of $(X^\dag_1,\ldots,X^\dag_4)$.
Consider the following data-generating configurations.
\begin{itemize} \addtolength{\itemsep}{-.05in}

\item[(C1)] Generate $T$ given $X$ from a Bernoulli distribution with \vspace{-.1in}
\begin{align*}
P(T=1 | X) = \{ 1+ \exp( X^\dag_1 - 0.5 X^\dag_2 + 0.25 X^\dag_3 + 0.1 X^\dag_4 ) \}^{-1},
\end{align*}
and, independently, generate $Y^1$ given $X$ from a Normal distribution with variance 1 and mean either (``Linear outcome configuration 1") \vspace{-.1in}
\begin{align*}
E( Y^1 | X) = X^\dag_1 + 0.5 X^\dag_2 + 0.5 X^\dag_3 + 0.5 X^\dag_4 ,
\end{align*}
or (``Linear outcome configuration 2") \vspace{-.1in}
\begin{align*}
E( Y^1 | X) = 0.25 X^\dag_1 + 0.5 X^\dag_2 + 0.5 X^\dag_3 + 0.5 X^\dag_4 .
\end{align*}
The main difference between the two outcome configurations is that
$X_1^\dag$ is both the most important variable influencing the propensity score and that influencing the outcome regression function in the first configuration.

\item[(C2)] Generate $T$ give $X$ as in (C1), but, independently, generate $Y^1$ given $X$ from a Normal distribution with variance 1 and mean either (``Linear outcome configuration 1") \vspace{-.1in}
\begin{align*}
E( Y^1 | X) = X_1 + 0.5 X_2 + 0.5 X_3 + 0.5 X_4,
\end{align*}
or (``Linear outcome configuration 2") \vspace{-.1in}
\begin{align*}
E( Y^1 | X) = 0.25 X_1 + 0.5 X_2 + 0.5 X_3 + 0.5 X_4 .
\end{align*}
As $X_1$ and $X_1^\dag$ are monotone transformations of each other, the variable $X_1^\dag$
remains roughly both the most important variable influencing the propensity score and that influencing the outcome regression function in the first configuration.

\item[(C3)] Generate  $Y^1$ given $X$ as in (C1), but, independently, generate $T$ given $X$ from a Bernoulli distribution with \vspace{-.1in}
\begin{align*}
P(T=1 | X) = \{ 1+ \exp( X_1 - 0.5 X_2 + 0.25 X_3 + 0.1 X_4 ) \}^{-1} .
\end{align*}
\end{itemize}

\begin{table}[t!]
\caption{Summary of results with linear outcome models ($n=800$, $p=200$)} \label{linearOR-n800-p200}  \vspace{-.02in}
\footnotesize
\begin{center}
\begin{tabular*}{1\textwidth}{@{\extracolsep\fill} c c ll c ll c ll} \hline
      && \multicolumn{2}{c}{cor PS, cor OR} && \multicolumn{2}{c}{cor PS, mis OR} && \multicolumn{2}{c}{mis PS, cor OR} \\ \cline{3-4}  \cline{6-7} \cline{9-10}
      && RML.RML & RCAL.RWL  && RML.RML  &  RCAL.RWL  && RML.RML   &  RCAL.RWL   \\ \hline
      && \multicolumn{8}{c}{Linear outcome configuration 1} \\
Bias  && $-.041$ & $-.022$ && $-.007$ & $-.008$ && $-.006$ & $-.002$ \\
$\sqrt{\mbox{Var}}$
      && .071  & .071 && .072 & .072  && .077 & .072$^\dag$ \\   %(77/72)^2 = 1.14
$\sqrt{\mbox{EVar}}$
      && .083  & .083 && .081 & .080  && .083 & .083 \\
Cov90 && .790  & .822$^*$ && .850 & .848  && .856 & .837 \\
Cov95 && .859  & .891$^*$ && .910 & .915  && .925 & .912  \\
      && \multicolumn{8}{c}{Linear outcome configuration 2} \\
Bias  && $-.038$ & $-.019$ && $-.040$ & $-.019$ && $-.006$ & $-.002$ \\
$\sqrt{\mbox{Var}}$
      && .063  & .063 && .062 & .064 && .069 & .063$^\dag$ \\  % (69/63)^2 = 1.20
$\sqrt{\mbox{EVar}}$
      && .072  & .073 && .070 & .069 && .074 & .074 \\
Cov90 && .782  & .826$^*$ && .786 & .865$^*$ && .855 & .838 \\
Cov95 && .858  & .885 && .866 & .918$^*$ && .926 & .901 \\ \hline
\end{tabular*}\\[.1in]
\parbox{1\textwidth}{\small Note: RML.RML denotes $\hat \mu^1(\hat m^1_{\mbox{\tiny RML}}, \hat\pi^1_{\mbox{\tiny RML}})$
and RCAL.RWL denotes $\hat \mu^1(\hat m^1_{\mbox{\tiny RWL}}, \hat\pi^1_{\mbox{\tiny RCAL}})$. Bias and
$\sqrt{\mbox{Var}}$ are respectively the Monte Carlo bias and standard deviation of the points estimates,
$\sqrt{\mbox{EVar}}$ is the square root of the mean of the variance estimates, and Cov90 or Cov95 is the coverage proportion of
the 90\% or 95\% confidence intervals, based on 1000 repeated simulations. $^\dag$ indicates a case where the Monte Carlo variance from the competitive method is
at least 10\% higher. $^*$ indicates a coverage proportion that is 3\% or higher than that from
the competitive method.}
\end{center}  \vspace{-.1in}
\end{table}

As in Section~\ref{sec:setup},
the observed data consist of independent and identically distributed observations $\{(T_i Y_i, T_i, X_i): i=1,\ldots, n\}$.
Consider logistic
propensity score model (\ref{model-PS}) and linear outcome model (\ref{lm-OR}), both with $f_j(X) = X_j^\dag$ for $j=1,\ldots,p$.
Then the two models can be classified as follows, depending on the data configuration above:
\begin{itemize} \addtolength{\itemsep}{-.05in}
\item[(C1)] PS and OR models both correctly specified;

\item[(C2)] PS model correctly specified, but OR model misspecified;

\item[(C3)] PS model misspecified, but OR model correctly specified.
\end{itemize}
%See the Supplementary Material for boxplots of $X^\dag_j$ within $\{T=1\}$ and $\{T=0\}$
%and scatterplots of $Y$ against $X^\dag_j$ within $\{T=1\}$.
As demonstrated in Kang \& Schafer (2007) for $p=4$, the PS model (\ref{model-PS}) in the scenario (C3), although
misspecified, appears adequate as examined by conventional techniques for logistic regression.
Similarly, the OR model (\ref{lm-OR}) in the misspecified case (C2) can also
be shown as ``nearly correct" by standard techniques for linear regression.
On the other hand, neither the PS model (\ref{model-PS}) in the correctly specified case (C1) or (C2)
nor the OR model (\ref{lm-OR}) in the correctly specified case (C1) or (C3) is used in Kang \& Schafer (2007),
where correct PS and misspecified OR model (or misspecified PS and correct OR) involve two completely different sets of regressors.
This aspect of the Kang--Schafer design needs to be modified in our study,
where the same vector of regressors $f(X)$ is used in models (\ref{model-PS}) and (\ref{lm-OR}).

\begin{figure}[t!]
\caption{\small QQ plots of the $t$-statistics against standard normal with linear outcome models ($n=800$, $p=200$),
based on the estimators $\hat \mu^1(\hat m^1_{\mbox{\tiny RML}}, \hat\pi^1_{\mbox{\tiny RML}})$ ($\circ$)
and $\hat \mu^1(\hat m^1_{\mbox{\tiny RWL}}, \hat\pi^1_{\mbox{\tiny RCAL}})$ ($\times$).
For readability, only a subset of 100 order statistics are shown as points on the QQ lines.}
\label{fig:linearOR-n800-p200} \vspace{.15in}
\begin{tabular}{c}
\includegraphics[width=6.2in, height=4.5in]{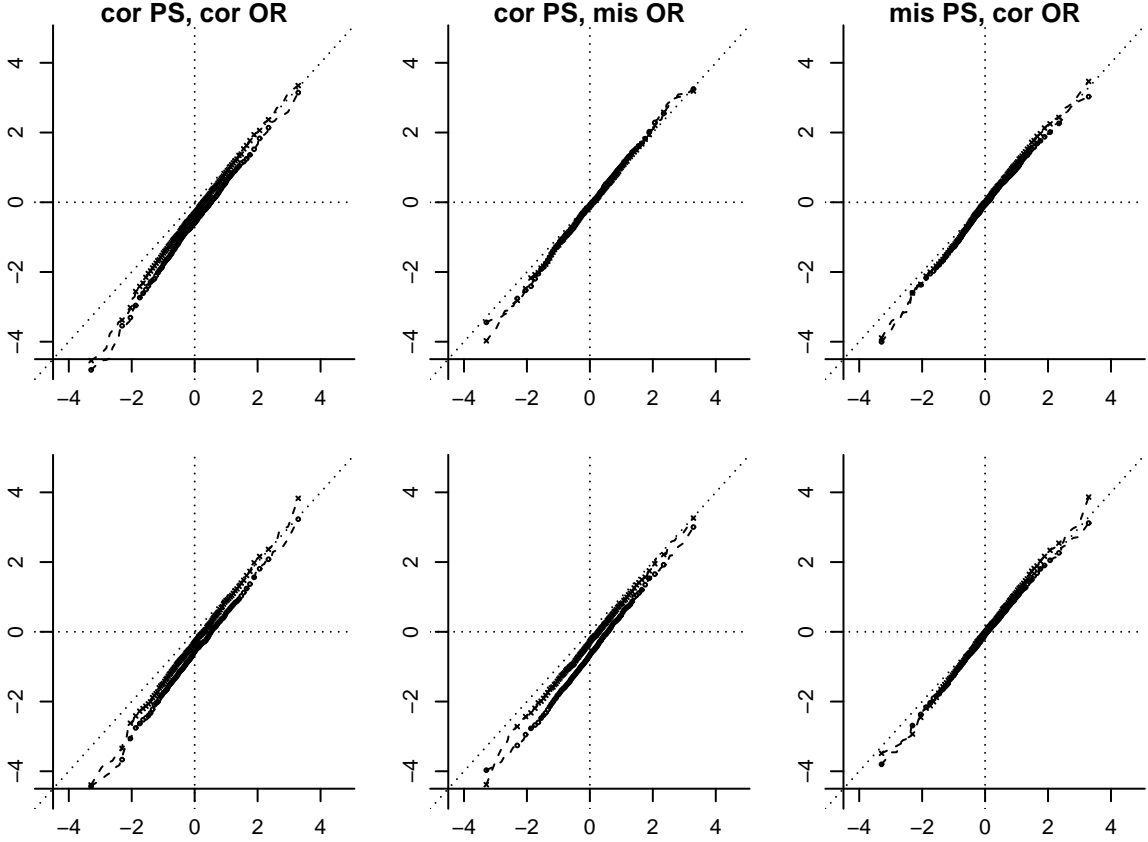} \vspace{-.25in}
\end{tabular}
\end{figure}

We conducted 1000 repeated simulations, each with the sample size $n=400$ or $800$ and the number of regressors $p=100$ or $200$.
For $n=800$ and $p=200$, Table~\ref{linearOR-n800-p200} summarizes the results about $\hat \mu^1(\hat m^1_{\mbox{\tiny RML}}, \hat\pi^1_{\mbox{\tiny RML}})$ and
$\hat \mu^1(\hat m^1_{\mbox{\tiny RWL}}, \hat\pi^1_{\mbox{\tiny RCAL}})$ and their associated confidence intervals,
and Figure~\ref{fig:linearOR-n800-p200} presents the QQ plots of the corresponding $t$-statistics.
See the Supplementary Material for similar results obtained with other values of $(n,p)$.

There are several advantages demonstrated from these results for the proposed method.
Compared with $\hat \mu^1(\hat m^1_{\mbox{\tiny RML}}, \hat\pi^1_{\mbox{\tiny RML}})$,
the estimator $\hat \mu^1(\hat m^1_{\mbox{\tiny RWL}}, \hat\pi^1_{\mbox{\tiny RCAL}})$ has consistently smaller biases in absolute values,
and either similar or noticeably smaller variances, for example, in the case of misspecified PS model and correct OR model.
The coverage proportions of confidence intervals based on $\hat \mu^1(\hat m^1_{\mbox{\tiny RWL}}, \hat\pi^1_{\mbox{\tiny RCAL}})$  are similar or noticeably higher than
those based on $\hat \mu^1(\hat m^1_{\mbox{\tiny RML}}, \hat\pi^1_{\mbox{\tiny RML}})$, although both coverage proportions are below the nominal probabilities to various degree.
From the QQ plots, the $t$-statistics based on $\hat \mu^1(\hat m^1_{\mbox{\tiny RWL}}, \hat\pi^1_{\mbox{\tiny RCAL}})$ also appear to be more
aligned with standard normal than those based on $\hat \mu^1(\hat m^1_{\mbox{\tiny RML}}, \hat\pi^1_{\mbox{\tiny RML}})$.

\begin{table}[t!]
\caption{Summary of results with logistic outcome models ($n=800$, $p=200$)} \label{logitOR-n800-p200}  \vspace{-.02in}
\footnotesize
\begin{center}
\begin{tabular*}{1\textwidth}{@{\extracolsep\fill} c c ll c ll c ll} \hline
      && \multicolumn{2}{c}{cor PS, cor OR} && \multicolumn{2}{c}{cor PS, mis OR} && \multicolumn{2}{c}{mis PS, cor OR} \\ \cline{3-4}  \cline{6-7} \cline{9-10}
      && RML.RML & RCAL.RWL  && RML.RML  &  RCAL.RWL  && RML.RML   &  RCAL.RWL   \\ \hline
      && \multicolumn{8}{c}{Logistic outcome configuration 1} \\
Bias  && $-.013$ & $-.004$ && $-.007$ & $-.003$ && $-.005$ & $-.001$ \\
$\sqrt{\mbox{Var}}$
      && .023  & .024 && .023 & .024  && .024 & .023 \\
$\sqrt{\mbox{EVar}}$
      && .026  & .026 && .025 & .026  && .027 & .027  \\
Cov90 && .814  & .868$^*$ && .841 & .872$^*$  && .845 & .859 \\
Cov95 && .876  & .920$^*$ && .916 & .928  && .914 & .912  \\
      && \multicolumn{8}{c}{Logistic outcome configuration 2} \\
Bias  && $-.009$ & $-.003$ && $-.007$ & $-.002$ && $-.002$ & $.001$ \\
$\sqrt{\mbox{Var}}$
      && .024  & .025 && .024 & .026 && .026 & .025  \\
$\sqrt{\mbox{EVar}}$
      && .026  & .026 && .026 & .027 && .027 & .027 \\
Cov90 && .849  & .876 && .864 & .879 && .870 & .865  \\
Cov95 && .909  & .936 && .925 & .933 && .931 & .927  \\ \hline
\end{tabular*}\\[.1in]
\parbox{1\textwidth}{\small Note: See the footnote of Table~\ref{linearOR-n800-p200}. For scenario (C5) (``mis OR"), the true value $\mu^1$ is $0.5$ by symmetry. For scenarios (C4) and (C6)
(``cor OR"), the true value $\mu^1$ is not analytically available but calculated using Monte Carlo integration,
as shown in the Supplementary Material.}
\end{center}  \vspace{-.1in}
\end{table}

\subsection{Logistic outcome models} \label{sec:sim-logitOR}

\begin{figure}[t!]
\caption{\small QQ plots of the $t$-statistics against standard normal with logistic outcome models ($n=800$, $p=200$),
based on the estimators $\hat \mu^1(\hat m^1_{\mbox{\tiny RML}}, \hat\pi^1_{\mbox{\tiny RML}})$ ($\circ$)
and $\hat \mu^1(\hat m^1_{\mbox{\tiny RWL}}, \hat\pi^1_{\mbox{\tiny RCAL}})$ ($\times$).
For readability, only a subset of 100 order statistics are shown as points on the QQ lines.}
\label{fig:logitOR-n800-p200} \vspace{.15in}
\begin{tabular}{c}
\includegraphics[width=6.2in, height=4.5in]{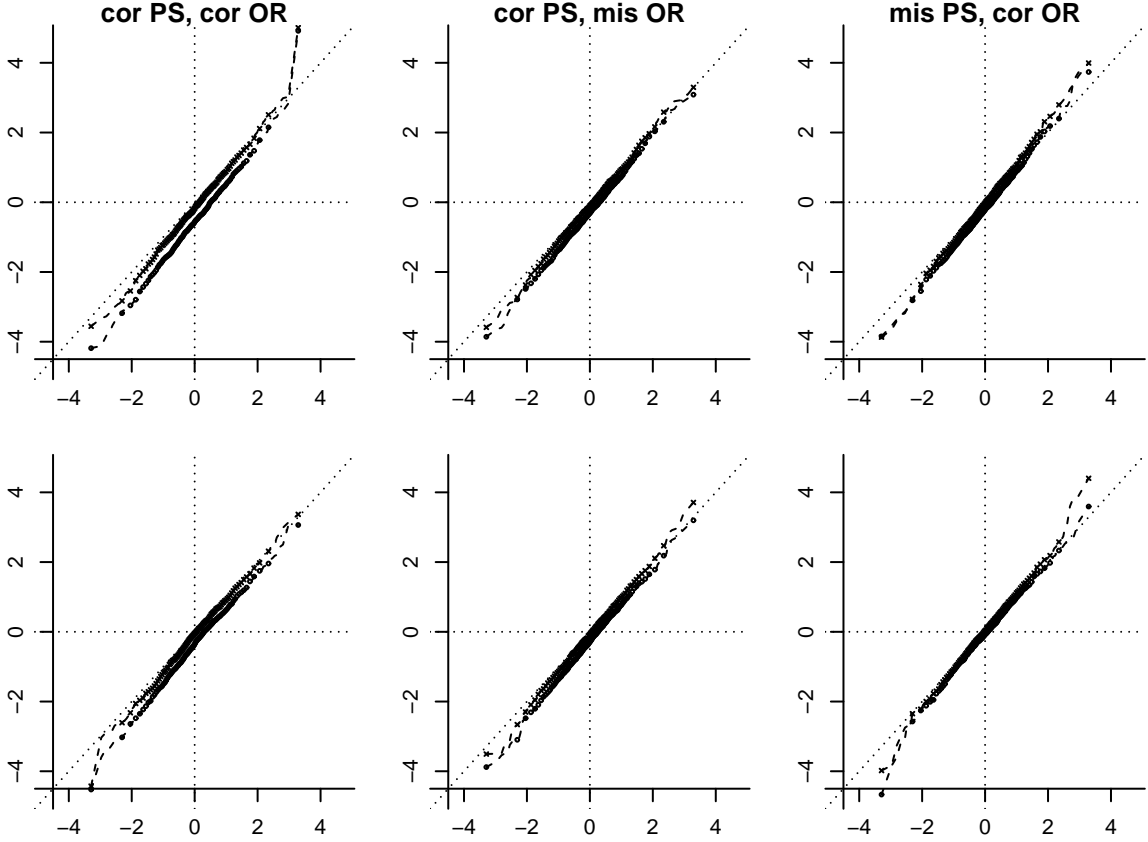} \vspace{-.25in}
\end{tabular}
\end{figure}

For simulations with binary outcomes, let $X$ and $X^\dag$ be as in Section~\ref{sec:sim-linearOR}.
Consider the following data-generating configurations, in parallel to (C1)--(C3).
\begin{itemize} \addtolength{\itemsep}{-.05in}

\item[(C4)] Generate $T$ given $X$ as in (C1) and, independently, generate $Y^1$ given $X$ from a Bernoulli distribution with probability (``Logistic outcome configuration 1") \vspace{-.1in}
\begin{align*}
P( Y^1=1 | X) = [ 1+ \exp\{ -(X^\dag_1 + 0.5 X^\dag_2 + 0.5 X^\dag_3 + 0.5 X^\dag_4)\}]^{-1} ,
\end{align*}
or (``Logistic outcome configuration 2") \vspace{-.1in}
\begin{align*}
P( Y^1=1 | X) = [1+ \exp\{-(0.25 X^\dag_1 + 0.5 X^\dag_2 + 0.5 X^\dag_3 + 0.5 X^\dag_4)\}]^{-1} .
\end{align*}

\item[(C5)] Generate $T$ give $X$ as in (C1), and, independently, generate $Y^1$ given $X$ from a Bernoulli distribution with probability  (``Logistic outcome configuration 1") \vspace{-.1in}
\begin{align*}
P( Y^1=1 | X) = [ 1+ \exp\{ -(X_1 + 0.5 X_2 + 0.5 X_3 + 0.5 X_4)\}]^{-1} ,
\end{align*}
or (``Logistic outcome configuration 2") \vspace{-.1in}
\begin{align*}
P( Y^1=1 | X) = [ 1+ \exp\{ -(0.25 X_1 + 0.5 X_2 + 0.5 X_3 + 0.5 X_4)\}]^{-1} .
\end{align*}

\item[(C6)] Generate  $Y^1$ given $X$ as in (C4), and, independently, generate $T$ given $X$  as in (C3).
\end{itemize}
Consider logistic
propensity score model (\ref{model-PS}) and logistic outcome model (\ref{glm-OR}), both with $f_j(X) = X_j^\dag$ for $j=1,\ldots,p$.
Then the two models are correctly specified in scenario (C4), only PS model (\ref{model-PS}) is correctly specified in scenario (C5),
and only OR model (\ref{glm-OR}) is correctly specified in scenario (C6), similarly as in Section~\ref{sec:sim-linearOR}.

For $n=800$ and $p=200$, Table~\ref{logitOR-n800-p200} and Figure~\ref{fig:logitOR-n800-p200} present the results from 1000 repeated simulations, about
$\hat \mu^1(\hat m^1_{\mbox{\tiny RML}}, \hat\pi^1_{\mbox{\tiny RML}})$ and
$\hat \mu^1(\hat m^1_{\mbox{\tiny RWL}}, \hat\pi^1_{\mbox{\tiny RCAL}})$ and their associated confidence intervals. Similar conclusions can be drawn as from Table~\ref{linearOR-n800-p200} and
Figure~\ref{fig:linearOR-n800-p200}. It is interesting that the coverage proportions of confidence intervals based on $\hat \mu^1(\hat m^1_{\mbox{\tiny RWL}}, \hat\pi^1_{\mbox{\tiny RCAL}})$
are noticeably higher (and closer to the nominal probabilities)
than those based on $\hat \mu^1(\hat m^1_{\mbox{\tiny RML}}, \hat\pi^1_{\mbox{\tiny RML}})$ in the case where both PS and OR models are correctly specified.
This difference can also be seen from the QQ plots.
The confidence intervals from both methods appear to yield reasonable coverage proportions when the PS model is misspecified
but the OR model is correctly specified, even though these results are not necessarily predicted by asymptotic theory.
See the Supplementary Material for additional results from simulations with other values of $(n,p)$.

\section{Application to a medical study} \label{sec:application}

We provide an empirical application to a medical study in Connors et al.~(1996) on the effects of right heart catheterization (RHC).
The study  included $n=5735$ critically ill patients  admitted  to  the  intensive  care  units  of  5  medical  centers.
For each patient, the data consist of treatment status $T$ ($=1$ if RHC was used
within  24  hours  of  admission  and  0  otherwise),  health  outcome $Y$
(survival  time  up  to  30  days),  and  a list of 75 covariates $X$
specified by medical specialists in critical care.
For previous analyses, propensity score and outcome regression models were employed
either with main effects only (Hirano \& Imbens 2002; Vermeulen \& Vansteelandt 2015) or
with interaction terms manually added (Tan 2006).

\begin{table}[b!]
\caption{Estimates of $30$-day survival probabilities and ATE} \label{rhc-ATE}  \vspace{-.02in}
\footnotesize
\begin{center}
\begin{tabular*}{.8\textwidth}{@{\extracolsep\fill} c cc c cc} \hline
         & \multicolumn{2}{c}{IPW} && \multicolumn{2}{c}{Augmented IPW} \\ \cline{2-3}  \cline{5-6}
         & RML  & RCAL   &&  RML.RML  & RCAL.RWL  \\ \hline
$\mu^1$  & $0.636 \pm 0.026$  & $0.634 \pm 0.023$  && $0.636 \pm 0.021$  & $0.635 \pm 0.021$ \\
$\mu^0$  & $0.690 \pm 0.017$  & $0.687 \pm 0.017$  && $0.691 \pm 0.016$  & $0.688 \pm 0.016$  \\
ATE      & $-0.054\pm 0.031$  & $-0.053\pm 0.029$  && $-0.055 \pm 0.025$ & $-0.053 \pm 0.025$ \\ \hline
\end{tabular*}\\[.1in]
\parbox{.8\textwidth}{\small Note: Estimate $\pm \; 2 \times$standard error, including nominal standard errors for IPW.}
\end{center}  \vspace{-.1in}
\end{table}

To explore dependency beyond main effects, we consider a logistic propensity score model (\ref{model-PS})
and a logistic outcome model (\ref{glm-OR}) for 30-day survival status $1\{Y>30\}$,
with the vector $f(X)$ including all main effects and two-way interactions of $X$ except those with the fractions of nonzero values less than 46 (i.e., 0.8\% of the sample size 5735).
The dimension of $f(X)$ is $p=1855$, excluding the constant. All variables in $f(X)$ are standardized with sample means 0 and variances 1.
We apply the augmented IPW estimators $\hat \mu^1(\hat m^1_{\mbox{\tiny RWL}}, \hat\pi^1_{\mbox{\tiny RCAL}})$ and $\hat \mu^0(\hat m^0_{\mbox{\tiny RWL}}, \hat\pi^0_{\mbox{\tiny RCAL}})$
using regularized calibrated (RCAL) estimation and the corresponding estimators such as $\hat \mu^1(\hat m^1_{\mbox{\tiny RML}}, \hat\pi_{\mbox{\tiny RML}})$
using regularized maximum likelihood (RML) estimation, similarly as in the simulation study.
The Lasso tuning parameter $\lambda$ is selected by cross validation over a discrete set $\{\lambda^* / 2^{j/4}: j=0,1,\ldots,24\}$,
where $\lambda^*$ is the value leading to a zero solution $\gamma_1=\cdots=\gamma_p=0$.
We also compute the (ratio) IPW estimators, such as $\hat \mu^1_{\mbox{\tiny rIPW}}$, along with nominal standard errors
obtained by ignoring data-dependency of the fitted propensity scores.

\begin{figure}[t!]
\caption{\small Boxplots of inverse probability weights within the treated (left) and untreated (middle) groups, each normalized to sum to the sample size $n$,
and QQ plots with a 45-degree line of the standardized sample influence functions based on $\varphi(Y,T,X;\cdot)$ in (\ref{influence-function}) for ATE (right).}
\label{fig:rhc-box-qq} \vspace{.15in}
\begin{tabular}{c}
\includegraphics[width=6.2in, height=2.5in]{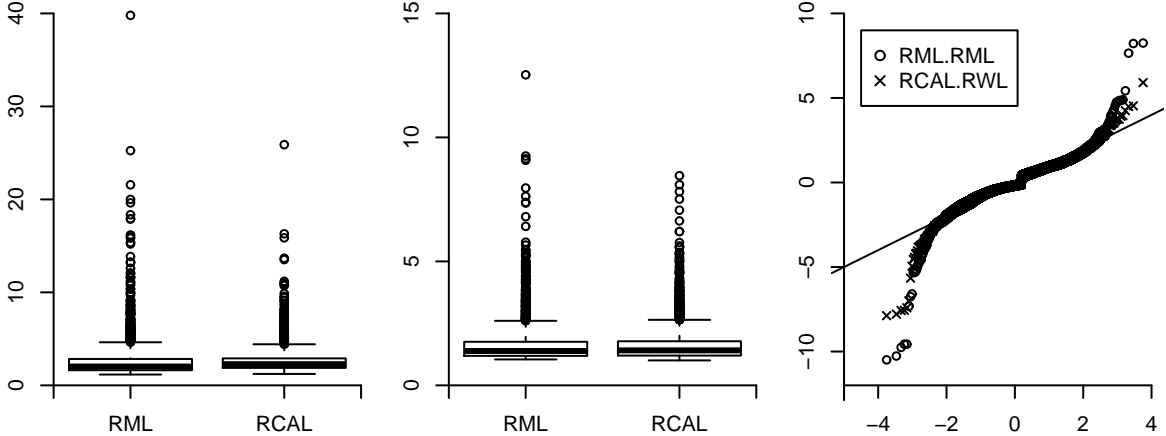} \vspace{-.25in}
\end{tabular}
\end{figure}

Table~\ref{rhc-ATE} shows various estimates of survival probabilities and ATE.\ The IPW estimates from RCAL estimation of propensity scores
have noticeably smaller nominal standard errors than RML estimation,
for example, with the relative efficiency $(0.026/0.023)^2=1.28$ for estimation of $\mu^1$.
This improvement can also be seen from Figure~\ref{fig:rhc-box-qq}, where the RCAL inverse probability weights are much less variable than RML weights.
See Tan (2017) for additional results on covariate balance and parameter sparsity from RML and RCAL estimation of propensity scores.

The augmented IPW estimates and confidence intervals are similar to each other from RCAL and RML estimation.
However, the validity of RML confidence intervals depends on both PS and OR models being correctly specified,
whereas that of RCAL confidence intervals holds even when the OR model is misspecified.
While assessment of this difference is difficult with real data, Figure~\ref{fig:rhc-box-qq} shows that
the sample influence functions for ATE using RCAL estimation appears to be more normally distributed especially in the tails
than RML estimation.

Finally, the augmented IPW estimates here are smaller in absolute values, and also with smaller standard errors, than previous estimates
based on main-effect models, about $-0.060 \pm 2 \times 0.015$ (Vermeulen \& Vansteelandt 2015).
The reduction in standard errors might be explained by the well-known property that an augmented IPW estimator has a smaller asymptotic variance when obtained using
a larger (correct) propensity score model.

\vspace{.25in}
\centerline{\bf\Large References}

\begin{description}\addtolength{\itemsep}{-.12in}

\item Athey, S., Imbens, G.W., and Wager, S. (2016) ``Approximate residual balancing: De-biased inference of average treatment effects in high dimensions," arXiv:1604.07125.

\item Belloni, A., Chernozhukov, V., Fernandez-Val, I., and Hansen, C. (2017) "Program evaluation and causal inference
with high-dimensional data," {\em Econometrica}, 85, 233--298.

\item Bickel, P., Ritov, Y., and Tsybakov, A.B. (2009) ``Simultaneous analysis of Lasso and Dantzig selector," {\em Annals
of Statistics}, 37, 1705--1732.

\item Buhlmann, P. and van de Geer, S. (2011) {\em Statistics for High-Dimensional Data: Methods, Theory and Applications}, New York: Springer.

\item Chan, K.C.G., Yam, S.C.P., and Zhang, Z. (2016)
``Globally efficient non-parametric inference of average treatment effects by empirical balancing
calibration weighting," {\em Journal of the Royal Statistical Society}, Ser. B,  78, 673--700.

\item Connors, A.F., Speroff, T., Dawson, N.V., et al. (1996) ``The effectiveness
of right heart catheterization in the initial care of critically ill patients,"
{\em Journal of the American Medical Association}, 276, 889--897.

\item Farrell, M.H. (2015) ``Robust inference on average treatment effects with possibly more covariates
than observations." {\em Journal of Econometrics}, 189, 1--23.

\item Folsom, R.E. (1991) ``Exponential and logistic weight adjustments for sampling and nonresponse error reduction,"
{\em Proceedings of the American Statistical Association}, Social Statistics Section, 197--202.

%\item Friedman, J., Hastie, T., Hoefling, H. and Tibshirani, R. (2007) ``Pathwise
%coordinate optimization," {\em Annals of Applied Statistics}, 1, 302--332.

\item Friedman, J., Hastie, T., and Tibshirani, R. (2010) ``Regularization paths for generalized linear models via coordinate descent,"
{\em Journal of Statistical Software}, 33, 1--22.

\item Graham, B.S., de Xavier Pinto, C.C., and Egel, D. (2012) ``Inverse probability
tilting for moment condition models with missing data," {\em Review of Economic Studies}, 79, 1053--1079.

\item Hirano, K., and Imbens,  G.W. (2002) ``Estimation of causal effects using
propensity score weighting: An application to data on right heart catheterization,"
{\em Health Services and Outcomes Research Methodology}, 2, 259--278.

\item Hahn, J. (1998) ``On the role of the propensity score in efficient semiparametric
estimation of average treatment effects," {\em Econometrica}, 66, 315--331.

\item Hainmueller, J. (2012) ``Entropy balancing for causal effects: Multivariate reweighting method
to produce balanced samples in observational studies," {\em Political Analysis}, 20, 25--46.

%\item Hastie, T., Tibshirani, R., and Wainwright, M.J. (2015) {\em Statistical Learning with Sparsity: The Lasso and Generalizations}, New York: Chapman \& Hall.

\item Huang, J. and Zhang, C.-H. (2012) ``Estimation and selection via absolute penalized convex minimization
and its multistage adaptive applications," {\em Journal of Machine Learning Research}, 13, 1839--1864.

\item Imai, K. and Ratkovic, M. (2014) ``Covariate balancing propensity score," {\em Journal of the Royal
Statistical Society}, Ser. B, 76, 243--263.

\item Javanmard, A. and Montanari, A. (2014) ``Confidence intervals and hypothesis testing
for high-dimensional regression," {\em Journal of Machine Learning Research}, 15, 2869--2909.

\item Kang, J.D.Y. and Schafer, J.L. (2007) ``Demystifying double robustness: A comparison of alternative strategies for estimating a population mean from incomplete data" (with discussion), {\em Statistical Science}, 523--539.

\item Kim, J.K. and Haziza, D. (2014) ``Doubly robust inference with missing data in survey sampling," {\em Statistica Sinica}, 24, 375--394.

\item Manski, C.F. (1988) {\em Analog Estimation Methods in Econometrics}, New York:
Chapman \& Hall.

\item McCullagh, P. and Nelder, J. (1989) {\em Generalized Linear Models} (2nd edition), New York: Chapman \& Hall.

\item Negahban, S.N., Ravikumar, P., Wainwright, M.J., and Yu, B. (2012) ``A unified
framework for high-dimensional analysis of M-estimators with decomposable
regularizers," {\em Statistical Science}, 27, 538--557.

\item Neyman, J. (1923) ``On the application of probability theory to agricultural
experiments: Essay on principles, Section 9," translated in {\em Statistical Science}, 1990, 5, 465--480.

\item Osborne, M., Presnell, B., and Turlach, B. (2000) ``A new approach to variable selection in least
squares problems." {\em IMA Journal of Numerical Analysis}, 20, 389--404.

\item Robins, J.M., Rotnitzky, A., and Zhao, L.P. (1994) ``Estimation of regression
coefficients when some regressors are not always observed," {\em Journal of the American
Statistical Association}, 89, 846--866.

\item Rosenbaum, P.R. and Rubin, D.B. (1983) ``The central role of the propensity score in observational studies for causal
effects," {\em Biometrika}, 70, 41-–55.

\item Rosenbaum, P.R. and Rubin, D.B. (1984) ``Reducing bias in observational studies using subclassification on the
propensity score," {\em Journal of the American Statistical Association}, 79, 516--524.

\item Rubin, D.B. (1976) ``Inference and missing data," {\em Biometrika}, 63, 581--590.

\item Sarndal, C.E., Swensson, B. and Wretman, J.H. (1992) {\em Model Assisted Survey Sampling}, New York: Springer.

%\item Sarndal, C.E. and Wright, R.L. (1984) `Cosmetic form of estimators in survey sampling," {\em Scandinavian Journal of Statistics}, 11, 146--156.

\item Tan, Z. (2006) ``A distributional approach for causal inference using propensity scores,"
{\em Journal of the American Statistical Association}, 101, 1619–-1637.

\item Tan, Z. (2007) ``Comment: Understanding OR, PS, and DR," {\em Statistical
Science}, 22, 560--568.

\item Tan, Z. (2010) ``Bounded, efficient, and doubly robust estimation with inverse weighting," {\em Biometrika}, 97, 661--682.

\item Tan, Z. (2017) ``Regularized calibrated estimation of propensity scores with model misspecification and high-dimensional data," arXiv:1710.08074.

\item Tibshirani, R. (1996) ``Regression shrinkage and selection via the Lasso," {\em Journal of the Royal Statistical Society}, Ser. B, 58, 267--288.

\item Tsiatis, A.A. (2006) {\em Semiparametric Theory and Missing Data}, New York: Springer.

\item van de Geer, S., Buhlmann, P., Ritov, Y., Dezeure, R. (2014) ``On asymptotically optimal
confidence regions and tests for high-dimensional models" {\em Annals of Statistics}, 42,
1166--1202.

\item Vermeulen. K. and Vansteelandt, S. (2015) ``Bias-reduced doubly
robust estimation," {\em Journal of the American Statistical Association}, 110, 1024--1036.

\item White, H. (1982) ``Maximum Likelihood Estimation of Misspecified Models," {\em Econometrica}, 50, 1--25.

\item Zhang, C.-H. and Zhang, S.S. (2014) ``Confidence intervals for low-dimensional parameters
with high-dimensional data," {\em Journal of the Royal Statistical Society}, Ser. B, 76, 217--242.
\end{description}

%\end{document}

\clearpage

\setcounter{page}{1}

\setcounter{section}{0}
\setcounter{equation}{0}

\setcounter{table}{0}
\setcounter{figure}{0}

\setcounter{pro}{0}
\renewcommand{\thepro}{S\arabic{pro}}

\setcounter{cor}{0}
\renewcommand{\thecor}{S\arabic{cor}}

\renewcommand{\theequation}{S\arabic{equation}}
\renewcommand{\thesection}{\Roman{section}}

\renewcommand\thefigure{S\arabic{figure}}
\renewcommand\thetable{S\arabic{table}}

\begin{center}
{\Large Supplementary Material for ``Model-assisted inference for treatment effects using regularized calibrated estimation with high-dimensional data"}

\vspace{.1in} Zhiqiang Tan
\end{center}
\vspace{.1in}

The Supplementary Material contains Appendices I--II.

\section{Additional results for simulation study}

\subsection{Results for simulation setup}

Denote by $\phi()$ the probability density function and $\Phi()$ the cumulative distribution function for $\N(0,1)$.
For $a=2.5$, let $Z$ be $\N(0,1)$ truncated to the interval $(-a,a)$, with the density function $\phi(z)/c$ if $z \in (-a,a)$ or 0 otherwise,
where $c = \Phi(a)- \Phi(-a) = 2\Phi(a)-1$. Then $E(Z)= 0$ and $\var(Z) = 1 - 2a\phi(a)/c$, denoted as $b^2$.

Let $(X_1,\ldots,X_4) = (Z_1,\ldots, Z_4)/b$,
where $(Z_1,\ldots, Z_4)$ are independent variables, each from $\N(0,1)$ truncated to $(-a,a)$.
The variables $(X^\dag_1,\ldots, X^\dag_4)$ are determined by standardization from $(X_1,\ldots,X_4)$ using the following results.
\begin{itemize}
\item $E( \me^{0.5 X_1} ) = \exp(\frac{1}{8b^2}) \{ \Phi (a-\frac{1}{2b}) - \Phi(-a- \frac{1}{2b}) \} /c$,

$\var( \me^{0.5 X_1} ) =  \exp(\frac{1}{2b^2}) \{ \Phi (a-\frac{1}{b} ) - \Phi(-a- \frac{1}{b}) \} /c -  E^2 ( \me^{0.5 X} )$.

\item $E( \frac{X_2 }{1+\me^{X_1 }} ) =0$,

$ \var ( \frac{X_2 }{1+\me^{X_1 }} ) = \frac{1}{c} \int_{-a}^a \frac{1}{(1+\me^{z/b})^2} \phi(z)\,\dif z \approx (0.54257865)^2$ by numerical integration.

\item $E \{( \frac{X_1X_3}{25 } + 0.6)^3\} = 3 /25^2 *(.6) + (.6)^3$,

 $E \{( \frac{X_1X_3}{25 } + 0.6)^6\} = m_6^2 /25^6 + 15 * m_4^2 /25^4 *(.6)^2 + 15 /25^2 *(.6)^4 + (.6)^6$,

 where $m_4 = \frac{1}{b^4 c} \int_{-a}^a z^4 \phi(z) \,\dif z = \frac{1}{b^4 c} \{ (3/2)*(2\Phi(z)-1) - z(z^2+3)\phi(z) \} |_{-a}^a $ and

 $m_6 =  \frac{1}{b^6 c} \int_{-a}^a z^6 \phi(z) \,\dif z =\frac{1}{b^6 c} \{ (15/2)*(2\Phi(z)-1) - z(z^4+5z^2+15) \phi(z) \}  |_{-a}^a $.

\item $E \{( X_2 + X_4 + 20)^2 \} = 2 + 20^2$,

 $E \{( X_2 + X_4 + 20)^4 \} = (2m_4+6) + 6 *2 *20^2 + 20^4$.
\end{itemize}

For binary outcomes in scenarios (C4) and (C6), the true value $\mu^1=E\{ m_1^*(X)\}$  is estimated
by Monte Carlo integration, using 100 repeated samples of $(X_1,\ldots,X_4)$ each of size $10^7$.
The estimates of $\mu^1$ are $ 0.4949676$ and $0.4992349$ in ``logistic outcome configuration 1"
and ``logistic outcome configuration 2."
The standard errors are smaller than $8\times 10^{-6}$.

\subsection{Additional simulation results}

\begin{figure}[b!]
\caption{\small Scatter plots of $(X^\dag_1,X^\dag_2, X^\dag_3,X^\dag_4)$ from a sample of size $n=800$.}
\label{fig:scatterplot} \vspace{.15in}
\begin{tabular}{c}
\includegraphics[width=6.2in, height=6in]{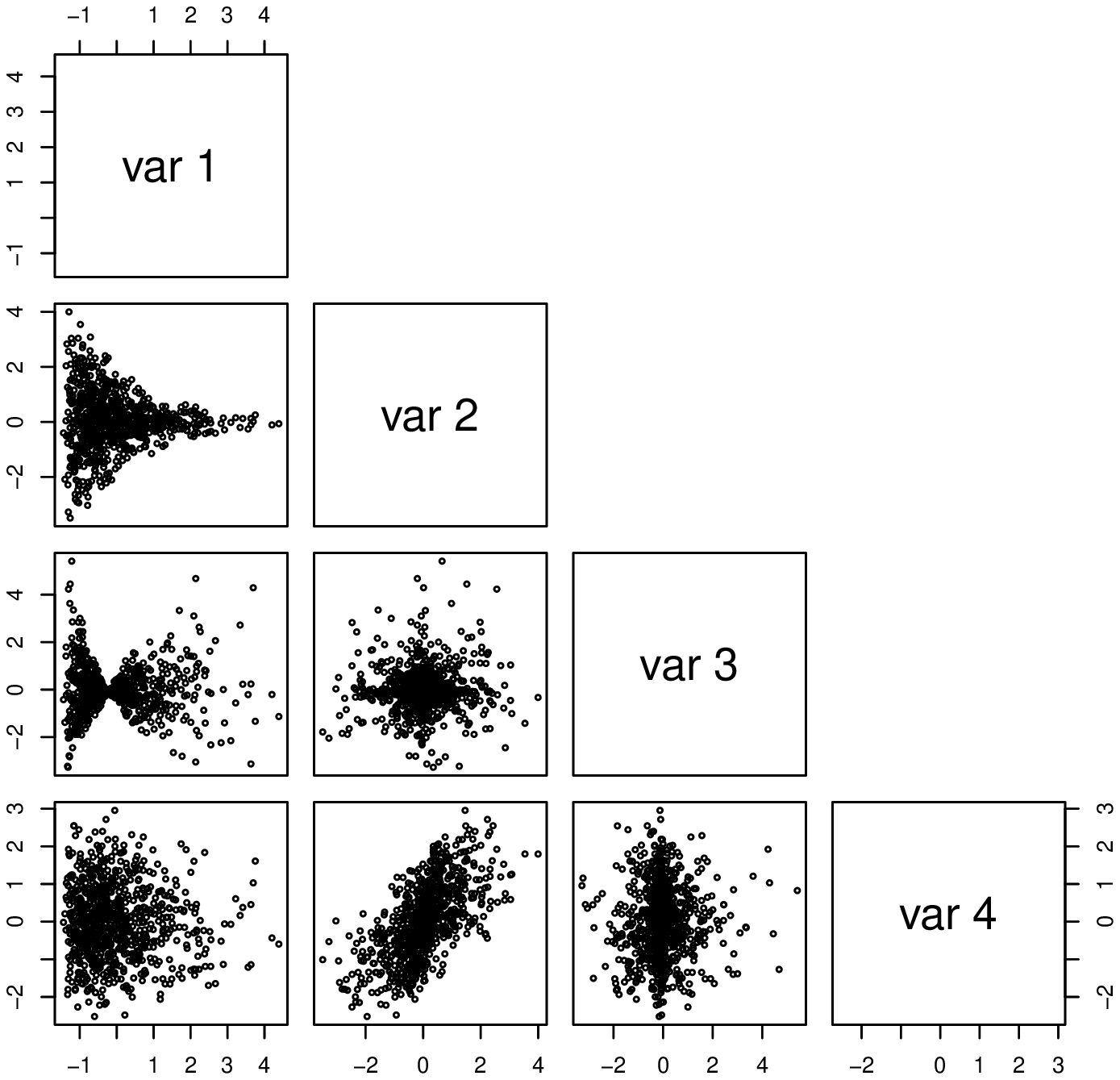} \vspace{-.25in}
\end{tabular}
\end{figure}

Figure~\ref{fig:scatterplot} shows the the scatter plots of the variables  $(X^\dag_1,X^\dag_2, X^\dag_3,X^\dag_4)$, which are correlated with each other as would be
found in real data.

\begin{figure}[t!]
\caption{\small QQ plots of the $t$-statistics against standard normal with linear outcome models ($n=400$, $p=100$),
based on the estimators $\hat \mu^1(\hat m^1_{\mbox{\tiny RML}}, \hat\pi^1_{\mbox{\tiny RML}})$ ($\circ$)
and $\hat \mu^1(\hat m^1_{\mbox{\tiny RWL}}, \hat\pi^1_{\mbox{\tiny RCAL}})$ ($\times$).
For readability, only a subset of 100 order statistics are shown as points on the QQ lines.}
\label{fig:linearOR-n400-p100} \vspace{.15in}
\begin{tabular}{c}
\includegraphics[width=6.2in, height=4.5in]{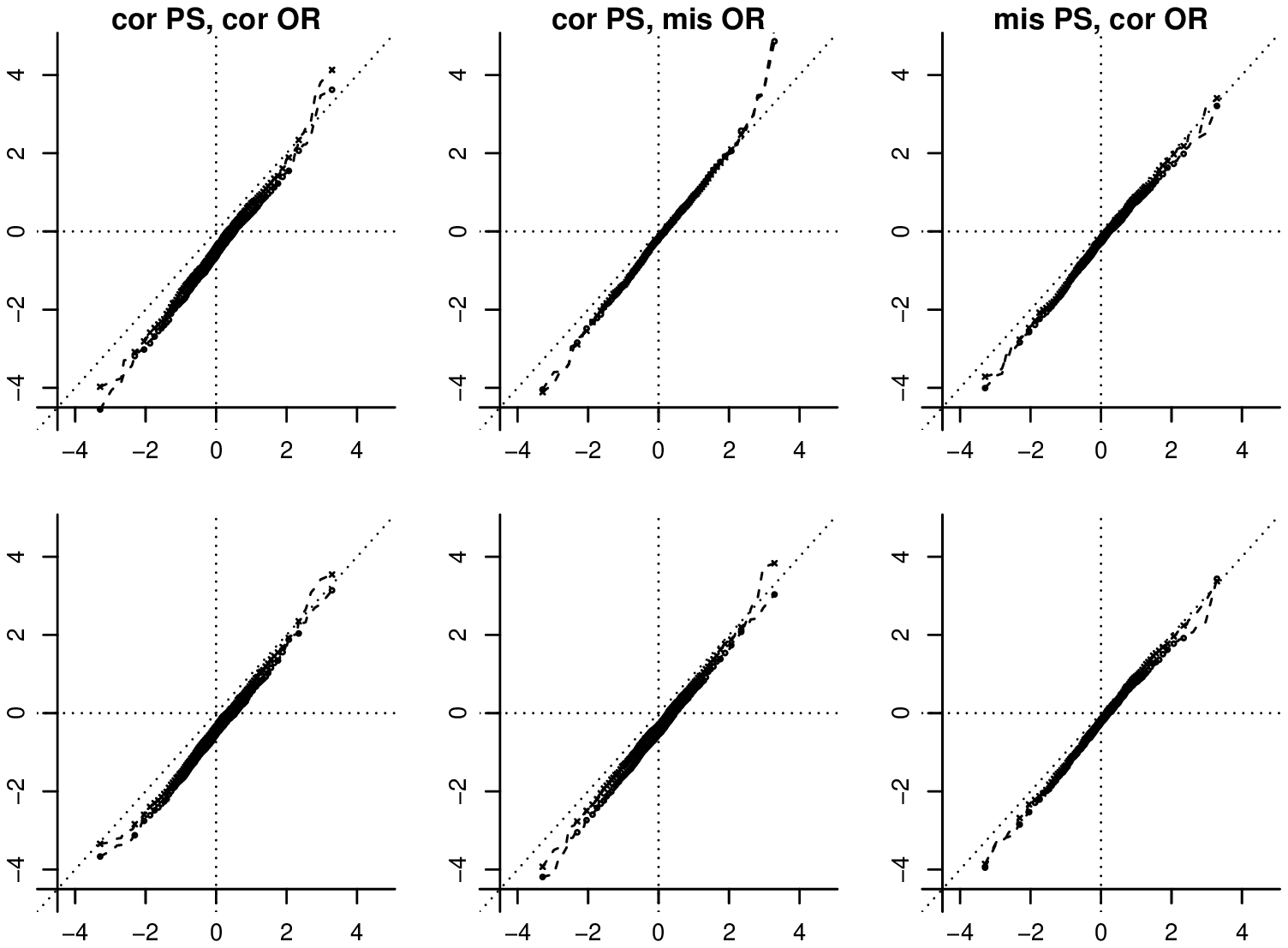} \vspace{-.25in}
\end{tabular}
\end{figure}

\begin{figure}[t!]
\caption{\small QQ plots of the $t$-statistics against standard normal with linear outcome models ($n=800$, $p=100$),
based on the estimators $\hat \mu^1(\hat m^1_{\mbox{\tiny RML}}, \hat\pi^1_{\mbox{\tiny RML}})$ ($\circ$)
and $\hat \mu^1(\hat m^1_{\mbox{\tiny RWL}}, \hat\pi^1_{\mbox{\tiny RCAL}})$ ($\times$).
For readability, only a subset of 100 order statistics are shown as points on the QQ lines.}
\label{fig:linearOR-n800-p100} \vspace{.15in}
\begin{tabular}{c}
\includegraphics[width=6.2in, height=4.5in]{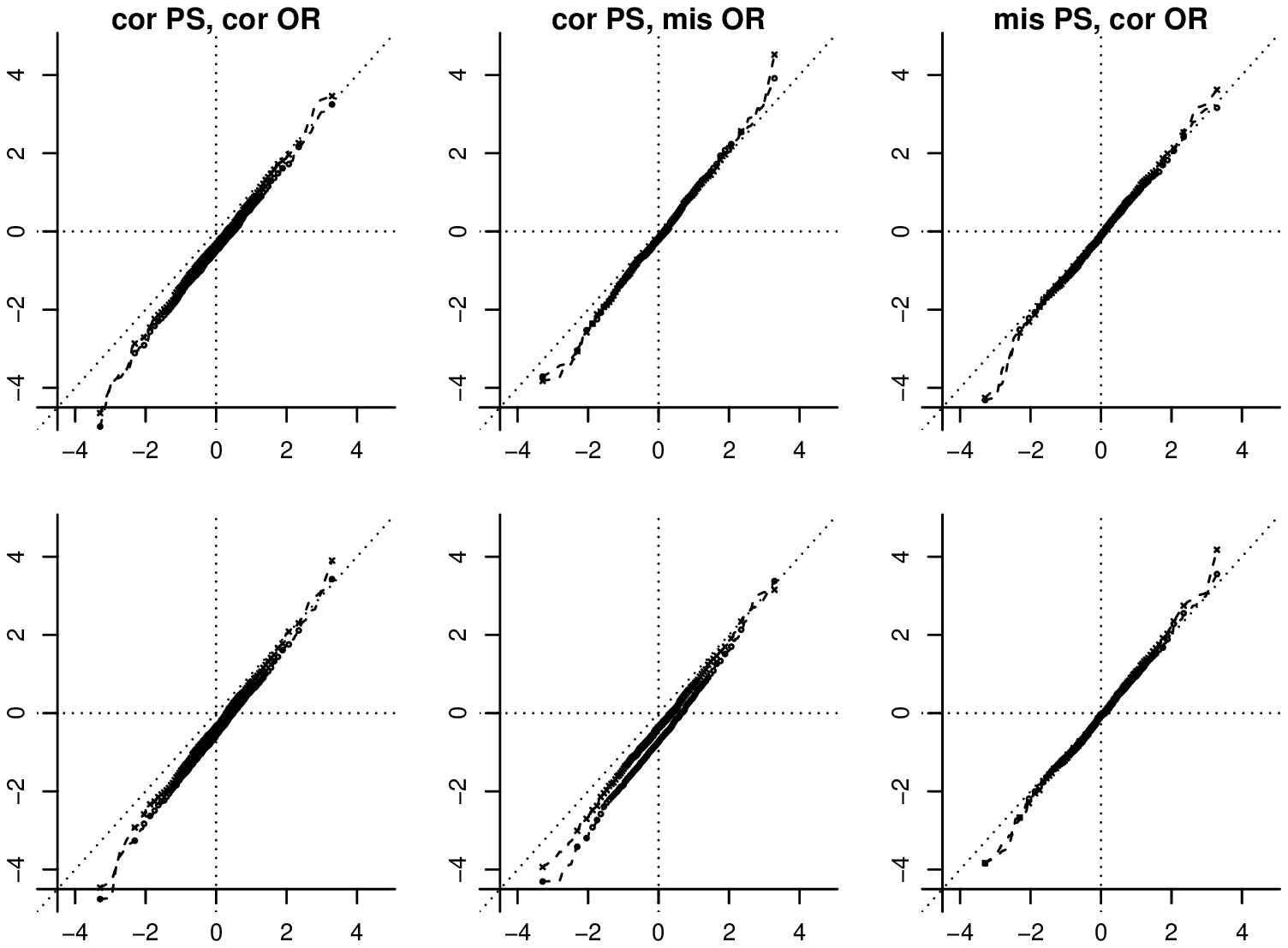} \vspace{-.25in}
\end{tabular}
\end{figure}

\begin{figure}[t!]
\caption{\small QQ plots of the $t$-statistics against standard normal with linear outcome models ($n=400$, $p=200$),
based on the estimators $\hat \mu^1(\hat m^1_{\mbox{\tiny RML}}, \hat\pi^1_{\mbox{\tiny RML}})$ ($\circ$)
and $\hat \mu^1(\hat m^1_{\mbox{\tiny RWL}}, \hat\pi^1_{\mbox{\tiny RCAL}})$ ($\times$).
For readability, only a subset of 100 order statistics are shown as points on the QQ lines.}
\label{fig:linearOR-n400-p200} \vspace{.15in}
\begin{tabular}{c}
\includegraphics[width=6.2in, height=4.5in]{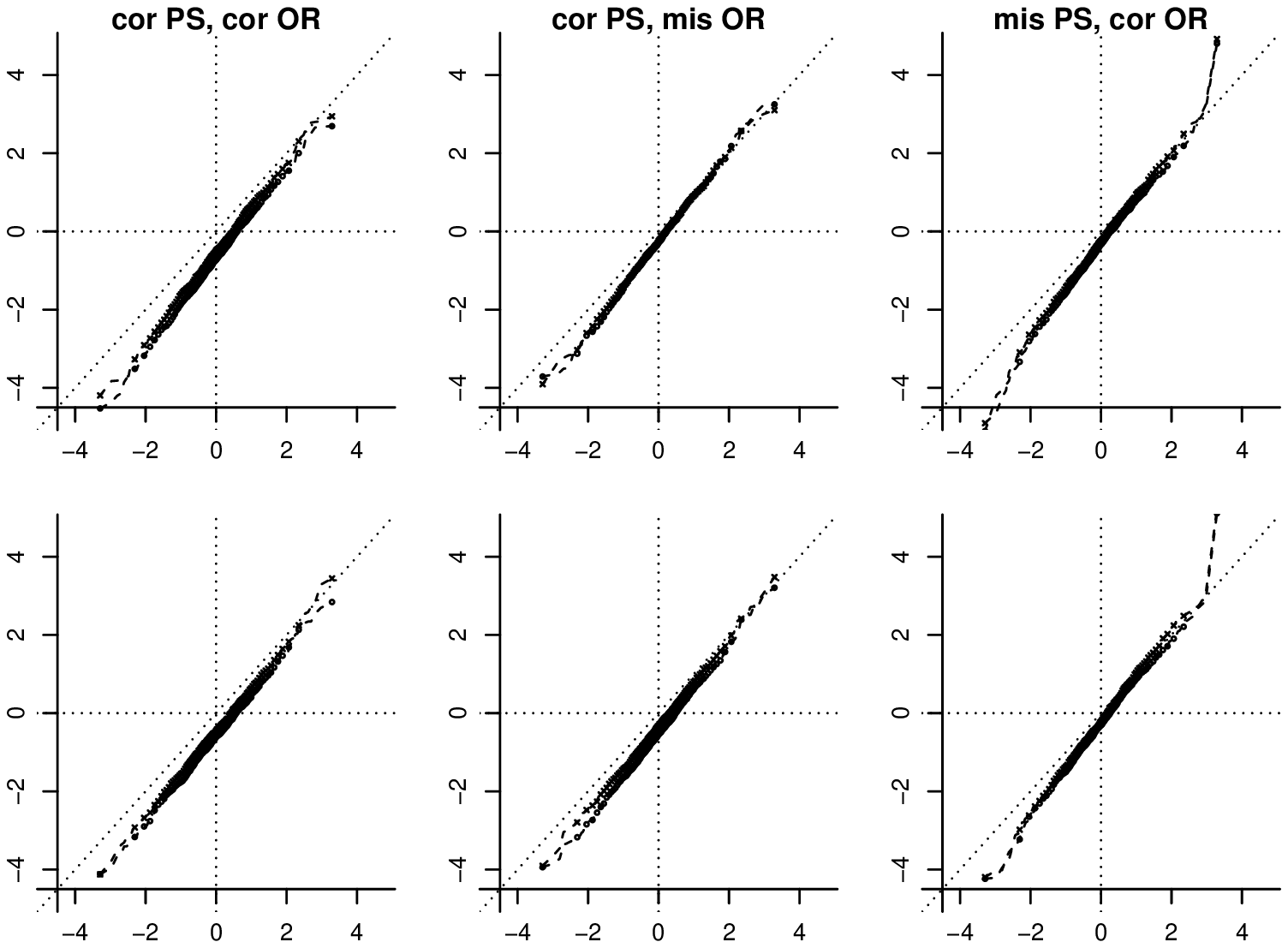} \vspace{-.25in}
\end{tabular}
\end{figure}

\begin{figure}[t!]
\caption{\small QQ plots of the $t$-statistics against standard normal with logistic outcome models ($n=400$, $p=100$),
based on the estimators $\hat \mu^1(\hat m^1_{\mbox{\tiny RML}}, \hat\pi^1_{\mbox{\tiny RML}})$ ($\circ$)
and $\hat \mu^1(\hat m^1_{\mbox{\tiny RWL}}, \hat\pi^1_{\mbox{\tiny RCAL}})$ ($\times$).
For readability, only a subset of 100 order statistics are shown as points on the QQ lines.}
\label{fig:logitOR-n400-p100} \vspace{.15in}
\begin{tabular}{c}
\includegraphics[width=6.2in, height=4.5in]{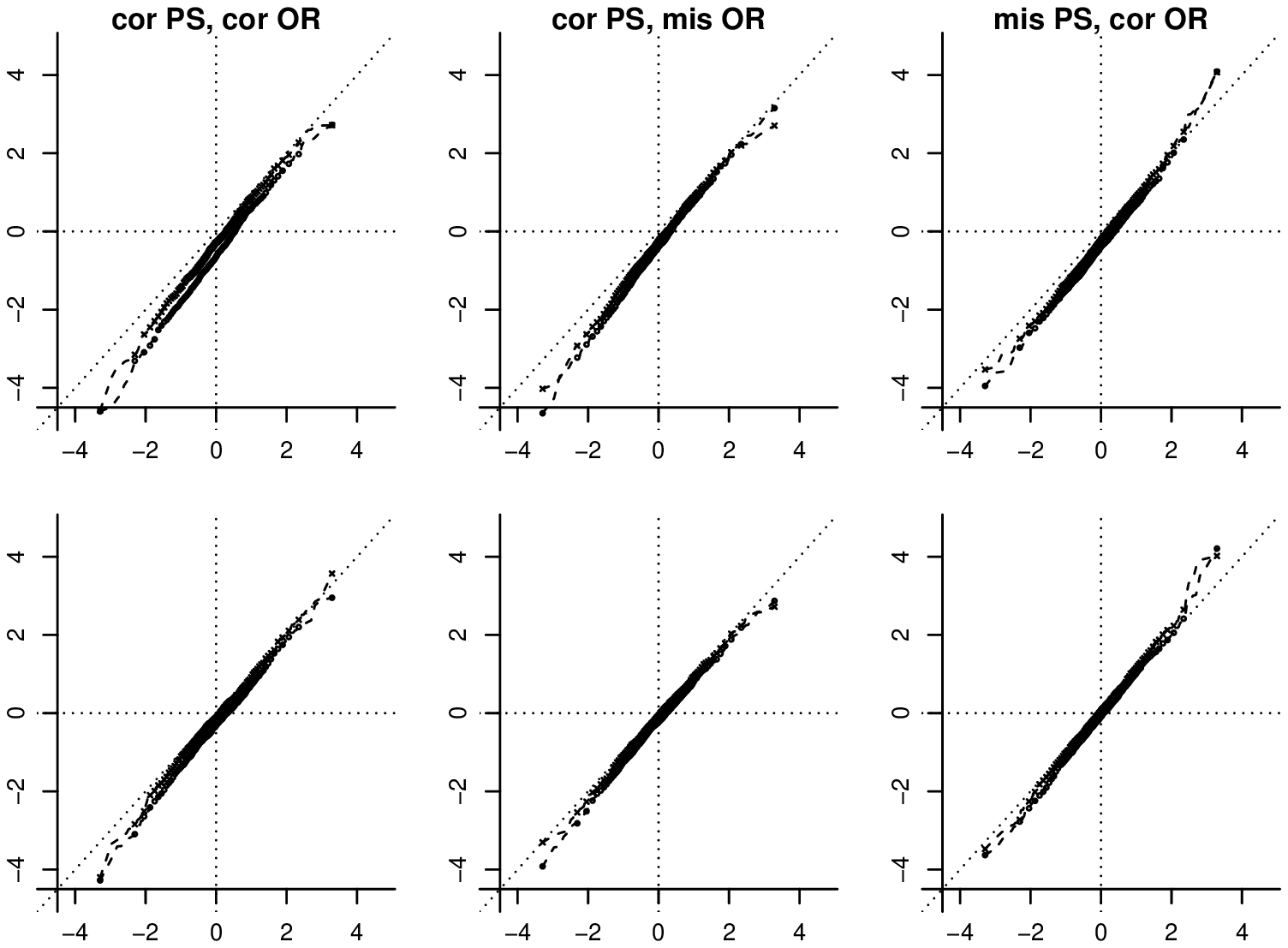} \vspace{-.25in}
\end{tabular}
\end{figure}

\begin{figure}[t!]
\caption{\small QQ plots of the $t$-statistics against standard normal with logistic outcome models ($n=800$, $p=100$),
based on the estimators $\hat \mu^1(\hat m^1_{\mbox{\tiny RML}}, \hat\pi^1_{\mbox{\tiny RML}})$ ($\circ$)
and $\hat \mu^1(\hat m^1_{\mbox{\tiny RWL}}, \hat\pi^1_{\mbox{\tiny RCAL}})$ ($\times$).
For readability, only a subset of 100 order statistics are shown as points on the QQ lines.}
\label{fig:logitOR-n800-p100} \vspace{.15in}
\begin{tabular}{c}
\includegraphics[width=6.2in, height=4.5in]{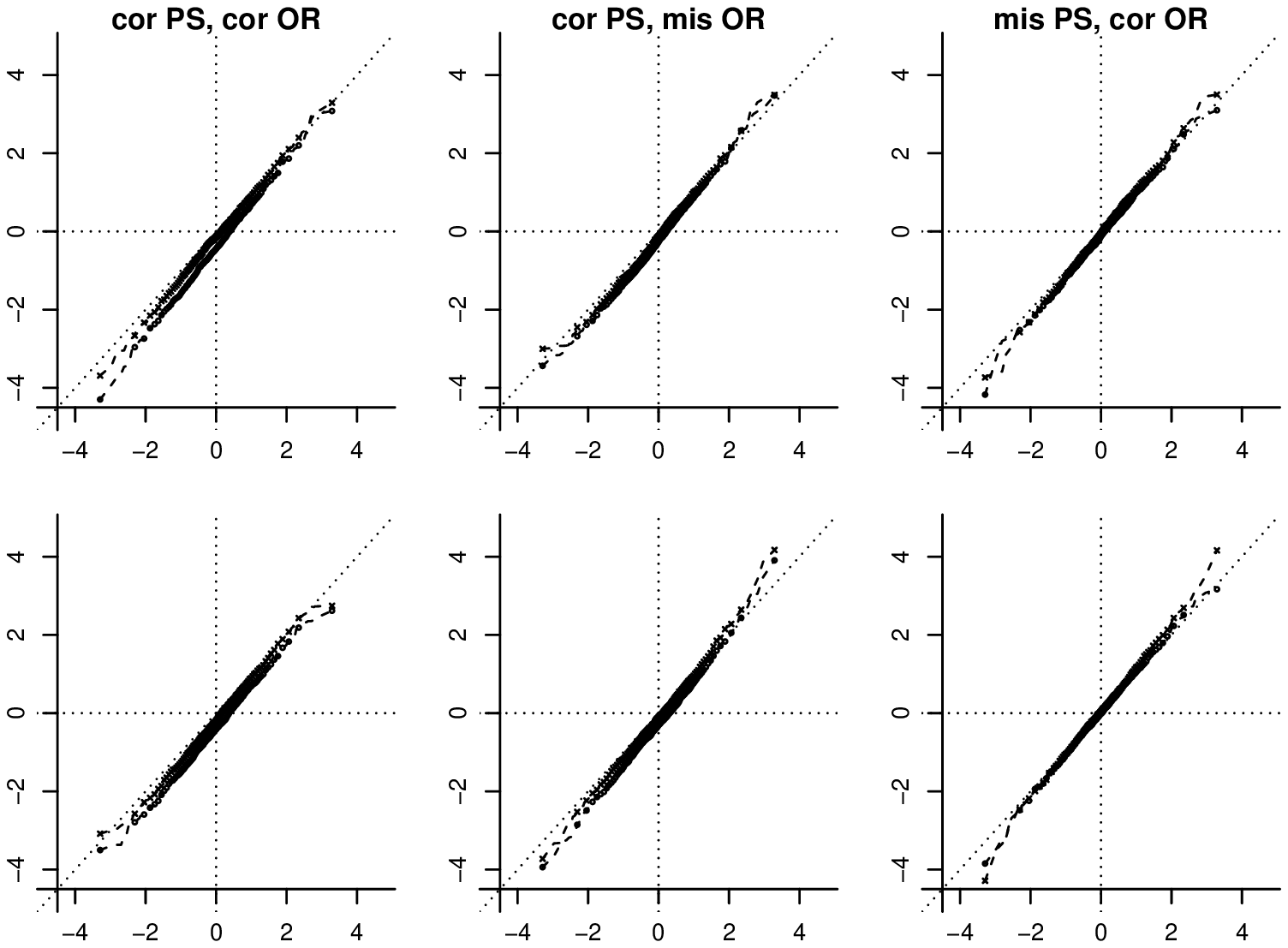} \vspace{-.25in}
\end{tabular}
\end{figure}

\begin{figure}[t!]
\caption{\small QQ plots of the $t$-statistics against standard normal with logistic outcome models ($n=400$, $p=200$),
based on the estimators $\hat \mu^1(\hat m^1_{\mbox{\tiny RML}}, \hat\pi^1_{\mbox{\tiny RML}})$ ($\circ$)
and $\hat \mu^1(\hat m^1_{\mbox{\tiny RWL}}, \hat\pi^1_{\mbox{\tiny RCAL}})$ ($\times$).
For readability, only a subset of 100 order statistics are shown as points on the QQ lines.}
\label{fig:logitOR-n400-p200} \vspace{.15in}
\begin{tabular}{c}
\includegraphics[width=6.2in, height=4.5in]{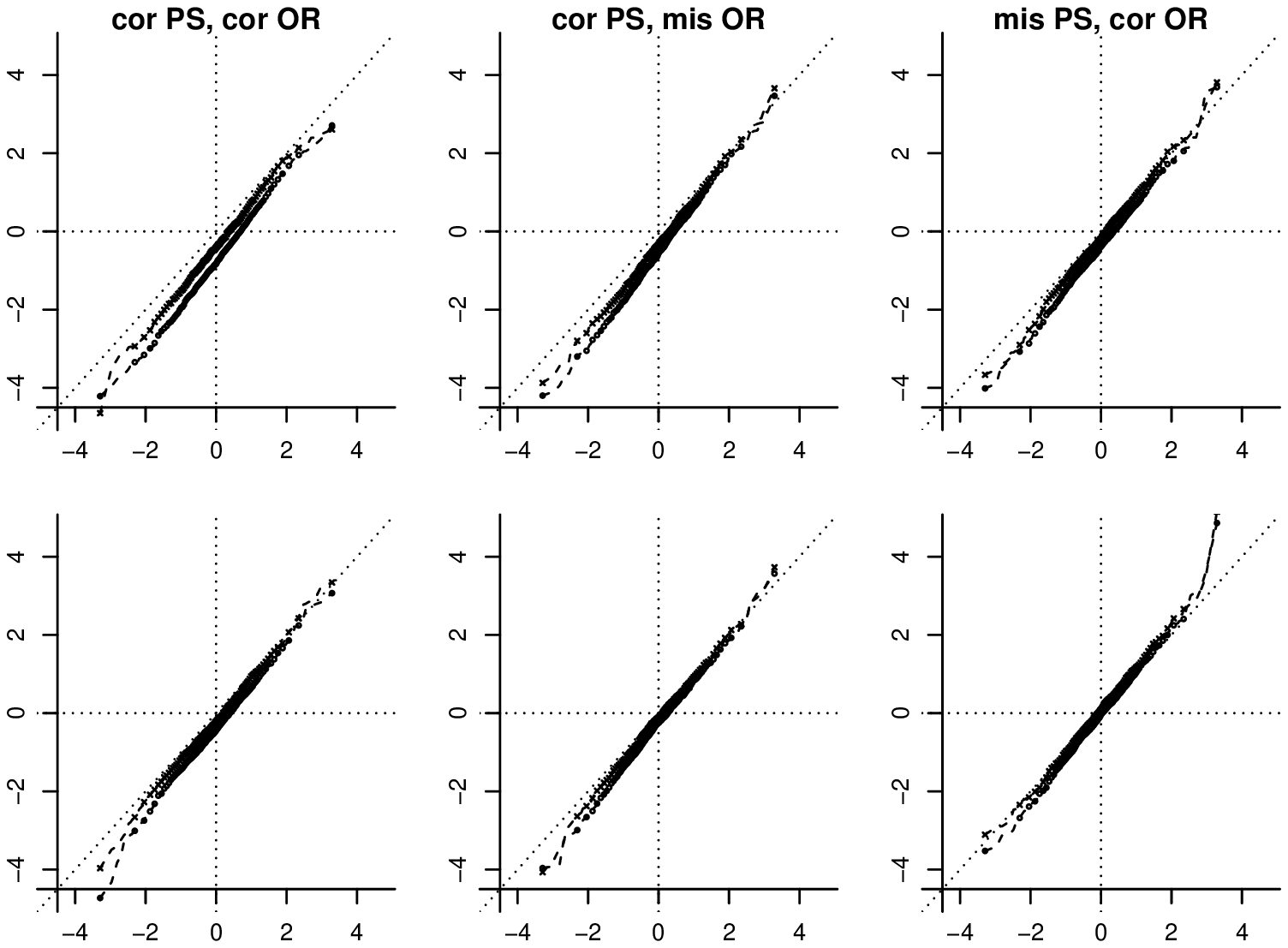} \vspace{-.25in}
\end{tabular}
\end{figure}

Tables~\ref{linearOR-n400-p100}--\ref{linearOR-n400-p200} and Figures~\ref{fig:linearOR-n400-p100}--\ref{fig:linearOR-n400-p200},
present additional simulation results from Section~\ref{sec:sim-linearOR} with linear outcome models, similarly as Table~\ref{linearOR-n800-p200} and Figure~\ref{fig:linearOR-n800-p200}
but for different values of $(n,p)$.
Tables~\ref{logitOR-n400-p100}--\ref{logitOR-n400-p200} and Figures~\ref{fig:logitOR-n400-p100}--\ref{fig:logitOR-n400-p200},
present additional simulation results from Section~\ref{sec:sim-logitOR} with logistic outcome models, similarly as Table~\ref{logitOR-n800-p200} and Figure~\ref{fig:logitOR-n800-p200}
but for different values of $(n,p)$. Similar conclusions can be drawn as discussed in Sections~\ref{sec:sim-linearOR}--\ref{sec:sim-logitOR}.

\clearpage

\begin{table}
\caption{Summary of results with linear outcome models ($n=400$, $p=100$)} \label{linearOR-n400-p100} \vspace{-.02in}
\footnotesize
\begin{center}
\begin{tabular*}{1\textwidth}{@{\extracolsep\fill} c c ll c ll c ll} \hline
      && \multicolumn{2}{c}{cor PS, cor OR} && \multicolumn{2}{c}{cor PS, mis OR} && \multicolumn{2}{c}{mis PS, cor OR} \\ \cline{3-4}  \cline{6-7} \cline{9-10}
      && RML.RML & RCAL.RWL  && RML.RML  &  RCAL.RWL  && RML.RML   &  RCAL.RWL   \\ \hline
      && \multicolumn{8}{c}{Linear outcome configuration 1} \\
Bias  && $-.061$ & $-.041$ && $-.019$ & $-.019$ && $-.031$ & $-.021$ \\
$\sqrt{\mbox{Var}}$
      && .097  & .097 && .097 & .099  && .105 & .0994$^\dag$ \\  %(105/99)^2 = 1.12
$\sqrt{\mbox{EVar}}$
      && .108  & .110 && .111 & .112  && .109 & .109 \\
Cov90 && .787  & .829$^*$ && .844 & .853  && .848 & .845 \\
Cov95 && .862  & .883 && .915 & .916  && .916 & .920  \\
      && \multicolumn{8}{c}{Linear outcome configuration 2} \\
Bias  && $-.048$ & $-.032$ && $-.044$ & $-.025$ && $-.025$ & $-.017$ \\
$\sqrt{\mbox{Var}}$
      && .086  & .087 && .085 & .088 && .093 & .088$^\dag$ \\  % (93/88)^2 = 1.12
$\sqrt{\mbox{EVar}}$
      && .094  & .096 && .092 & .093 && .096 & .096 \\
Cov90 && .799  & .833$^*$ && .828 & .859$^*$ && .864 & .860  \\
Cov95 && .879  & .898 && .896 & .926$^*$ && .927 & .919  \\ \hline
\end{tabular*}\\[.1in]
\parbox{1\textwidth}{\small Note: See the footnote of Table~\ref{linearOR-n800-p200}.}
\end{center}
\end{table}

\begin{table}
\caption{Summary of results with linear outcome models ($n=800$, $p=100$)} \label{linearOR-n800-p100} \vspace{-.02in}
\footnotesize
\begin{center}
\begin{tabular*}{1\textwidth}{@{\extracolsep\fill} c c ll c ll c ll} \hline
      && \multicolumn{2}{c}{cor PS, cor OR} && \multicolumn{2}{c}{cor PS, mis OR} && \multicolumn{2}{c}{mis PS, cor OR} \\ \cline{3-4}  \cline{6-7} \cline{9-10}
      && RML.RML & RCAL.RWL  && RML.RML  &  RCAL.RWL  && RML.RML   &  RCAL.RWL   \\ \hline
      && \multicolumn{8}{c}{Linear outcome configuration 1} \\
Bias  && $-.034$ & $-.022$ && $-.010$ & $-.011$ && $-.006$ & $-.006$ \\
$\sqrt{\mbox{Var}}$
      && .071  & .071 && .073 & .072  && .078 & .072$^\dag$ \\  % (78/72)^2 = 1.17
$\sqrt{\mbox{EVar}}$
      && .079  & .080 && .084 & .083  && .081 & .080 \\
Cov90 && .829  & .836 && .845 & .852  && .889 & .881 \\
Cov95 && .889  & .901 && .905 & .909  && .938 & .929  \\
      && \multicolumn{8}{c}{Linear outcome configuration 2} \\
Bias  && $-.034$ & $-.021$ && $-.044$ & $-.023$ && $-.004$ & $-.004$ \\
$\sqrt{\mbox{Var}}$
      && .064  & .063 && .064 & .064 && .070 & .064$^\dag$ \\  % (70/64)^2 = 1.20
$\sqrt{\mbox{EVar}}$
      && .071  & .072 && .072 & .070 && .072 & .071 \\
Cov90 && .814  & .830 && .782 & .850$^*$ && .896 & .875 \\
Cov95 && .880  & .893 && .862 & .912$^*$ && .941 & .924 \\ \hline
\end{tabular*}\\[.1in]
\parbox{1\textwidth}{\small Note: See the footnote of Table~\ref{linearOR-n800-p200}.}
\end{center}
\end{table}

\begin{table}
\caption{Summary of results with linear outcome models ($n=400$, $p=200$)}  \label{linearOR-n400-p200} \vspace{-.02in}
\footnotesize
\begin{center}
\begin{tabular*}{1\textwidth}{@{\extracolsep\fill} c c ll c ll c ll} \hline
      && \multicolumn{2}{c}{cor PS, cor OR} && \multicolumn{2}{c}{cor PS, mis OR} && \multicolumn{2}{c}{mis PS, cor OR} \\ \cline{3-4}  \cline{6-7} \cline{9-10}
      && RML.RML & RCAL.RWL  && RML.RML  &  RCAL.RWL  && RML.RML   &  RCAL.RWL   \\ \hline
      && \multicolumn{8}{c}{Linear outcome configuration 1} \\
Bias  && $-.068$ & $-.049$ && $-.026$ & $-.024$ && $-.035$ & $-.024$ \\
$\sqrt{\mbox{Var}}$
      && .095  & .096 && .096 & .099  && .102 & .098 \\
$\sqrt{\mbox{EVar}}$
      && .109  & .110 && .112 & .113  && .118 & .116 \\
Cov90 && .770  & .820$^*$ && .819 & .834  && .823 & .829 \\
Cov95 && .845  & .884$^*$ && .895 & .903  && .893 & .896  \\
      && \multicolumn{8}{c}{Linear outcome configuration 2} \\
Bias  && $-.053$ & $-.038$ && $-.045$ & $-.028$ && $-.026$ & $-.017$ \\
$\sqrt{\mbox{Var}}$
      && .085  & .086 && .084 & .087 && .091 & .087 \\
$\sqrt{\mbox{EVar}}$
      && .092  & .094 && .094 & .095 && .103 & .103 \\
Cov90 && .788  & .833$^*$ && .814 & .853$^*$ && .842 & .836 \\
Cov95 && .877  & .905 && .884 & .913 && .904 & .896 \\ \hline
\end{tabular*}\\[.1in]
\parbox{1\textwidth}{\small Note: See the footnote of Table~\ref{linearOR-n800-p200}.}
\end{center}
\end{table}

%%%%%%%%%%%%%%%%%%

\begin{table}
\caption{Summary of results with logistic outcome models ($n=400$, $p=100$)} \label{logitOR-n400-p100} \vspace{-.02in}
\footnotesize
\begin{center}
\begin{tabular*}{1\textwidth}{@{\extracolsep\fill} c c ll c ll c ll} \hline
      && \multicolumn{2}{c}{cor PS, cor OR} && \multicolumn{2}{c}{cor PS, mis OR} && \multicolumn{2}{c}{mis PS, cor OR} \\ \cline{3-4}  \cline{6-7} \cline{9-10}
      && RML.RML & RCAL.RWL  && RML.RML  &  RCAL.RWL  && RML.RML   &  RCAL.RWL   \\ \hline
      && \multicolumn{8}{c}{Logistic outcome configuration 1} \\
Bias  && $-.021$ & $-.011$ && $-.014$ & $-.010$ && $-.014$ & $-.008$ \\
$\sqrt{\mbox{Var}}$
      && .032  & .033 && .032 & .034  && .033 & .033 \\
$\sqrt{\mbox{EVar}}$
      && .038  & .038 && .038 & .038  && .037 & .037 \\
Cov90 && .776  & .835$^*$ && .804 & .845$^*$  && .834 & .852 \\
Cov95 && .864  & .911$^*$ && .877 & .901  && .899 & .904  \\
      && \multicolumn{8}{c}{Logistic outcome configuration 2} \\
Bias  && $-.011$ & $-.004$ && $-.008$ & $-.003$ && $-.003$ & $.002$ \\
$\sqrt{\mbox{Var}}$
      && .033  & .034 && .034 & .035 && .035 & .034 \\
$\sqrt{\mbox{EVar}}$
      && .035  & .036 && .035 & .035 && .038 & .038 \\
Cov90 && .852  & .876 && .876 & .883 && .863 & .864  \\
Cov95 && .910  & .931 && .931 & .942 && .921 & .924  \\ \hline
\end{tabular*}\\[.1in]
\parbox{1\textwidth}{\small Note: See the footnote of Table~\ref{logitOR-n800-p200}.}
\end{center}
\end{table}

\begin{table}
\caption{Summary of results with logistic outcome models ($n=800$, $p=100$)} \label{logitOR-n800-p100}   \vspace{-.02in}
\footnotesize
\begin{center}
\begin{tabular*}{1\textwidth}{@{\extracolsep\fill} c c ll c ll c ll} \hline
      && \multicolumn{2}{c}{cor PS, cor OR} && \multicolumn{2}{c}{cor PS, mis OR} && \multicolumn{2}{c}{mis PS, cor OR} \\ \cline{3-4}  \cline{6-7} \cline{9-10}
      && RML.RML & RCAL.RWL  && RML.RML  &  RCAL.RWL  && RML.RML   &  RCAL.RWL   \\ \hline
      && \multicolumn{8}{c}{Logistic outcome configuration 1} \\
Bias  && $-.010$ & $-.004$ && $-.006$ & $-.003$ && $-.004$ & $-.001$ \\
$\sqrt{\mbox{Var}}$
      && .023  & .024 && .024 & .025  && .024 & .024 \\
$\sqrt{\mbox{EVar}}$
      && .026  & .026 && .026 & .027  && .026 & .026 \\
Cov90 && .816  & .868$^*$ && .851 & .868  && .877 & .870 \\
Cov95 && .896  & .924 && .924 & .929  && .928 & .925  \\
      && \multicolumn{8}{c}{Logistic outcome configuration 2} \\
Bias  && $-.010$ & $-.004$ && $-.007$ & $-.002$ && $.000$ & $.001$ \\
$\sqrt{\mbox{Var}}$
      && .024  & .025 && .025 & .026 && .027 & .025 \\
$\sqrt{\mbox{EVar}}$
      && .026  & .026 && .027 & .028 && .028 & .027 \\
Cov90 && .842  & .870 && .841 & .867 && .881 & .862  \\
Cov95 && .913  & .926 && .911 & .923 && .940 & .925  \\ \hline
\end{tabular*}\\[.1in]
\parbox{1\textwidth}{\small Note: See the footnote of Table~\ref{logitOR-n800-p200}.}
\end{center}
\end{table}

\begin{table}
\caption{Summary of results with logistic outcome models ($n=400$, $p=200$)} \label{logitOR-n400-p200}  \vspace{-.02in}
\footnotesize
\begin{center}
\begin{tabular*}{1\textwidth}{@{\extracolsep\fill} c c ll c ll c ll} \hline
      && \multicolumn{2}{c}{cor PS, cor OR} && \multicolumn{2}{c}{cor PS, mis OR} && \multicolumn{2}{c}{mis PS, cor OR} \\ \cline{3-4}  \cline{6-7} \cline{9-10}
      && RML.RML & RCAL.RWL  && RML.RML  &  RCAL.RWL  && RML.RML   &  RCAL.RWL   \\ \hline
      && \multicolumn{8}{c}{Logistic outcome configuration 1} \\
Bias  && $-.025$ & $-.013$ && $-.019$ & $-.012$ && $-.013$ & $-.005$ \\
$\sqrt{\mbox{Var}}$
      && .032  & .033 && .032 & .034  && .033 & .033 \\
$\sqrt{\mbox{EVar}}$
      && .037  & .037 && .039 & .038  && .036 & .036  \\
Cov90 && .754  & .826$^*$ && .773 & .827$^*$  && .834 & .866 \\
Cov95 && .833  & .907$^*$ && .852 & .897$^*$  && .898 & .917  \\
      && \multicolumn{8}{c}{Logistic outcome configuration 2} \\
Bias  && $-.013$ & $-.006$ && $-.010$ & $-.004$ && $-.002$ & $.003$ \\
$\sqrt{\mbox{Var}}$
      && .032  & .034 && .033 & .035 && .034 & .034 \\
$\sqrt{\mbox{EVar}}$
      && .035  & .036 && .037 & .037 && .037 & .037 \\
Cov90 && .858  & .884 && .848 & .864 && .858 & .857  \\
Cov95 && .915  & .936 && .904 & .932 && .916 & .926  \\ \hline
\end{tabular*}\\[.1in]
\parbox{1\textwidth}{\small Note: See the footnote of Table~\ref{logitOR-n800-p200}.}
\end{center}
\end{table}

%%%%%%%%%%%%%%%%%%

\clearpage
\section{Technical details}

\subsection{Inside Theorem~\ref{thm-RCAL}}
The following result (ii) is taken from Tan (2017), Lemma 1(ii), and
result (i) can be shown similarly using Lemma~\ref{lem-max-bounded} in Section~\ref{sec:tech-tool} and the union bound.

\begin{lem} \label{lem-prob0}
(i) Denoted by $\Omega_0$ the event that
\begin{align*}
\sup_{j=0,1,\ldots,p} \left| \tilde E \left[ \left\{-T \me^{-\bar h^1_{\mbox{\tiny CAL}}(X)} + (1-T) \right\} f_j(X) \right] \right| \le \lambda_0 .
\end{align*}
Under Assumption~\ref{ass-RCAL}(i)--(ii), if $\lambda_0 \ge  \sqrt{2} (\me^{-B_0}+1) C_0 \sqrt{\log\{(1+p)/\epsilon\}/n}$, then $P(\Omega_0) \ge 1- 2\epsilon$.\\
(ii) Denote by $\Omega_1$ the event that
\begin{align}
\sup_{j,k=0,1,\ldots,p} | (\tilde \Sigma_\gamma)_{jk} - (\Sigma_\gamma)_{jk} | \le \lambda_0, \label{lem-prob0-eq}
\end{align}
Under Assumption~\ref{ass-RCAL}(i)--(ii), if $\lambda_0 \ge (4 \me^{-B_0} C_0^2) \sqrt{\log\{(1+p)/\epsilon\}/n}$, then $P( \Omega_1 ) \ge 1- 2 \epsilon^2$.
\end{lem}

Take $\lambda_0 =C_{01} \sqrt{\log\{(1+p)/\epsilon\}/n}$ with
$$
C_{01} = \max\left\{ \sqrt{2} (\me^{-B_0}+1) C_0, 4 \me^{-B_0} C_0^2 \right\}.
$$
Then under the conditions of Theorem~\ref{thm-RCAL}, inequality (\ref{thm-RCAL-eq}) holds in the event $\Omega_0\cap \Omega_1$, with probability at least $1-4\epsilon$, by the proof of Tan (2017, Corollary 2).

\subsection{Probability lemmas}

\begin{lem} \label{lem-prob}
Denote by $\Omega_2$ the event that
\begin{align}
\sup_{j=0,1,\ldots,p} \left| \tilde E \left[ T w(X; \bar\gamma^1_{\mbox{\tiny CAL}})\{Y - \bar m^1_{\mbox{\tiny WL}}(X)\} f_j(X) \right] \right| \le \lambda_1 . \label{lem-prob-eq}
\end{align}
Under Assumptions~\ref{ass-RCAL}(i)--(ii) and \ref{ass-RWL}(i), if $\lambda_1 \ge (\me^{-B_0} C_0 )\sqrt{ 8 (D_0^2 + D_1^2) }\sqrt{\log\{(1+p)/\epsilon\}/n}$, then $P( \Omega_2 ) \ge 1- 2 \epsilon$.
\end{lem}

\begin{prf}
Let $Z_j = T w(X; \bar\gamma^1_{\mbox{\tiny CAL}}) \{ Y - \bar m^1_{\mbox{\tiny WL}}(X)\}f_j(X) $ for $j=0,1\ldots,p$.
Then $E(Z_j)=0$ by the definition of $\bar\alpha^1_{\mbox{\tiny WL}}$. Under Assumption~\ref{ass-RCAL}(i)--(ii), $|Z_j| \le \me^{-B_0} C_0 | T\{ Y- \bar m^1_{\mbox{\tiny WL}}(X)\}|$.
By Assumption~\ref{ass-RWL}(i), the variables $(Z_0,Z_1,\ldots,Z_p)$ are uniformly
sub-gaussian: $\max_{j=0,1,\ldots,p} D_2^2 $ $E\{ \exp(Z_j^2/ D_2^2)  -1 \} \le D_3^2$, with $D_2= \me^{-B_0} C_0 D_0$
and $D_3=\me^{-B_0} C_0 D_1$.
Therefore, Lemma~\ref{lem-prob}(i) holds by Lemma~\ref{lem-max-subG} in Section~\ref{sec:tech-tool} and the union bound.
\end{prf}

Denote  $\Sigma_{\alpha2} = E [ T w(X; \bar\gamma^1_{\mbox{\tiny CAL}}) \{Y- \bar m^1_{\mbox{\tiny WL}}(X) \}^2 f(X) f^\T (X)] $,
and $\tilde \Sigma_{\alpha2} = \tilde E [ T w(X; \bar\gamma^1_{\mbox{\tiny CAL}}) \{Y- \bar m^1_{\mbox{\tiny WL}}(X) \}^2 f(X) f^\T (X)] $,
the sample version of $\tilde \Sigma_{\alpha2}$.

\begin{lem} \label{lem-prob2}
Denote by  $\Omega_3$ the event that
\begin{align}
\sup_{j,k=0,1,\ldots,p} | (\tilde \Sigma_{\alpha2})_{jk} - (\Sigma_{\alpha2})_{jk} | \le (D_0^2 + D_0D_1) \lambda_0, \label{lem-prob2-eq}
\end{align}
Under Assumptions~\ref{ass-RCAL}(i)--(ii) and \ref{ass-RWL}(i), if
$$
(D_0^2 + D_0D_1) \lambda_0 \ge  4 \me^{-B_0} C_0^2 \left[ D_0^2\log\{(1+p)/\epsilon\}/n +  D_0 D_1 \sqrt{\log\{(1+p)/\epsilon\}/n} \right],
$$
then $P( \Omega_3 ) \ge 1- 2 \epsilon^2$.
\end{lem}

\begin{prf}
For any $j,k=0,1,\ldots,p$, the variable $T w(X; \bar\gamma^1_{\mbox{\tiny CAL}}) \{Y- \bar m^1_{\mbox{\tiny WL}}(X) \}^2 f_j(X) f_k(X)$ is the
product of $w(X; \bar\gamma^1_{\mbox{\tiny CAL}})f_j(X) f_k(X)$ and $ T\{Y- \bar m^1_{\mbox{\tiny WL}}(X)\}^2$,
where $| w(X; \bar\gamma^1_{\mbox{\tiny CAL}})f_j(X)$ $ f_k(X)| \le \me^{-B_0} C_0^2$ by Assumptions~\ref{ass-RCAL}(i)--(ii) and $ T\{Y- \bar m^1_{\mbox{\tiny WL}}(X)\} $
is sub-gaussian by Assumption~\ref{ass-RWL}(i).
Applying Lemmas~\ref{lem-max-subE} and \ref{lem-subG-sq} in Section~\ref{sec:tech-tool} yields
$$
P \left\{  | (\tilde \Sigma_{\alpha2})_{jk} - (\Sigma_{\alpha2})_{jk} | > 2\me^{-B_0} C_0^2 D_0^2 t + 2\me^{-B_0} C_0^2 D_0 D_1 t \sqrt{2t} \right \} \le 2 \frac{\epsilon^2}{(1+p)^2},
$$
for $j,k=0,1,\ldots,p$, where $t = \log\{(1+p)^2/\epsilon^2\}/n $. The result then follows from the union bound.
\end{prf}

Denote  $\Sigma_{\alpha1} = E [ T w(X; \bar\gamma^1_{\mbox{\tiny CAL}}) | Y- \bar m^1_{\mbox{\tiny WL}}(X) | f(X) f^\T (X)] $,
and $\tilde \Sigma_{\alpha1} = \tilde E [ T w(X; \bar\gamma^1_{\mbox{\tiny CAL}}) | Y- \bar m^1_{\mbox{\tiny WL}}(X) | f(X) f^\T (X)] $,
the sample version of $\Sigma_{\alpha1}$.

\begin{lem} \label{lem-prob3}
Denote by $\Omega_4$ the event that
\begin{align}
\sup_{j,k=0,1,\ldots,p} | (\tilde \Sigma_{\alpha1})_{jk} - (\Sigma_{\alpha1})_{jk} | \le \sqrt{D_0^2 + D_1^2} \lambda_0, \label{lem-prob3-eq}
\end{align}
Under Assumptions~\ref{ass-RCAL}(i)--(ii) and \ref{ass-RWL}(i), if $\lambda_0 \ge  4 \me^{-B_0} C_0^2 \sqrt{\log\{(1+p)/\epsilon\}/n} $,
then $P( \Omega_3 ) \ge 1- 2 \epsilon^2$.
\end{lem}

\begin{prf}
The variables $T w(X; \bar\gamma^1_{\mbox{\tiny CAL}}) |Y- \bar m^1_{\mbox{\tiny WL}}(X) | f_j(X) f_k(X)$ for
$j,k=0,1,\ldots,p$ are uniformly sub-gaussian, because $| w(X; \bar\gamma^1_{\mbox{\tiny CAL}})f_j(X)$ $ f_k(X)| \le \me^{-B_0} C_0^2$ by Assumptions~\ref{ass-RCAL}(i)--(ii) and $ T |Y- \bar m^1_{\mbox{\tiny WL}}(X) | $
is sub-gaussian by Assumption~\ref{ass-RWL}(i). Applying Lemma~\ref{lem-max-subG} yields
$$
P \left\{  | (\tilde \Sigma_{\alpha1})_{jk} - (\Sigma_{\alpha1})_{jk} | > t \right \} \le 2 \frac{\epsilon^2}{(1+p)^2},
$$
for $j,k=0,1,\ldots,p$, where $t = \me^{-B_0} C_0^2 \sqrt{ 8 (D_0^2 + D_1^2) } \sqrt{\log\{(1+p)^2/\epsilon^2\}/n } $. The result then follows from the union bound.
\end{prf}

Denote  $\Sigma_0 = E [ f(X) f^\T (X)] $ and $\tilde \Sigma_0 = \tilde E [ f(X) f^\T (X)] $, the sample version of $\Sigma_0$.

\begin{lem} \label{lem-prob4}
Denote by $\Omega_5$ the event that
\begin{align}
\sup_{j,k=0,1,\ldots,p} | (\tilde \Sigma_0)_{jk} - (\Sigma_0)_{jk} | \le \me^{B_0} \lambda_0, \label{lem-prob4-eq}
\end{align}
Under Assumption~\ref{ass-RCAL}(i), if $\lambda_0 \ge  4 \me^{-B_0}  C_0^2 \sqrt{\log\{(1+p)/\epsilon\}/n} $,
then $P( \Omega_5 ) \ge 1- 2 \epsilon^2$.
\end{lem}

\begin{prf}
This result follows directly from Lemma~\ref{lem-max-bounded} and the union bound,
with $ |f_j(X) f_k(X)| \le C_0^2 $ and hence
$ |f_j(X) f_k(X) - (\Sigma_0)_{jk} | \le 2 C_0^2$.
\end{prf}

\subsection{Proof of Theorems~\ref{thm-RWL} and \ref{thm-RWL2}}

Thoughout this section, suppose that Assumption~\ref{ass-RCAL} holds.
The proof of Theorem~\ref{thm-RWL2} is completed by combining Lemmas~\ref{lem-prob}--\ref{lem-prob2} and \ref{lem-basic-ineq}--\ref{lem-orac-ineq}.
Theorem~\ref{thm-RWL} is a special case of Theorem~\ref{thm-RWL2}, where Assumptions~\ref{ass-RWL2}(ii)--(iv) are satisfied with $C_2=1$ and $C_3=\eta_2=\eta_3=0$.

\begin{lem} \label{lem-basic-ineq}
For any coefficient vector $\alpha^1$ and $\myeta (X) = \alpha^{\one\T} f(X)$, we have
\begin{align}
& D^\dag_{\mbox{\tiny WL}} (\hat \myeta^1_{\mbox{\tiny RWL}} , \myeta; \hat\gamma^1_{\mbox{\tiny RCAL}} ) + \lambda \| \hat\alpha^1_{\mbox{\tiny RWL}, 1:p}\|_1  \nonumber \\
& \le (\hat\alpha^1_{\mbox{\tiny RWL}} -\alpha^1)^\T
\tilde E \left[ T w(X; \hat\gamma^1_{\mbox{\tiny RCAL}}) \{Y -m_1(X;\alpha^1)\}  f(X) \right] +\lambda \| \alpha^1_{1:p}\|_1 . \label{lem1-ineq1}
\end{align}
\end{lem}

\begin{prf}
For any $u \in (0,1]$, the definition of $\hat\alpha^1_{\mbox{\tiny RWL}}$ implies
\begin{align*}
& \ell_{\mbox{\tiny RWL}} (\hat\alpha^1_{\mbox{\tiny RWL}} ; \hat \gamma^1_{\mbox{\tiny RCAL}}) + \lambda  \| \hat\alpha^1_{\mbox{\tiny RWL}, 1:p}\|_1  \\
& \le \ell_{\mbox{\tiny RWL}} \{(1-u)\hat\alpha^1_{\mbox{\tiny RWL}} + u\alpha^1 ; \hat \gamma^1_{\mbox{\tiny RCAL}} \} + \lambda  \|(1-u)\hat\alpha^1_{\mbox{\tiny RWL}, 1:p} + u\alpha^1_{1:p}\|_1 ,
\end{align*}
which, by the convexity of $\|\cdot\|_1$, gives
\begin{align*}
\ell_{\mbox{\tiny RWL}} (\hat\alpha^1_{\mbox{\tiny RWL}} ; \hat \gamma^1_{\mbox{\tiny RCAL}} ) - \ell_{\mbox{\tiny RWL}}
\{ (1-u) \hat\alpha^1_{\mbox{\tiny RWL}} + u \alpha^1; \hat \gamma^1_{\mbox{\tiny RCAL}} \} + \lambda u \| \hat\alpha^1_{\mbox{\tiny RWL}, 1:p}\|_1 \le \lambda u \| \alpha^1_{1:p}\|_1 .
\end{align*}
Dividing both sides of the preceding inequality by $u$ and letting $u \to 0+$ yields
\begin{align*}
-\tilde E \left[ T w(X; \hat\gamma^1_{\mbox{\tiny RCAL}}) \{ Y - \hat m^1_{\mbox{\tiny RWL}} (X)\} \{ \hat \myeta^1_{\mbox{\tiny RWL}} (X)- \myeta (X) \}  \right]  + \lambda \| \hat\alpha^1_{\mbox{\tiny RWL}, 1:p}\|_1
\le \lambda \| \alpha^1_{1:p}\|_1 ,
\end{align*}
which leads to (\ref{lem1-ineq1}) after simple rearrangement using (\ref{symD-WL}).
\end{prf}

\begin{lem} \label{lem-add-ineq}
In the event $\Omega_0 \cap \Omega_1$, we have
\begin{align}
\tilde E \left[ T w(X;\bar\gamma^1_{\mbox{\tiny CAL}}) \{ \hat \myeta^1_{\mbox{\tiny RCAL}}(X) - \bar \myeta^1_{\mbox{\tiny CAL}}(X)\}^2 \right] \le \me^{\eta_{01}} M_0 |S_\gamma| \lambda_0^2 , \label{lem-add-ineq-eq1}
\end{align}
and for any function $\myeta (X)$,
\begin{align}
 D^\dag_{\mbox{\tiny WL}} ( \hat \myeta^1_{\mbox{\tiny RWL}}, \myeta; \hat\gamma^1_{\mbox{\tiny RCAL}} )
 \ge \me^{-\eta_{01}} D^\dag_{\mbox{\tiny WL}} (\hat \myeta^1_{\mbox{\tiny RWL}} , \myeta ; \bar\gamma^1_{\mbox{\tiny CAL}}), \label{lem-add-ineq-eq2}
\end{align}
where $\eta_{01} = (A_0-1)^{-1} M_0  \eta_0 C_0$. \\
\end{lem}

\begin{prf}
By direct calculation from the definition of $D_{\mbox{\tiny CAL}}()$, we find
\begin{align*}
D_{\mbox{\tiny CAL}} ( \hat\myeta^1_{\mbox{\tiny RCAL}}, \bar\myeta^1_{\mbox{\tiny CAL}})
& = - \tilde E \left[ T \left\{ \me^{-\hat\myeta^1_{\mbox{\tiny RCAL}}(X)} - \me^{-\bar\myeta^1_{\mbox{\tiny CAL}}(X)} \right\} \{\hat\myeta^1_{\mbox{\tiny RCAL}}(X)- \bar\myeta^1_{\mbox{\tiny CAL}}(X) \} \right] \\
& =  \tilde E \left[ T \me^{-u (  \hat\gamma^1_{\mbox{\tiny RCAL}}- \bar\gamma^1_{\mbox{\tiny CAL}} )^\T f(X)} w(X;\bar\gamma^1_{\mbox{\tiny CAL}}) \{\hat\myeta^1_{\mbox{\tiny RCAL}}(X)- \bar\myeta^1_{\mbox{\tiny CAL}}(X) \}^2 \right]
\end{align*}
for some $u \in (0,1)$, where the second step uses the mean value theorem,
\begin{align}
\me^{-\hat\myeta^1_{\mbox{\tiny RCAL}}(X)} - \me^{-\bar\myeta^1_{\mbox{\tiny CAL}}(X)} =
-\me^{-u \hat\myeta^1_{\mbox{\tiny RCAL}}(X) - (1-u)\bar\myeta^1_{\mbox{\tiny CAL}}(X) } (  \hat\gamma^1_{\mbox{\tiny RCAL}}- \bar\gamma^1_{\mbox{\tiny CAL}} )^\T f(X) . \label{mean-val}
\end{align}
In the event $\Omega_0\cap \Omega_1$ that (\ref{thm-RCAL-eq}) holds, we have
\begin{align}
\|  \hat\gamma^1_{\mbox{\tiny RCAL}}- \bar\gamma^1_{\mbox{\tiny CAL}}\|_1 \le (A_0-1)^{-1} M_0 |S_\gamma| \lambda_0
\le (A_0-1)^{-1} M_0  \eta_0 , \label{gamma-L1-bound}
\end{align}
by Assumption~\ref{ass-RCAL}(iv), $|S_\gamma| \lambda_0\le \eta_0$, and hence
\begin{align*}
M_0 |S_\gamma| \lambda_0^0 \ge D_{\mbox{\tiny CAL}} ( \hat\myeta^1_{\mbox{\tiny RCAL}}, \bar\myeta^1_{\mbox{\tiny CAL}})
\ge \me^{-\eta_{01}} \tilde E \left[ T w(X;\bar\gamma^1_{\mbox{\tiny CAL}}) \{\hat\myeta^1_{\mbox{\tiny RCAL}}(X)- \bar\myeta^1_{\mbox{\tiny CAL}}(X) \}^2 \right],
\end{align*}
which gives the desired inequality (\ref{lem-add-ineq-eq1}). In addition, we write
\begin{align*}
& D^\dag_{\mbox{\tiny WL}} (\hat \myeta^1_{\mbox{\tiny RWL}} , \myeta; \hat\gamma^1_{\mbox{\tiny RCAL}} ) \\
& = \tilde E \left( T w(X; \hat\gamma^1_{\mbox{\tiny RCAL}}) \left[ \psi\{\hat \myeta^1_{\mbox{\tiny RWL}} (X)\}- \psi\{\myeta (X)\} \right]
\{ \hat \myeta^1_{\mbox{\tiny RWL}} (X)- \myeta (X) \} \right) \\
& =\tilde E \left( T \me^{-(  \hat\gamma^1_{\mbox{\tiny RCAL}}- \bar\gamma^1_{\mbox{\tiny CAL}} )^\T f(X)} w(X;\bar\gamma^1_{\mbox{\tiny CAL}})
 \left[ \psi\{\hat \myeta^1_{\mbox{\tiny RWL}} (X)\}- \psi\{\myeta (X)\} \right] \{ \hat m^1_{\mbox{\tiny RWL}} (X)- \myeta (X) \} \right),
\end{align*}
which, in the event $\Omega_0\cap \Omega_1$, yields inequality (\ref{lem-add-ineq-eq2}) by (\ref{gamma-L1-bound}) and Assumption~\ref{ass-RCAL}(i).
\end{prf}

For two functions $\myeta(x)$ and $\myeta^\prime(x)$, denote
$$
Q_{\mbox{\tiny WL}} ( \myeta^\prime , \myeta; \gamma ) =
\tilde E \left[T w(X; \gamma) \{ \myeta^\prime (X)- \myeta(X) \}^2 \right].
$$

\begin{lem} \label{lem-add-ineq2}
Take $\alpha^1 = \bar \alpha^1_{\mbox{\tiny WL}}$ and $\myeta (X) = \bar \alpha^{\one\T}_{\mbox{\tiny WL}} f(X)$.
Suppose that Assumption~\ref{ass-RWL}(i) holds. Then in the event $\Omega_0 \cap \Omega_1 \cap \Omega_3$, (\ref{lem1-ineq1}) implies
\begin{align*}
& \me^{-\eta_{01}} D^\dag_{\mbox{\tiny WL}} (\hat \myeta^1_{\mbox{\tiny RWL}} , \myeta; \bar\gamma^1_{\mbox{\tiny CAL}} ) + \lambda \| \hat\alpha^1_{\mbox{\tiny RWL}, 1:p}\|_1  \\
& \le (\hat\alpha^1_{\mbox{\tiny RWL}} -\alpha^1)^\T
\tilde E \left[ T w(X; \bar\gamma^1_{\mbox{\tiny CAL}}) \{ Y -m_1(X;\alpha^1)\} f(X) \right] +\lambda \| \alpha^1_{1:p} \|_1 \\
& \quad +  \me^{\eta_{01}} \left( M_{01} |S_\gamma| \lambda_0^2 \right)^{1/2}
\{ Q_{\mbox{\tiny WL}} (  \hat\myeta^1_{\mbox{\tiny RWL}} , \myeta; \bar\gamma^1_{\mbox{\tiny CAL}} ) \}^{1/2},
\end{align*}
where $M_{01} =(D_0^2 + D_1^2) (\me^{\eta_{01}} M_0 +  \eta_{02}) + (D_0^2 + D_0 D_1) \eta_{02}$, and $\eta_{02} = (A_0-1)^{-2} M_0^2  \eta_0 $.
\end{lem}

\begin{prf}
Consider the following decomposition,
\begin{align}
& (\hat\alpha^1_{\mbox{\tiny RWL}} - \alpha^1)^\T \tilde E \left[ T w(X; \hat\gamma^1_{\mbox{\tiny RCAL}}) \{ Y -m_1(X;\alpha^1)\} f(X) \right] \nonumber \\
& = (\hat\alpha^1_{\mbox{\tiny RWL}} - \alpha^1)^\T \tilde E \left[ T w(X; \bar\gamma^1_{\mbox{\tiny CAL}}) \{ Y - m_1(X;\alpha^1)\} f(X) \right] \nonumber \\
& \quad + \tilde E \left[ T \left\{ \me^{-\hat\myeta^1_{\mbox{\tiny RCAL}}(X)} - \me^{-\bar\myeta^1_{\mbox{\tiny CAL}}(X)} \right\}\{ Y - m_1(X;\alpha^1)\}
 \{ \hat \myeta^1_{\mbox{\tiny RWL}} (X)- \myeta (X)\} \right] ,\label{lem-add-ineq2-eq1}
\end{align}
denoted as $\Delta_1 + \Delta_2$. By the mean value equation (\ref{mean-val}) and the Cauchy--Schwartz inequality, the second term $\Delta_2$ can be bounded from above as
\begin{align}
\Delta_2 &\le \me^{ C_0 \|  \hat\gamma^1_{\mbox{\tiny RCAL}}- \bar\gamma^1_{\mbox{\tiny CAL}}\|_1}
\tilde E^{1/2} \left[ T w(X;\bar\gamma^1_{\mbox{\tiny CAL}})\{ \hat \myeta^1_{\mbox{\tiny RWL}} (X)- \myeta (X)\}^2 \right] \nonumber \\
& \quad \times \tilde E^{1/2} \left[ T w(X;\bar\gamma^1_{\mbox{\tiny CAL}}) \{ Y -m_1(X;\alpha^1)\} ^2
\{\hat\myeta^1_{\mbox{\tiny RCAL}}(X)- \bar\myeta^1_{\mbox{\tiny CAL}}(X) \}^2 \right]. \label{lem-add-ineq2-eq2}
\end{align}
We upper-bound the third term on the right hand side in several steps.
First, in the event $\Omega_3$, we have by inequality (\ref{lem-prob2-eq}),
\begin{align*}
& (\tilde E -E)\left[ T w(X;\bar\gamma^1_{\mbox{\tiny CAL}}) \{ Y -m_1(X;\alpha^1)\} ^2
\{\hat\myeta^1_{\mbox{\tiny RCAL}}(X)- \bar\myeta^1_{\mbox{\tiny CAL}}(X) \}^2 \right]  \\
& \le (D_0^2 + D_0D_1)\lambda_0  \|  \hat\gamma^1_{\mbox{\tiny RCAL}}- \bar\gamma^1_{\mbox{\tiny CAL}}\|_1 ^2,
\end{align*}
where, by some abuse of notation, $(\tilde E-E)(Z)$ denotes $n^{-1}\sum_{i=1}^n \{Z_i - E(Z)\}$ for a variable $Z$ that is a function of $(T,Y,X)$.
Second, by Assumption~\ref{ass-RWL}(i) and Lemma~\ref{lem-subG-moment},  $E[ \{ Y^1 -m_1(X;\alpha^1)\} ^2 |X] \le D_0^2 + D_1^2$ and hence
\begin{align*}
& E \left[ T w(X;\bar\gamma^1_{\mbox{\tiny CAL}}) \{ Y - m_1(X;\alpha^1)\} ^2
\{\hat\myeta^1_{\mbox{\tiny RCAL}}(X)- \bar\myeta^1_{\mbox{\tiny CAL}}(X) \}^2 \right] \\
& \le (D_0^2 + D_1^2) E \left[ T w(X;\bar\gamma^1_{\mbox{\tiny CAL}})
\{\hat\myeta^1_{\mbox{\tiny RCAL}}(X)- \bar\myeta^1_{\mbox{\tiny CAL}}(X) \}^2 \right] .
\end{align*}
Third, in the event $\Omega_1$, we have by inequality (\ref{lem-prob0-eq}),
\begin{align*}
& (E -\tilde E)\left[ T w(X;\bar\gamma^1_{\mbox{\tiny CAL}})
\{\hat\myeta^1_{\mbox{\tiny RCAL}}(X)- \bar\myeta^1_{\mbox{\tiny CAL}}(X) \}^2 \right]   \le \lambda_0 \|  \hat\gamma^1_{\mbox{\tiny RCAL}}- \bar\gamma^1_{\mbox{\tiny CAL}}\|_1 ^2.
\end{align*}
Combining the preceding inequalities, we have in the event $\Omega_1 \cap \Omega_3$,
\begin{align}
& \tilde E \left[ T w(X;\bar\gamma^1_{\mbox{\tiny CAL}}) \{ Y -m_1(X;\alpha^1)\} ^2
\{\hat\myeta^1_{\mbox{\tiny RCAL}}(X)- \bar\myeta^1_{\mbox{\tiny CAL}}(X) \}^2 \right] \nonumber \\
& \le  (D_0^2 + D_0D_1) \lambda_0  \|  \hat\gamma^1_{\mbox{\tiny RCAL}}- \bar\gamma^1_{\mbox{\tiny CAL}}\|_1 ^2 \nonumber \\
& \quad + (D_0^2 + D_1^2)
\left\{  \lambda_0 \|  \hat\gamma^1_{\mbox{\tiny RCAL}}- \bar\gamma^1_{\mbox{\tiny CAL}}\|_1 ^2 +
\tilde E \left[ T w(X;\bar\gamma^1_{\mbox{\tiny CAL}}) \{ \hat \myeta^1_{\mbox{\tiny RCAL}}(X) - \bar \myeta^1_{\mbox{\tiny CAL}}(X)\}^2 \right] \right\}. \label{lem-add-ineq2-eq3}
\end{align}
The desired result follows by collecting inequalities (\ref{lem-add-ineq2-eq1})--(\ref{lem-add-ineq2-eq3}) and applying (\ref{lem-add-ineq-eq1}), (\ref{lem-add-ineq-eq2}) and (\ref{gamma-L1-bound}) in the event
$\Omega_0\cap\Omega_1$.
\end{prf}

\begin{lem} \label{lem2-basic-ineq}
Denote $b=  \hat\alpha^1_{\mbox{\tiny RWL}} - \bar\alpha^1_{\mbox{\tiny WL}}$.
Suppose that Assumption~\ref{ass-RWL}(i) holds. In the event $\Omega_0 \cap \Omega_1\cap\Omega_2 \cap \Omega_3$, we have
\begin{align}
& \me^{-\eta_{01}} D^\dag_{\mbox{\tiny WL}} (  \hat \myeta^1_{\mbox{\tiny RWL}} , \bar \myeta^1_{\mbox{\tiny WL}};\bar\gamma^1_{\mbox{\tiny CAL}} ) +
(A_1-1) \lambda_1 \| b\|_1 \nonumber  \\
& \le \me^{\eta_{01}} \left( M_{01} |S_\gamma| \lambda_0^2 \right)^{1/2}
\{ Q_{\mbox{\tiny WL}} (  \hat \myeta^1_{\mbox{\tiny RWL}} , \bar \myeta^1_{\mbox{\tiny WL}};\bar\gamma^1_{\mbox{\tiny CAL}} ) \}^{1/2}
+ 2 A_1 \lambda_1 \sum_{j\in S_\alpha} |b_j|. \label{lem2-ineq2}
\end{align}
\end{lem}

\begin{prf}
In the event $\Omega_2$, we have
\begin{align*}
  b^\T
 \tilde E \left[ T w(X; \bar\gamma^1_{\mbox{\tiny CAL}}) \{ Y -\bar m^1_{\mbox{\tiny WL}} (X)\} f(X) \right]  \le \lambda_1  \|  b \|_1 .
\end{align*}
From this bound and Lemma~\ref{lem-add-ineq2} with $\alpha^1 = \bar\alpha^1_{\mbox{\tiny WL} }$, we have in the event $\Omega_0 \cap \Omega_1\cap\Omega_2 \cap \Omega_3$,
\begin{align*}
& \me^{-\eta_{01}} D^\dag_{\mbox{\tiny WL}} (  \hat \myeta^1_{\mbox{\tiny RWL}} , \bar \myeta^1_{\mbox{\tiny WL}};\bar\gamma^1_{\mbox{\tiny CAL}} ) + A_1\lambda_1 \| \hat\alpha^1_{\mbox{\tiny RWL}, 1:p} \|_1\\
& \le \lambda_1  \|  b \|_1  + A_1\lambda_1 \| \bar\alpha^1_{\mbox{\tiny WL}, 1:p} \|_1  +  \me^{\eta_{01}} \left( M_{01} |S_\gamma| \lambda_0^2 \right)^{1/2}
\{ Q_{\mbox{\tiny WL}} (  \hat \myeta^1_{\mbox{\tiny RWL}} , \bar \myeta^1_{\mbox{\tiny WL}};\bar\gamma^1_{\mbox{\tiny CAL}} ) \}^{1/2} .
\end{align*}
Applying to the preceding inequality the identity $| \hat\alpha^1_{\mbox{\tiny RWL}, j} | = |\hat\alpha^1_{\mbox{\tiny RWL},j} - \bar\alpha^1_{\mbox{\tiny WL},j} | $ for $j\not\in S_\alpha$
and the triangle inequality
\begin{align*}
| \hat\alpha^1_{\mbox{\tiny RWL}, j} | & \ge |\bar\alpha^1_{\mbox{\tiny WL},j}| - |\hat\alpha^1_{\mbox{\tiny RWL},j} - \bar\alpha^1_{\mbox{\tiny WL},j} |  , \quad j  \in S_\alpha \backslash \{0\},
\end{align*}
and rearranging the result gives
\begin{align*}
& \me^{-\eta_{01}} D^\dag_{\mbox{\tiny WL}} (  \hat \myeta^1_{\mbox{\tiny RWL}} , \bar \myeta^1_{\mbox{\tiny WL}};\bar\gamma^1_{\mbox{\tiny CAL}} ) +
(A_1-1) \lambda_1 \|  b_{1:p} \|_1  \\
& \le \lambda_1 | b_0| + 2 A_1 \lambda_1 \sum_{j\in S_\alpha \backslash\{0\}} | b_j| +  \me^{\eta_{01}} \left( M_{01} |S_\gamma| \lambda_0^2 \right)^{1/2}
\{ Q_{\mbox{\tiny WL}} (  \hat \myeta^1_{\mbox{\tiny RWL}} , \bar \myeta^1_{\mbox{\tiny WL}};\bar\gamma^1_{\mbox{\tiny CAL}} ) \}^{1/2}  .
\end{align*}
The conclusion follows by adding $(A_0-1)\lambda_0 | b_0|$ to both sides above.
\end{prf}

Denote $\tilde \Sigma_{\alpha} = \tilde E [ T w(X; \bar\gamma^1_{\mbox{\tiny CAL}}) \psi_2\{\bar \myeta^1_{\mbox{\tiny WL}}(X)\} f(X) f^\T (X)] $.

\begin{lem} \label{lem-hessian}
Suppose that Assumption~\ref{ass-RWL2}(iii) holds. Then for any $\myeta=\alpha^{\one\T} f$ and $\myeta^\prime = {\alpha^{\one\prime}}^\T f$,
\begin{align*}
D^\dag_{\mbox{\tiny WL}} (\myeta, \myeta^\prime; \bar\gamma^1_{\mbox{\tiny CAL}}) \ge \frac{1-\me^{-C_4 \|b\|_1 } }{C_4 \|b\|_1} \left( b^\T \tilde\Sigma_\alpha b \right) ,
\end{align*}
where $b=\alpha^{\one\prime}- \alpha^1$ and $C_4=C_0C_3$. Throughout, set $(1-\me^{-c})/c=1$ for $c=0$.
\end{lem}

\begin{prf}
Set $\gamma = \bar\gamma^1_{\mbox{\tiny CAL}}$. By direct calculation, we have
\begin{align*}
& D^\dag_{\mbox{\tiny WL}} (\myeta, \myeta^\prime; \gamma )= \tilde E \left(T w(X;\gamma) \left[ \psi \{\myeta^\prime(X)\} - \psi \{\myeta(X)\} \right] \left\{ \myeta^\prime(X) - \myeta(X) \right\} \right) \\
& = \tilde E \left[ T w(X;\gamma) \left( \int_0^1 \psi_2 \left[ \myeta(X) + u \left\{ \myeta^\prime(X)-\myeta(X) \right\} \right]  \dif u \right)  \left\{ \myeta^\prime(X) - \myeta(X) \right\}^2 \right] .
\end{align*}
By Assumption~\ref{ass-RWL2}(iii) and the fact that $ | \myeta^\prime(X) - \myeta(X) | \le \{ \sup_{j=0,1,\ldots,p} |f_j(X) | \}\, \|\alpha^{\one\prime}-\alpha^1\|_1 \le C_0 \|\alpha^{\one\prime}-\alpha^1\|_1$ by Assumption~\ref{ass-RCAL}(i),
it follows that
\begin{align*}
& D^\dag_{\mbox{\tiny WL}} (\myeta, \myeta^\prime; \gamma)
\ge \tilde E \left[ T w(X;\gamma) \left( \int_0^1  \psi_2 \left\{ \myeta(X) \right\} \me^{-C_3 u | \myeta^\prime(X)-\myeta(X)|} \dif u \right) \left\{ \myeta^\prime(X) - \myeta(X) \right\}^2 \right] \\
& \ge \tilde E \left[ T w(X;\gamma)  \psi_2 \left\{ \myeta(X) \right\}  \left\{ \myeta^\prime(X) - \myeta(X) \right\}^2 \right] \left( \int_0^1  \me^{-C_4 u \|\alpha^{\one\prime}-\alpha^1\|_1} \dif u \right),
\end{align*}
which gives the desired result because $\int_0^1 \me^{-cu}\,\dif u = (1-\me^{-c})/c$ for $c \ge 0$.
\end{prf}

\begin{lem} \label{lem-compat}
Suppose that Assumption~\ref{ass-RWL}(iii) holds. In the event $\Omega_1$, Assumption~\ref{ass-RWL}(ii)
implies a compatibility condition for $\tilde\Sigma_\gamma$:
for any vector $b=(b_0,b_1,\ldots,b_p)^\T \in \mathbb R^{1+p} $ such that $\sum_{j\not \in S_\alpha} |b_j| \le \xi_1 \sum_{j\in S_\alpha} |b_j|$, we have
\begin{align}
(1-\eta_1) \nu_1^2  \left(\sum_{j\in S_\alpha} |b_j|\right)^2 \le |S_\alpha| \left( b^\T \tilde \Sigma_\gamma b  \right) . \label{lem-compat-eq}
\end{align}
\end{lem}

\begin{prf}
In the event $\Omega_1$, we have $ |b^\T (\tilde \Sigma_\gamma - \Sigma_\gamma) b | \le \lambda_1 \|b\|_1^2$ by (\ref{lem-prob0-eq}).
Then Assumption~\ref{ass-RWL}(ii) implies that for any vector $b=(b_0,b_1,\ldots,b_p)^\T $ satisfying $\sum_{j\not \in S_\alpha} |b_j| \le \xi_1 \sum_{j\in S_\alpha} |b_j|$,
\begin{align*}
&\nu_1^2 \|b_{S_\alpha} \|_1^2 \le |S_\alpha| (b^\T \Sigma_\gamma b) \le |S_\alpha| \left(b^\T \tilde \Sigma_\gamma b + \lambda_0 \|b\|_1^2 \right) \\
& \le |S_\alpha| (b^\T \tilde \Sigma_\gamma b ) + |S_\alpha| \lambda_1 (1+\xi_1)^2 \|b_{S_\alpha} \|_1^2 ,
\end{align*}
where $\|b_{S_\alpha} \|_1=\sum_{j\in S_\alpha} |b_j|$. The last inequality uses $\|b\|_1 \le (1+\xi_1) \| b_{S_\alpha}\|_1$.
Then (\ref{lem-compat-eq}) follows because $(1+\xi_1)^2 \nu_1^{-2} |S_\alpha| \lambda_1 \le \eta_1 \,(<1)$ by Assumption~\ref{ass-RWL}(iii).
\end{prf}

\begin{lem} \label{lem-orac-ineq}
Suppose that Assumptions~\ref{ass-RWL} and \ref{ass-RWL2} hold, and $A_1 > (\xi_1+1)/(\xi_1-1)$.
In the event $\Omega_0 \cap \Omega_1\cap\Omega_2 \cap \Omega_3$, inequality (\ref{thm-RWL-eq}) holds as in Theorem~\ref{thm-RWL}.
\end{lem}

\begin{prf}
Denote $b=  \hat\alpha^1_{\mbox{\tiny RWL}} - \bar\alpha^1_{\mbox{\tiny WL}}$,
$D^\dag_{\mbox{\tiny WL}} = D^\dag_{\mbox{\tiny WL}} (  \hat \myeta^1_{\mbox{\tiny RWL}} , \bar \myeta^1_{\mbox{\tiny WL}};  \bar \gamma^1_{\mbox{\tiny CAL}})$,
$Q_{\mbox{\tiny WL}} = Q_{\mbox{\tiny WL}} (  \hat \myeta^1_{\mbox{\tiny RWL}} , \bar \myeta^1_{\mbox{\tiny WL}} ; \bar \gamma^1_{\mbox{\tiny CAL}})$, and
$$
D^\ddag_{\mbox{\tiny WL}} = \me^{-\eta_{01}} D^\dag_{\mbox{\tiny WL}} +
(A_1-1) \lambda_1 \| b \|_1 .
$$
In the event $\Omega_0 \cap \Omega_1\cap\Omega_2 \cap \Omega_3$, inequality (\ref{lem2-ineq2})
from Lemma~\ref{lem2-basic-ineq} with Assumption~\ref{ass-RWL}(i) leads to two possible cases: either
\begin{align}
\xi_2 D^\ddag_{\mbox{\tiny WL}}  \le \me^{\eta_{01}} \left( M_{01} |S_\gamma| \lambda_0^2\right)^{1/2}
( Q_{\mbox{\tiny WL}} )^{1/2}, \label{lem3-ineq1}
\end{align}
or $(1-\xi_2) D^\ddag_{\mbox{\tiny WL}} \le 2 A_1 \lambda_1 \sum_{j\in S_\alpha} |b_j|$, that is,
\begin{align}
D^\ddag_{\mbox{\tiny WL}}
\le (\xi_1+1) (A_1-1) \lambda_1 \sum_{j\in S_\alpha} |b_j| = \xi_3 \lambda_1 \sum_{j\in S_\alpha} |b_j| , \label{lem3-ineq2}
\end{align}
where $\xi_2 = 1-2A_1 /\{ (\xi_1+1)(A_1-1)\} \in (0,1]$ because $A_1 > (\xi_1+1)/(\xi_1-1)$ and $\xi_3= (\xi_1+1)(A_1-1)$.
We deal with the two cases separately as follows.

If (\ref{lem3-ineq2}) holds, then $\sum_{j \not\in S_\alpha} |b_j| \le \xi_1 \sum_{j\in S_\alpha} |b_j|$,
which, by Lemma~\ref{lem-compat} and Assumption~\ref{ass-RWL}(ii)--(iii), implies (\ref{lem-compat-eq}), that is,
\begin{align}
\sum_{j\in S_\alpha} |b_j| \le (1-\eta_1)^{-1/2}\nu_1^{-1} |S_\alpha|^{1/2} \left( b^\T \tilde \Sigma_\gamma b  \right)^{1/2}. \label{lem3-ineq3}
\end{align}
By Assumption~\ref{ass-RWL2}(ii) and Lemma~\ref{lem-hessian} with Assumption~\ref{ass-RWL2}(iii), we have
\begin{align}
D^\dag_{\mbox{\tiny WL}} \ge \frac{1-\me^{-C_4 \|b\|_1 } }{C_4 \|b\|_1} \left( b^\T \tilde \Sigma_\alpha b\right)
\ge \frac{1-\me^{-C_4 \|b\|_1 } }{C_4 \|b\|_1} C_2 \left( b^\T \tilde \Sigma_\gamma b\right) . \label{lem3-ineq4}
\end{align}
Combining (\ref{lem3-ineq2}), (\ref{lem3-ineq3}), and (\ref{lem3-ineq4}) and using $D^\dag_{\mbox{\tiny WL}} \le \me^{\eta_{01}} D^\ddag_{\mbox{\tiny WL}}$ yields
\begin{align}
D^\ddag_{\mbox{\tiny WL}} \le  \me^{\eta_{01}} \xi_3^2(1-\eta_1)^{-1} \nu_1^{-2} C_2^{-1}|S_\alpha|  \lambda_1^2  \frac{C_4 \|b\|_1}{1-\me^{-C_4 \|b\|_1 } }. \label{lem3-ineq5}
\end{align}
But $(A_1-1)\lambda_1 \|b\|_1 \le D^\ddag_{\mbox{\tiny WL}}$. Inequality (\ref{lem3-ineq5}) along with Assumption~\ref{ass-RWL2}(iv) implies that
$1 - \me^{- C_4 \|b\|_1} \le C_4 (A_1-1)^{-1} \xi_3^2 (1-\eta_1)^{-1} \nu_1^{-2} C_2^{-1} |S_\alpha| \lambda_1 \le \eta_2 \, (<1)$.
As a result, $C_4 \|b\|_1 \le - \log(1- \eta_2)$ and hence
\begin{align*}
\frac{1-\me^{-C_4 \|b\|_1 } }{C_4 \|b\|_1} = \int_0^1 \me^{-C_4 \|b\|_1 u} \dif u \ge \me^{-C_4 \|b\|_1} \ge 1-\eta_2.
\end{align*}
From this bound, inequality (\ref{lem3-ineq5}) then leads to $D^\ddag_{\mbox{\tiny WL}} \le  \me^{\eta_{01}}  \xi_3^2 \nu_3^{-2}|S_\alpha| \lambda_1^2  $.

If (\ref{lem3-ineq1}) holds, then simple manipulation using $D^\dag_{\mbox{\tiny WL}} \le \me^{\eta_{01}} D^\ddag_{\mbox{\tiny WL}}$ and (\ref{lem3-ineq4}) together with
$Q_{\mbox{\tiny WL}} =  b^\T \tilde \Sigma_\gamma b$ gives
\begin{align}
D^\ddag_{\mbox{\tiny WL}}  \le  \me^{3\eta_{01}} \xi_2^{-2} C_2^{-1}  \left(  M_{01} |S_\gamma| \lambda_0^2 \right) \frac{C_4 \|b\|_1}{1-\me^{-C_4 \|b\|_1 } } . \label{lem3-ineq1b}
\end{align}
Similarly as above, using $(A_1-1)\lambda_1 \|b\|_1 \le D^\ddag_{\mbox{\tiny WL}}$ and inequality (\ref{lem3-ineq1b}) along with Assumption~\ref{ass-RWL2}(iv), we find
$1 - \me^{- C_4 \|b\|_1} \le C_4 \me^{3\eta_{01}} (A_1-1)^{-1}\xi_2^{-2} C_2^{-1} ( M_{01} |S_\gamma| \lambda_0 )\le \eta_3 \, (<1)$.
As a result, $C_4 \|b\|_1 \le - \log(1- \eta_3)$ and hence
\begin{align*}
\frac{1-\me^{-C_4 \|b\|_1 } }{C_4 \|b\|_1} = \int_0^1 \me^{-C_4 \|b\|_1 u} \dif u \ge \me^{-C_4 \|b\|_1} \ge 1-\eta_3.
\end{align*}
From this bound, inequality (\ref{lem3-ineq1b}) then leads to $D^\ddag_{\mbox{\tiny WL}} \le  \me^{3\eta_{01}}  \xi_4^{-2}  ( M_{01} |S_\gamma| \lambda_0^2 )$.
Therefore, (\ref{thm-RWL2-eq}) holds through (\ref{lem3-ineq1}) and (\ref{lem3-ineq2}) in the event $\Omega_0 \cap \Omega_1\cap\Omega_2 \cap \Omega_3$.
\end{prf}

\subsection{Proof of Theorem~\ref{thm-mu}} \label{sec:prf-thm-mu}

Denote $\hat\varphi = \varphi(T,Y,X; \hat m^1_{\mbox{\tiny RWL}}, \hat \pi^1_{\mbox{\tiny RCAL}})$
and $\bar\varphi = \varphi(T,Y,X; \bar m^1_{\mbox{\tiny WL}}, \bar \pi^1_{\mbox{\tiny CAL}})$. Then
\begin{align*}
\hat \mu^1(\hat m^1_{\mbox{\tiny RWL}}, \hat \pi^1_{\mbox{\tiny RCAL}}) = \bar \mu^1(\bar m^1_{\mbox{\tiny WL}}, \bar \pi^1_{\mbox{\tiny CAL}})
+ \tilde E ( \hat \varphi - \bar\varphi ) .
\end{align*}
Consider the following decomposition,
\begin{align}
\hat \varphi - \bar\varphi & =  \{\hat m^1_{\mbox{\tiny RWL}}(X) - \bar m^1_{\mbox{\tiny WL}}(X)\} \left\{ 1 - \frac{T}{\bar \pi^1_{\mbox{\tiny CAL}}(X)} \right\} \nonumber \\
& \quad +  T\{Y-\bar m^1_{\mbox{\tiny WL}}(X)\} \left\{ \frac{1}{\hat \pi^1_{\mbox{\tiny RCAL}}(X)} - \frac{1}{\bar \pi^1_{\mbox{\tiny CAL}}(X)} \right\} \nonumber \\
& \quad +  \{\hat m^1_{\mbox{\tiny RWL}}(X) - \bar m^1_{\mbox{\tiny WL}}(X)\} \left\{ \frac{T}{\bar \pi^1_{\mbox{\tiny CAL}}(X)} - \frac{T}{\hat \pi^1_{\mbox{\tiny RCAL}}(X)} \right\} ,
\label{phi-decomp}
\end{align}
denoted as $\delta_1 + \delta_2 + \delta_3$.

We show that in the event $\Omega_0 \cap \Omega_1\cap\Omega_2 \cap \Omega_3 \cap \Omega_4$, inequality (\ref{thm-mu-eq}) holds as in Theorem~\ref{thm-mu}.
The decomposition (\ref{mu-decomp}) for $ \hat \mu^1(\hat m^1_{\mbox{\tiny RWL}}, \hat \pi^1_{\mbox{\tiny RCAL}})$ amounts to
\begin{align*}
\hat \mu^1(\hat m^1_{\mbox{\tiny RWL}}, \hat \pi^1_{\mbox{\tiny RCAL}}) & = \bar \mu^1(\bar m^1_{\mbox{\tiny WL}}, \bar \pi^1_{\mbox{\tiny CAL}}) +
\Delta_1 + \Delta_2,
\end{align*}
where
\begin{align*}
\Delta_1 & = \tilde E(\delta_1 + \delta_3)= ( \hat \alpha^1_{\mbox{\tiny RWL}}  - \bar \alpha^1_{\mbox{\tiny WL}} ) ^\T \tilde E \left[\left\{ 1 - \frac{T}{\hat \pi^1_{\mbox{\tiny RCAL}}(X)} \right\} f(X)  \right]  , \\
\Delta_2 & = \tilde E (\delta_2) = \tilde E \left[ T\{Y- \bar m^1_{\mbox{\tiny WL}}(X)\} \left\{ \frac{1}{\hat \pi^1_{\mbox{\tiny RCAL}}(X)} - \frac{1}{\bar \pi^1_{\mbox{\tiny CAL}}(X)} \right\} \right] .
\end{align*}
In the event $\Omega_0 \cap \Omega_1\cap\Omega_2 \cap \Omega_3$, we have
\begin{align}
| \Delta_1 | \le (A_1-1)^{-1} M_1 (|S_\gamma| \lambda_0 + |S_\alpha| \lambda_1) \times A_0 \lambda_0 , \label{prf-thm-mu-eq1}
\end{align}
by inequality (\ref{thm-RWL-eq}) and the Karush--Kuhn--Tucker conditions (\ref{ineq-CAL-1})--(\ref{ineq-CAL-2}).
Moreover, a Taylor expansion for $\Delta_2$ yields for some $u \in (0,1)$,
\begin{align*}
& \Delta_2 = -( \hat \gamma^1_{\mbox{\tiny RCAL}}  - \bar \gamma^1_{\mbox{\tiny CAL}} ) ^\T \tilde E \left[ T\{Y- \bar m^1_{\mbox{\tiny WL}}(X)\}
\me^{-\bar \myeta^1_{\mbox{\tiny CAL}} (X) } f(X) \right]  \\
& \quad + ( \hat \gamma^1_{\mbox{\tiny RCAL}}  - \bar \gamma^1_{\mbox{\tiny CAL}} ) ^\T \tilde E \left[ T\{Y- \bar m^1_{\mbox{\tiny WL}}(X)\}
\me^{-u \hat \myeta^1_{\mbox{\tiny RCAL}}(X) - (1-u) \bar \myeta^1_{\mbox{\tiny CAL}} (X) } f(X) f^\T(X) \right]( \hat \gamma^1_{\mbox{\tiny RCAL}}  - \bar \gamma^1_{\mbox{\tiny CAL}} )/2,
\end{align*}
denoted as $\Delta_{21} + \Delta_{22}$.
In the event $(\Omega_0 \cap\Omega_1) \cap \Omega_2$, we have
\begin{align}
| \Delta_{21} | \le (A_0-1)^{-1} M_0 |S_\gamma| \lambda_0 \times \lambda_1,  \label{prf-thm-mu-eq2}
\end{align}
by inequalities (\ref{thm-RCAL-eq}) and (\ref{lem-prob-eq}).
The term $\Delta_{22}$ can be bounded as
\begin{align}
 | \Delta_{22} | \le \me^{\|\hat \gamma^1_{\mbox{\tiny RCAL}}-\bar \gamma^1_{\mbox{\tiny CAL}}\|_1 C_0}  \tilde E \left[ T w(X; \bar \gamma^1_{\mbox{\tiny CAL}}) |Y- \bar m^1_{\mbox{\tiny WL}}(X) |
 \{ \hat \myeta^1_{\mbox{\tiny RCAL}}(X)  - \bar \myeta^1_{\mbox{\tiny CAL}} (X) \}^2 \right] /2 . \label{prf-thm-mu-eq3}
\end{align}
In the event $\Omega_1 \cap \Omega_4$, we have
\begin{align}
& \tilde E \left[ T w(X;\bar\gamma^1_{\mbox{\tiny CAL}}) | Y -\bar m^1_{\mbox{\tiny WL}} (X)|
\{\hat\myeta^1_{\mbox{\tiny RCAL}}(X)- \bar\myeta^1_{\mbox{\tiny CAL}}(X) \}^2 \right] \nonumber \\
& \le  \sqrt{D_0^2+D_1^2} \lambda_0  \|  \hat\gamma^1_{\mbox{\tiny RCAL}}- \bar\gamma^1_{\mbox{\tiny CAL}}\|_1 ^2 \nonumber \\
& \quad +  \sqrt{D_0^2+D_1^2}
\left\{  \lambda_0 \|  \hat\gamma^1_{\mbox{\tiny RCAL}}- \bar\gamma^1_{\mbox{\tiny CAL}}\|_1 ^2 +
\tilde E \left[ T w(X;\bar\gamma^1_{\mbox{\tiny CAL}}) \{ \hat \myeta^1_{\mbox{\tiny RCAL}}(X) - \bar \myeta^1_{\mbox{\tiny CAL}}(X)\}^2 \right] \right\},  \label{prf-thm-mu-eq4}
\end{align}
by inequalities (\ref{lem-prob0-eq}) and (\ref{lem-prob3-eq}) and  similar steps as in the proof of (\ref{lem-add-ineq2-eq3}).
Then (\ref{thm-mu-eq}) follows by collecting inequalities  (\ref{prf-thm-mu-eq1})--(\ref{prf-thm-mu-eq4}) and applying (\ref{lem-add-ineq-eq1}) and (\ref{gamma-L1-bound}) in the event
$\Omega_0 \cap \Omega_1$.

\subsection{Proof of Theorem~\ref{thm-var}}

Using $a^2-b^2 = 2(a-b)b + (a-b)^2$ and the Cauchy--Schwartz inequality, we find
\begin{align}
& \left| \tilde E \left(\hat \varphi_c^2 - \bar \varphi_c^2 \right) \right|
\le 2 \tilde E^{1/2} \left( \bar \varphi_c^2 \right) \tilde E^{1/2} \left\{ ( \hat \varphi_c -  \bar \varphi_c)^2 \right\} + \tilde E \left\{ (\hat \varphi_c -  \bar \varphi_c)^2 \right\}. \label{prf-thm-var-eq1}
\end{align}
Using $\hat\varphi_c = \hat \varphi -\hat \mu^1(\hat m^1_{\mbox{\tiny RWL}}, \hat\pi^1_{\mbox{\tiny RCAL}}) $ and
$\bar\varphi_c = \bar \varphi -\hat \mu^1(\bar m^1_{\mbox{\tiny WL}}, \bar\pi^1_{\mbox{\tiny CAL}}) $, we find
\begin{align}
 \tilde E \{ (\hat \varphi_c -  \bar \varphi_c)^2 \} \le 2  \tilde E \{ (\hat \varphi-  \bar \varphi)^2 \} +
 2 \left| \hat \mu^1(\hat m^1_{\mbox{\tiny RWL}}, \hat\pi^1_{\mbox{\tiny RCAL}}) - \hat \mu^1(\bar m^1_{\mbox{\tiny WL}}, \bar\pi^1_{\mbox{\tiny CAL}})  \right|^2 . \label{prf-thm-var-eq2}
\end{align}
To control $\tilde E \{ (\hat \varphi-  \bar \varphi)^2 \}$, we use the decomposition (\ref{phi-decomp}), denoted as $\delta_1 + \delta_2 + \delta_3$.

First, by the mean value equation (\ref{mean-val}) and Assumption~\ref{ass-RCAL}(i)--(ii), we have
\begin{align}
& \tilde E (\delta_2^2 ) = \tilde E \left[  T\{Y-\bar m^1_{\mbox{\tiny WL}}(X)\}^2 \left\{ \frac{1}{\hat \pi^1_{\mbox{\tiny RCAL}}(X)} - \frac{1}{\bar \pi^1_{\mbox{\tiny CAL}}(X)} \right\}^2 \right] \nonumber \\
& \le \me^{ -B_0 + 2 \|\hat \gamma^1_{\mbox{\tiny RCAL}} - \bar \gamma^1_{\mbox{\tiny CAL}}\|_1 C_0}
 \tilde E \left[  T w(X; \bar \gamma^1_{\mbox{\tiny CAL}}) \{Y-\bar m^1_{\mbox{\tiny WL}}(X)\}^2 \{ \hat h^1_{\mbox{\tiny RCAL}}(X) - \bar h^1_{\mbox{\tiny CAL}}(X) \}^2 \right] . \label{prf-thm-var-eq3}
\end{align}
Second, writing $\{\hat \pi^1_{\mbox{\tiny RCAL}}(X)\}^{-1} - \{\bar \pi^1_{\mbox{\tiny CAL}}(X) \}^{-1}=
\me^{-\bar h^1_{\mbox{\tiny CAL}}(X)} \{ \me^{-\hat h^1_{\mbox{\tiny RCAL}}(X)+\bar h^1_{\mbox{\tiny CAL}}(X)} -1\}$
and using Assumption~\ref{ass-RCAL}(i)--(ii), we have
\begin{align}
& \tilde E (\delta_3^2 ) = \tilde E \left[  T\{\hat m^1_{\mbox{\tiny RWL}}(X) - \bar m^1_{\mbox{\tiny WL}}(X)\}^2
\left\{ \frac{1}{\hat \pi^1_{\mbox{\tiny RCAL}}(X)} - \frac{1}{\bar \pi^1_{\mbox{\tiny CAL}}(X)} \right\}^2 \right]  \nonumber \\
& \le \me^{ -B_0} \left(1+\me^{\|\hat \gamma^1_{\mbox{\tiny RCAL}} - \bar \gamma^1_{\mbox{\tiny CAL}}\|_1 C_0} \right)^2
\tilde E \left[  T w(X; \bar \gamma^1_{\mbox{\tiny CAL}}) \{\hat m^1_{\mbox{\tiny RWL}}(X) - \bar m^1_{\mbox{\tiny WL}}(X)\}^2 \right]. \label{prf-thm-var-eq4}
\end{align}
Third, using Assumption~\ref{ass-RCAL}(i)--(ii), we also have
\begin{align}
& \tilde E (\delta_1^2 ) = \tilde E \left[ \{\hat m^1_{\mbox{\tiny RWL}}(X) - \bar m^1_{\mbox{\tiny WL}}(X)\}^2 \left\{ 1 - \frac{T}{\bar \pi^1_{\mbox{\tiny CAL}}(X)} \right\}^2 \right] \nonumber \\
& \le (1 + \me^{-B_0})^2 \tilde E \left[ \{\hat \myeta^1_{\mbox{\tiny RWL}}(X) - \bar \myeta^1_{\mbox{\tiny WL}}(X)\}^2 \right] \label{prf-thm-var-eq5} \\
& \le (1 + \me^{-B_0})^2 C_0^2 \| \hat \alpha^1_{\mbox{\tiny RWL}} - \bar \alpha^1_{\mbox{\tiny WL}} \|_1^2 . \label{prf-thm-var-eq6}
\end{align}
Inequality (\ref{thm-var-eq1}) follows by collecting inequalities (\ref{prf-thm-var-eq1})--(\ref{prf-thm-var-eq6}) and
applying (\ref{thm-RWL-eq}), (\ref{thm-mu-eq}), (\ref{gamma-L1-bound}), and (\ref{lem-add-ineq2-eq3})
in the event $\Omega_0 \cap \Omega_1 \cap \Omega_2 \cap \Omega_3 \cap \Omega_4$.
If condition (\ref{compat-add}) holds, then we have in the event $\Omega_1 \cap \Omega_5$,
\begin{align}
& \tilde E \left[\{\hat \myeta^1_{\mbox{\tiny RWL}}(X) - \bar \myeta^1_{\mbox{\tiny WL}}(X)\}^2 \right]
\le  \me^{B_0} \lambda_0  \|  \hat\alpha^1_{\mbox{\tiny RWL}}- \bar\alpha^1_{\mbox{\tiny WL}}\|_1 ^2  \nonumber \\
& \quad + \tau_0^{-1}
\left\{  \lambda_0 \|  \hat\alpha^1_{\mbox{\tiny RWL}}- \bar\alpha^1_{\mbox{\tiny WL}}\|_1 ^2 +
\tilde E \left[T w(X; \bar \alpha^1_{\mbox{\tiny WL}}) \{\hat \myeta^1_{\mbox{\tiny RWL}}(X) - \bar \myeta^1_{\mbox{\tiny WL}}(X)\}^2 \right] \right\}, \label{prf-thm-var-eq7}
\end{align}
by inequalities (\ref{lem-prob0-eq}) and (\ref{lem-prob4-eq}) and  similar steps as in the proof of (\ref{lem-add-ineq2-eq3}).
Inequality (\ref{thm-var-eq2}) follows, similarly as (\ref{thm-var-eq1}), by combining inequalities (\ref{prf-thm-var-eq1})--(\ref{prf-thm-var-eq5}) and (\ref{prf-thm-var-eq7}).

\subsection{Proof of Theorem~\ref{thm-mu2}}

We use the decomposition (\ref{phi-decomp}) and  handle $\delta_1$, $\delta_2$, and $\delta_3$ separately.
The term $\tilde E(\delta_2)$ can be bounded by (\ref{prf-thm-mu-eq2})--(\ref{prf-thm-mu-eq4}) as in the proof of Theorem~\ref{thm-mu}.
By the mean value equation (\ref{mean-val}) and the Cauchy--Schwartz inequality,  $\tilde E(\delta_3)$ can be bounded as
\begin{align}
\left| \tilde E(\delta_3) \right| &\le \me^{  C_0  \| \hat\gamma^1_{\mbox{\tiny RCAL}}- \bar\gamma^1_{\mbox{\tiny CAL}}\|_1}
\tilde E^{1/2} \left[ T w(X;\bar\gamma^1_{\mbox{\tiny CAL}})
\{\hat\myeta^1_{\mbox{\tiny RCAL}}(X)- \bar\myeta^1_{\mbox{\tiny CAL}}(X) \}^2 \right] \nonumber \\
& \quad \times \tilde E^{1/2} \left[ T w(X;\bar\gamma^1_{\mbox{\tiny CAL}})\{ \hat m^1_{\mbox{\tiny RWL}} (X)- \bar m^1_{\mbox{\tiny WL}} (X)\}^2 \right] .    \label{prf-thm-mu2-eq1}
\end{align}
Similarly as in Lemma~\ref{lem-hessian} but arguing in the reverse direction by Assumptions~\ref{ass-RCAL}(i) and \ref{ass-RWL2}(iii), we find
\begin{align}
&\tilde E \left[ T w(X;\bar\gamma^1_{\mbox{\tiny CAL}})\{ \hat m^1_{\mbox{\tiny RWL}} (X)- \bar m^1_{\mbox{\tiny WL}} (X)\}^2 \right] \le
\me^{C_4 \| \hat \alpha^1_{\mbox{\tiny RWL}} - \bar \alpha^1_{\mbox{\tiny WL}}\|_1}\nonumber  \\
& \quad \times \tilde E \left[ T w(X;\bar\gamma^1_{\mbox{\tiny CAL}}) \psi_2\{\bar \myeta^1_{\mbox{\tiny WL}} (X)\} \{ \hat m^1_{\mbox{\tiny RWL}} (X)- \bar m^1_{\mbox{\tiny WL}} (X)\}
\{ \hat \myeta^1_{\mbox{\tiny RWL}} (X)- \bar \myeta^1_{\mbox{\tiny WL}} (X)\} \right] \nonumber \\
& \le \me^{C_4 \| \hat \alpha^1_{\mbox{\tiny RWL}} - \bar \alpha^1_{\mbox{\tiny WL}}\|_1} C_1 D^\dag_{\mbox{\tiny WL}}(\hat m^1_{\mbox{\tiny RWL}}, \bar m^1_{\mbox{\tiny WL}}; \bar \gamma^1_{\mbox{\tiny CAL}}), \label{prf-thm-mu2-eq2}
\end{align}
where the second inequality follows from Assumption~\ref{ass-RWL2}(i).
In the following, we derive two different bounds on $\tilde E(\delta_1)$, leading to (\ref{thm-mu2-eq1}) and (\ref{thm-mu2-eq2}) respectively.

First, suppose that condition (\ref{compat-add}) holds. Consider the following decomposition
\begin{align}
\tilde E ( \delta_1 )& = \tilde E
\left[ \psi_2\{\bar \myeta^1_{\mbox{\tiny WL}}(X)\}  \{\hat \myeta^1_{\mbox{\tiny RWL}}(X) - \bar \myeta^1_{\mbox{\tiny WL}}(X)\} \left\{ 1 - \frac{T}{\bar \pi^1_{\mbox{\tiny CAL}}(X)} \right\}  \right] \nonumber \\
& \quad + \tilde E
\left[ \tilde\psi_2(X) \{\hat \myeta^1_{\mbox{\tiny RWL}}(X) - \bar \myeta^1_{\mbox{\tiny WL}}(X)\} \left\{ 1 - \frac{T}{\bar \pi^1_{\mbox{\tiny CAL}}(X)} \right\}  \right], \label{prf-thm-mu2-eq3}
\end{align}
denoted as $\Delta_{11} + \Delta_{12}$, where
\begin{align*}
\tilde\psi_2(X) = \int_0^1 \left( \psi_2[ \bar \myeta^1_{\mbox{\tiny WL}}(X) + u\{ \hat \myeta^1_{\mbox{\tiny RWL}}(X)-\bar \myeta^1_{\mbox{\tiny WL}}(X)\}] -\psi_2\{\bar \myeta^1_{\mbox{\tiny WL}}(X)\} \right)\,\dif u .
\end{align*}
Denote by $\Omega_6$ the event that
\begin{align*}
\sup_{j=0,1,\ldots,p} \left| (\tilde E - E )\left[ \psi_2\{\bar \myeta^1_{\mbox{\tiny WL}}(X)\} f_j(X) \left\{ 1 - \frac{T}{\bar \pi^1_{\mbox{\tiny CAL}}(X)} \right\}  \right] \right| \le 2 C_1 \lambda_0.
\end{align*}
Then $P(\Omega_6 ) \ge 1- 2 \epsilon$ similarly as in Lemma~\ref{lem-prob0}(i). In the event $\Omega_6$, we have
\begin{align}
& | \Delta_{11} | \le \| \hat \alpha^1_{\mbox{\tiny RWL}} -\bar \alpha^1_{\mbox{\tiny WL}} \|_1
\sup_{j=0,1,\ldots,p}  \left| \tilde E \left[ \psi_2\{\bar \myeta^1_{\mbox{\tiny WL}}(X)\} f_j(X) \left\{ 1 - \frac{T}{\bar \pi^1_{\mbox{\tiny CAL}}(X)} \right\}  \right] \right| \nonumber \\
& \le \| \hat \alpha^1_{\mbox{\tiny RWL}} -\bar \alpha^1_{\mbox{\tiny WL}} \|_1 ( \Lambda_1 + 2 C_1 \lambda_0 ) . \label{prf-thm-mu2-eq4}
\end{align}
To bound $\Delta_{12}$, we have by Assumption~\ref{ass-RWL2}(iii),
%\begin{align*}
%\psi_2\{\bar \myeta^1_{\mbox{\tiny WL}}(X)\} \left( \me^{-C_4 | \hat \myeta^1_{\mbox{\tiny RWL}}(X) -\bar \myeta^1_{\mbox{\tiny WL}}(X)|} - 1 \right)\le
%\tilde \psi_2(X) \le \psi_2\{\bar \myeta^1_{\mbox{\tiny WL}}(X)\} \left( \me^{C_3 | \hat \myeta^1_{\mbox{\tiny RWL}}(X) -\bar \myeta^1_{\mbox{\tiny WL}}(X)| } - 1 \right)
%\end{align*}
%and hence
\begin{align}
& | \tilde \psi_2(X) | \le \psi_2\{\bar \myeta^1_{\mbox{\tiny WL}}(X)\} \left( \me^{C_3 | \hat \myeta^1_{\mbox{\tiny RWL}}(X) -\bar \myeta^1_{\mbox{\tiny WL}}(X)| } - 1 \right) \nonumber \\
& \le  \psi_2\{\bar \myeta^1_{\mbox{\tiny WL}}(X)\} C_3 | \hat \myeta^1_{\mbox{\tiny RWL}}(X) -\bar \myeta^1_{\mbox{\tiny WL}}(X)| \me^{C_3 | \hat \myeta^1_{\mbox{\tiny RWL}}(X) -\bar \myeta^1_{\mbox{\tiny WL}}(X)| } , \label{prf-thm-mu2-eq5}
\end{align}
where the second inequality follows because $(\me^c - 1)/c = \int_0^1 \me^{uc} \,\dif u \le \me^c$ for $c\ge 0$.
As a result, we find from Assumptions~\ref{ass-RCAL}(i) and \ref{ass-RWL2}(i),
\begin{align}
| \Delta_{12} | \le (1+\me^{-B_0}) C_1 C_3 \me^{C_4 \| \hat \alpha^1_{\mbox{\tiny RWL}} -\bar \alpha^1_{\mbox{\tiny WL}} \|_1 }
\tilde E \left[ \{\hat \myeta^1_{\mbox{\tiny RWL}}(X) - \bar \myeta^1_{\mbox{\tiny WL}}(X)\}^2  \right]. \label{prf-thm-mu2-eq6}
\end{align}
By condition (\ref{compat-add}),
$\tilde E [ \{\hat \myeta^1_{\mbox{\tiny RWL}}(X) - \bar \myeta^1_{\mbox{\tiny WL}}(X)\}^2  ]$ can be bounded as (\ref{prf-thm-var-eq7}) in the event $\Omega_1\cap \Omega_5$.
Then (\ref{thm-mu2-eq2}) follows by collecting inequalities (\ref{prf-thm-mu-eq2})--(\ref{prf-thm-mu-eq4})
and (\ref{prf-thm-mu2-eq1})--(\ref{prf-thm-mu2-eq6})
and applying (\ref{thm-RWL2-eq}), (\ref{lem-add-ineq-eq1}), and (\ref{gamma-L1-bound}) in the event $\Omega_0 \cap \Omega_1 \cap \Omega_2 \cap \Omega_3 \cap \Omega_4$.

Now suppose that (\ref{compat-add}) may not hold. Denote $h(X;\alpha^1) = \alpha^{\one\T} f(X)$. Then $\tilde E( \delta_1)$ can be decomposed as
\begin{align*}
& \tilde E(\delta_1 ) = (\tilde E - E )\left( \left[ \psi\{\hat\myeta^1_{\mbox{\tiny RWL}} (X)\}- \psi\{\bar \myeta^1_{\mbox{\tiny WL}} (X)\}\right]
\left\{ 1 - \frac{T}{\bar \pi^1_{\mbox{\tiny CAL}}(X)} \right\} \right) \\
& \quad + E \left( \left[ \psi\{\hat\myeta^1_{\mbox{\tiny RWL}} (X)\}- \psi\{\bar \myeta^1_{\mbox{\tiny WL}} (X)\}\right]
\left\{ 1 - \frac{T}{\bar \pi^1_{\mbox{\tiny CAL}}(X)} \right\} \right),
\end{align*}
denoted as $\Delta_{13} + \Delta_{14}$.
In the event  $\Omega_0 \cap \Omega_1 \cap \Omega_2 \cap \Omega_3 \cap \Omega_4$, we have
$\| \hat \alpha^1_{\mbox{\tiny RWL}} -\bar \alpha^1_{\mbox{\tiny WL}} \|_1  \le \eta_{11}$ from (\ref{thm-RWL2-eq}) and hence by
the mean value theorem,
\begin{align}
& | \Delta_{14} | \le \eta_{11}  \sup_{j=0,1,\ldots,p}
\left| E \left[ \psi_2\{\myeta (X; \tilde \alpha^1)\} f_j(X) \left\{ 1 - \frac{T}{\bar \pi^1_{\mbox{\tiny CAL}}(X)} \right\}  \right] \right|
\le \eta_{11} \Lambda_0 (\eta_{11}),\label{prf-thm-mu2-eq7}
\end{align}
where $\tilde\alpha^1$ lies between $\hat \alpha^1_{\mbox{\tiny RWL}}$ and $\bar \alpha^1_{\mbox{\tiny WL}} $.
Moreover, in the event $(\Omega_0 \cap \Omega_1 \cap \Omega_2 \cap \Omega_3 \cap \Omega_4) \cap \Omega_7$, applying Lemma~\ref{lem-prob-contract} below yields
\begin{align}
| \Delta_{13} | \le 2C_1 (1 + C_3 \me^{C_4 \eta_{11}}) \eta_{11} \lambda_0 . \label{prf-thm-mu2-eq8}
\end{align}
Then (\ref{thm-mu2-eq2}) follows by combining (\ref{prf-thm-mu2-eq7})--(\ref{prf-thm-mu2-eq8}) and other aforementioned inequalities.

\begin{lem} \label{lem-prob-contract}
For $r \ge 0$, denote by $\Omega_7$ the event that
\begin{align*}
\sup_{\|\alpha^1 - \bar\alpha^1_{\mbox{\tiny WL}} \|_1 \le r} \left| (\tilde E - E )\left( \left[\psi\{\myeta(X;\alpha^1)\}- \psi\{\bar \myeta^1_{\mbox{\tiny WL}} (X)\}\right]
\left\{ 1 - \frac{T}{\bar \pi^1_{\mbox{\tiny CAL}}(X)} \right\}  \right) \right| \le \sqrt{8} C_1 (1 + C_3 \me^{C_4 r}) r \lambda_0.
\end{align*}
Under Assumptions~\ref{ass-RCAL}(i)--(ii), \ref{ass-RWL2}(i) and \ref{ass-RWL2}(iii),
if $\lambda_0 \ge  \sqrt{2} (\me^{-B_0}+1) C_0 \sqrt{\log\{(1+p)/\epsilon\}/n}$, then $P(\Omega_6 ) \ge 1-2\epsilon$.
\end{lem}

\begin{prf}
Denote
$$
g (T,X;\alpha^1) = \left[ \psi\{\myeta(X;\alpha^1)\}- \psi\{\bar \myeta^1_{\mbox{\tiny WL}} (X)\} \right] \left\{ 1 - \frac{T}{\bar \pi^1_{\mbox{\tiny CAL}}(X)} \right\}   .
$$
For $\|\alpha^1 - \bar\alpha^1_{\mbox{\tiny WL}} \|_1 \le r$,  similar manipulation as in (\ref{prf-thm-mu2-eq3}) and (\ref{prf-thm-mu2-eq5}) using
Assumptions~\ref{ass-RCAL}(i), \ref{ass-RWL2}(i) and \ref{ass-RWL2}(iii) yields
\begin{align}
& \left| \psi\{\myeta(X;\alpha^1)\}- \psi\{\bar \myeta^1_{\mbox{\tiny WL}} (X)\} \right| \le
\psi_2\{\bar \myeta^1_{\mbox{\tiny WL}}(X)\} | \myeta(X; \alpha^1) - \bar \myeta^1_{\mbox{\tiny WL}}(X)| \nonumber \\
&\quad + \psi_2\{\bar \myeta^1_{\mbox{\tiny WL}}(X)\} C_3 | \myeta(X;\alpha^1) -\bar \myeta^1_{\mbox{\tiny WL}}(X)| \me^{C_3 | \myeta(X;\alpha^1) -\bar \myeta^1_{\mbox{\tiny WL}}(X)| } \nonumber \\
& \le C_1 (1 + C_3 \me^{C_4 r}) | \myeta(X;\alpha^1) -\bar \myeta^1_{\mbox{\tiny WL}}(X)| , \label{Lip-cond}
\end{align}
that is, $\psi()$ satisfies a Lipschitz condition.
By the symmetrization and contraction theorems (e.g., Buhlmann \& van de Geer 2011, Theorems 14.3 and 14.4), we have
\begin{align*}
& E \left[ \sup_{\|\alpha^1 - \bar\alpha^1_{\mbox{\tiny WL}}\|_1 \le r}  \left| (\tilde E - E)\{ g(T,X;\alpha^1)\} \right| \right]
\le 2 E  \sup_{\|\alpha^1 - \bar\alpha^1_{\mbox{\tiny WL}}\|_1 \le r}  \left| \frac{1}{n} \sum_{i=1}^n \sigma_i g(T_i,X_i;\alpha^1)  \right| \\
& \le 2  C_1 (1 + C_3 \me^{C_4 r})
\times E \sup_{\|\alpha^1 - \bar\alpha^1_{\mbox{\tiny WL}}\|_1 \le r}
\left| \frac{1}{n} \sum_{i=1}^n \sigma_i \{ \myeta(X_i;\alpha^1) -\bar \myeta^1_{\mbox{\tiny WL}}(X_i) \} \left\{ 1 - \frac{T_i}{\bar \pi^1_{\mbox{\tiny CAL}}(X_i)} \right\}  \right| \\
& \le  2  C_1 (1 + C_3 \me^{C_4 r}) r
\times E \sup_{j=0,1,\ldots,p}
\left| \frac{1}{n} \sum_{i=1}^n \sigma_i f_j(X_i) \left\{ 1 - \frac{T_i}{\bar \pi^1_{\mbox{\tiny CAL}}(X_i)} \right\}  \right| ,
\end{align*}
where $(\sigma_1,\ldots,\sigma_n)$ are independent Rademacher variables with $P(\sigma_i =1) =P(\sigma_i=-1)=1/2$ for each $i$.
By Hoeffding's moment inequality (Buhlmann \& van de Geer 2011, Lemma 14.14), we find from the preceding inequality
\begin{align*}
& E \left[ \sup_{\|\alpha^1 - \bar\alpha^1_{\mbox{\tiny WL}}\|_1 \le r}  \left| (\tilde E - E)\{ g(T,X;\alpha^1)\} \right| \right] \\
& \le  2  C_1 (1 + C_3 \me^{C_4 r}) r \times C_0 (\me^{B_0}+1) \sqrt{\frac{2 \log(2+2p)}{n}} ,
\end{align*}
by Assumption~\ref{ass-RCAL}(i)--(ii).
For $\|\alpha^1 - \bar\alpha^1_{\mbox{\tiny WL}} \|_1 \le r$, inequality (\ref{Lip-cond}) also shows that $|g (T_i,X_i;\alpha^1)| \le C_1 (1 + C_3 \me^{C_4 r}) C_0 (\me^{B_0}+1) r$.
By  Massart's inequality (Buhlmann \& van de Geer 2001, Theorem~14.2), we have with probability at least $1-2\epsilon$,
\begin{align*}
& \sup_{\|\alpha^1 - \bar\alpha^1_{\mbox{\tiny WL}}\|_1 \le r}  \left| (\tilde E - E)\{ g(T,X;\alpha^1)\} \right| \\
& \le C_0 (\me^{B_0}+1) C_1 (1 + C_3 \me^{C_4 r}) r \left\{ 2\sqrt{\frac{2 \log(2+2p)}{n}} + \sqrt{\frac{8 \log\{1/(2\epsilon)\}}{n}} \right\}  \\
& \le C_0 (\me^{B_0}+1) C_1 (1 + C_3 \me^{C_4 r}) r \sqrt{\frac{16\log\{(1+p)/\epsilon\}}{n}} ,
\end{align*}
where the second inequality uses $\sqrt{a} + \sqrt{b} \le \sqrt{2(a+b)}$.
\end{prf}

\subsection{Proof of Theorem~\ref{thm-var2}}

The proof is similar to that of Theorem~\ref{thm-var}.
First, (\ref{prf-thm-var-eq3}) for $\tilde E(\delta_2^2)$ remains valid.
Second, combining (\ref{prf-thm-var-eq4}) and (\ref{prf-thm-mu2-eq2}) yields
\begin{align*}
& \tilde E (\delta_3^2 )
\le \me^{ -B_0} \left(1+\me^{\|\hat \gamma^1_{\mbox{\tiny RCAL}} - \bar \gamma^1_{\mbox{\tiny CAL}}\|_1 C_0} \right)^2
 \me^{C_4 \| \hat \alpha^1_{\mbox{\tiny RWL}} - \bar \alpha^1_{\mbox{\tiny WL}}\|_1} C_1 D^\dag_{\mbox{\tiny WL}}(\hat \myeta^1_{\mbox{\tiny RWL}}, \bar \myeta^1_{\mbox{\tiny WL}}; \bar \alpha^1_{\mbox{\tiny CAL}} ) .
\end{align*}
Third, similarly as in (\ref{prf-thm-var-eq6}) and (\ref{prf-thm-mu2-eq2}), we have
\begin{align}
& \tilde E (\delta_1^2 ) = \tilde E \left[ \{\hat m^1_{\mbox{\tiny RWL}}(X) - \bar m^1_{\mbox{\tiny WL}}(X)\}^2 \left\{ 1 - \frac{T}{\bar \pi^1_{\mbox{\tiny CAL}}(X)} \right\}^2 \right] \nonumber \\
& \le (1 + \me^{-B_0})^2 \me^{2C_4 \| \hat \alpha^1_{\mbox{\tiny RWL}} - \bar \alpha^1_{\mbox{\tiny WL}}\|_1} C_1^2 \tilde E \left[
\{\hat \myeta^1_{\mbox{\tiny RWL}}(X) - \bar \myeta^1_{\mbox{\tiny WL}}(X)\}^2 \right] \label{prf-thm-var2-eq1} \\
& \le (1 + \me^{-B_0})^2 \me^{2C_4 \| \hat \alpha^1_{\mbox{\tiny RWL}} - \bar \alpha^1_{\mbox{\tiny WL}}\|_1} C_0^2 C_1^2 \| \hat \alpha^1_{\mbox{\tiny RWL}} - \bar \alpha^1_{\mbox{\tiny WL}} \|_1^2 .  \nonumber
\end{align}
Inequality (\ref{thm-var2-eq1}) follows by collecting the aforementioned inequalities and
applying (\ref{thm-RWL2-eq}), (\ref{thm-mu2-eq1}), (\ref{gamma-L1-bound}), and (\ref{lem-add-ineq2-eq3})
in the event $\Omega_0 \cap \Omega_1 \cap \Omega_2 \cap \Omega_3 \cap \Omega_4$.
If condition (\ref{compat-add}) holds, then  in the event $\Omega_1 \cap \Omega_5$, combining (\ref{prf-thm-var-eq7}) and (\ref{lem3-ineq4})
and using $(1-\me^{-c})/c \ge \me^{-c}$ for $c \ge 0$ yields
\begin{align}
& \tilde E \left[\{\hat \myeta^1_{\mbox{\tiny RWL}}(X) - \bar \myeta^1_{\mbox{\tiny WL}}(X)\}^2 \right]
\le  \me^{B_0} \lambda_0  \|  \hat\alpha^1_{\mbox{\tiny RWL}}- \bar\alpha^1_{\mbox{\tiny WL}}\|_1 ^2  \nonumber \\
& \quad + \tau_0^{-1}
\left\{  \lambda_0 \|  \hat\alpha^1_{\mbox{\tiny RWL}}- \bar\alpha^1_{\mbox{\tiny WL}}\|_1 ^2 +
\me^{C_4 \| \hat \alpha^1_{\mbox{\tiny RWL}} - \bar \alpha^1_{\mbox{\tiny WL}}\|_1} C_2^{-1}
D^\dag_{\mbox{\tiny WL}}( \hat \myeta^1_{\mbox{\tiny RWL}}, \bar \myeta^1_{\mbox{\tiny WL}}; \bar \gamma^1_{\mbox{\tiny CAL}} ) \right\} . \label{prf-thm-var2-eq2}
\end{align}
Inequality (\ref{thm-var2-eq2}) follows by combining (\ref{prf-thm-var2-eq1})--(\ref{prf-thm-var2-eq2}) and other aforementioned inequalities except that
inequality (\ref{thm-mu2-eq1}) is replaced by (\ref{thm-mu2-eq2}).

\newpage

\subsection{Technical tools} \label{sec:tech-tool}

For completeness, we state the following concentration inequalities, which can be obtained from Buhlmann \& van de Geer (2011), Lemmas 14.11,  14.16, and 14.9.

\begin{lem} \label{lem-max-bounded}
Let $(Y_1,\ldots,Y_n)$ be independent variables such that $E(Y_i)=0$ for $i=1,\ldots,n$ and
$\max_{i=1,\ldots,n } |Y_i| \le c_0$
for some constant $c_0$. Then for any $t >0$,
\begin{align*}
P \left(  \left| \frac{1}{n} \sum_{i=1}^n Y_i \right| > t \right) \le 2 \exp \left(-\frac{n t^2}{2 c_0^2} \right).
\end{align*}
\end{lem}

\begin{lem} \label{lem-max-subG}
Let $(Y_1,\ldots,Y_n)$ be independent variables such that $E(Y_i)=0$ for $i=1,\ldots,n$ and $(Y_1,\ldots,Y_n)$ are uniformly sub-gaussian:
$\max_{i=1,\ldots,n } c_1^2  E \{\exp(Y_i^2/c_1^2) - 1 \} \le c_2^2$
for some constants $(c_1,c_2)$. Then for any $t >0$,
\begin{align*}
P \left(  \left| \frac{1}{n} \sum_{i=1}^n Y_i \right| > t \right) \le 2 \exp \left\{-\frac{n t^2}{8(c_1^2 + c_2^2)} \right\}.
\end{align*}
\end{lem}

\begin{lem} \label{lem-max-subE}
Let $(Y_1,\ldots,Y_n)$ be independent variables such that $E(Y_i)=0$ for $i=1,\ldots,n$ and
% $\max_{i=1,\ldots,n } 2 c_3^2  E \{\exp(|Y_i|/c_3) - 1 - |Y_i|/c_3 \} \le c_4^2$
$$
\frac{1}{n} \sum_{i=1}^n E ( |Y_i|^k ) \le \frac{k!}{2} c_3^{k-2} c_4^2, \quad k=2,3,\ldots,
$$
for some constants $(c_3,c_4)$. Then for any $t >0$,
\begin{align*}
P \left(  \left| \frac{1}{n} \sum_{i=1}^n Y_i \right| > c_3 t + c_4 \sqrt{2 t} \right) \le 2 \exp (-nt).
\end{align*}
\end{lem}

The following results about sub-gaussian variables can be deduced from Buhlmann \& van de Geer (2011), Lemmas 14.3 and 14.5.

\begin{lem} \label{lem-subG-moment}
Suppose that $Y$ is sub-gaussian: $c_1^2  E \{\exp(X^2/c_1^2) - 1 \} \le c_2^2$ for some constants $(c_1,c_2)$.
Then
$$
E ( |Y|^k ) \le \Gamma \left(\frac{k}{2}+1 \right) (c_1^2 + c_2^2) c_1^{k-2} , \quad k=2,3,\ldots.
$$
\end{lem}

\begin{lem} \label{lem-subG-sq}
Suppose that $X$ is bounded: $|X| \le c_0$ for a constant $c_0$, and $Y$ is sub-gaussian: $c_1^2  E \{\exp(X^2/c_1^2) - 1 \} \le c_2^2$ for some constants $(c_1,c_2)$.
Then $Z =X Y^2 $ satisfies
$$
E \left\{ |Z - E(Z)|^k \right\} \le \frac{k!}{2} c_3^{k-2} c_4^2, \quad k=2,3,\ldots,
$$
for $c_3 = 2 c_0 c_1^2$ and $c_4 = 2 c_0 c_1 c_2$.
\end{lem}

\end{document}